\documentstyle[bezier]{article}
\date{March 29, 1995}
\newtheorem{theo}{Theorem}
\newtheorem{defn}[theo]{Definition}

\newtheorem{lemma}{Lemma}
\newtheorem{prop}[theo]{Proposition}
\newtheorem{coro}[theo]{Corollary}

\def\sg{{\scriptstyle {\cal G}}}
\def\End{{\mbox{\it End\/}}}
\def\Aut{{\mbox{\it Aut\/}}}
\def\Hom{{\mbox{\it Hom\/}}}
\def\Hol{{\mbox{\it Hol\/}}}
\def\Inv{{\mbox{\it Inv\/}}}
\def\tr{{\mbox{\it tr\/}}}
\def\be{\begin{equation}}
\def\ba{\begin{eqnarray}}
\def\ee{\end{equation}}
\def\ea{\end{eqnarray}}
\def\o{\otimes }
\def\bo{\mbox{\,\raisebox{-0.65mm}{$\Box$} \hspace{-4.7mm}
${\scriptstyle\times}$ \/}}
\def\D{\Delta }

\def\k{{\cal k}}

\def\MC{{\cal M}}
\def\T{{\cal T}}
\def\lM{{\cal M}}
\def\A{{\cal A}}

\def\G{{\cal G}}
\def\V{{\cal V}}

\def\N{{\cal N}}
\def\S{{\cal S}}
\def\cO{{\cal O}}
\def\De{{\cal D}}

\def\C{{\cal C}}

\def\L{{\cal L}}
\def\ti{\times }
\def\vp{\varphi }
\def\s{\sigma }

\def\om{\omega}
\def\a{\alpha }
\def\b{\beta }

\def\d{\delta }
\def\e{\epsilon }

\def\k{\kappa }

\def\ha{\hat a}
\def\hb{\hat b}

\def\hd{\hat d}
\def\he{\hat e}

\def\heta{\hat \eta}

\def\hvrho{\hat \varrho}

\def\sd{\ti_S}
\def\hh{\hat h}

\def\vth{\vartheta}
\def\th{\theta}
\def\vac{|0 \rangle}

\def\t{\tau}
\def\nn{\nonumber}

\newcommand{\M}[2]{\stackrel{\scriptscriptstyle #1}{M}
            \hspace*{-4mm} \phantom{M}^{#2}}
\newcommand{\Ae}[2]{\stackrel{\scriptscriptstyle #1}{A}
            \hspace*{-4mm} \phantom{A}^{#2}}
\newcommand{\Be}[2]{\stackrel{\scriptscriptstyle #1}{B}
            \hspace*{-4mm} \phantom{B}^{#2}}

\begin{document}
\begin{titlepage}
\title{Representation Theory of Chern-Simons Observables}
\author{{\sc Anton Yu. Alekseev}
\thanks{On leave of absence from Steklov Mathematical Institute,
Fontanka 27, St.Petersburg, Russia; e-mail:
alekseev@itp.phys.ethz.ch}
\\Institut f\"ur Theoretische Physik,
\\ ETH -- H\"onggerberg,
CH-8093 Z\"urich, Switzerland \\
and \\
Department of Mathematics, Harvard University
\thanks{Tha visit of A.A. to Harvard University has been supported
by Alfred P. Sloan Foundation Grant No 90-10-13}\\
Cambridge, MA 02138, U.S.A. \\[4mm]
{\sc Volker Schomerus \thanks{Supported in
part by DOE Grant No DE-FG02-88ER25065; e-mail:
mathphys@huhepl. harvard.edu}} \\
Harvard University, Department of Physics \\ Cambridge, MA 02138,
U.S.A.}
\maketitle \thispagestyle{empty}
%--------------------------- Abstract---------------------------

\begin{abstract}
In \cite{AGS1}, \cite{AGS2} we suggested a new quantum algebra,
the moduli algebra, which
is conjectured to be a quantum algebra of observables of
the Hamiltonian Chern-Simons
theory.
This algebra  provides the quantization
of the algebra of functions on the moduli space of flat connections
on a 2-dimensional surface. In this paper we classify unitary representations
of this new algebra and identify the corresponding
representation spaces with
the spaces of conformal blocks of the WZW model. The mapping class
group of the surface is proved to act on the moduli algebra by
inner automorphisms. The generators of these automorphisms are
unitary elements of the moduli algebra. They are constructed
explicitly and proved to satisfy the relations of the (unique)
central extension of the mapping class group.
\end{abstract}
\vspace*{-17cm}
\hspace*{9cm}
{\large \tt  HUTMP 95-B342} \\
\hspace*{9cm}
{\large \tt ETH-TH/95-12}         \\
\hspace*{9cm}
{\large \tt ESI 214 (1995)}   \\
\hspace*{9cm}
{\large \tt hep-th 9503016}
\end{titlepage}

\tableofcontents
\def\tt{\tilde{\tau}}

\section{Introduction}

This paper is devoted to operator approach to quantization
of the Chern-Simons theory, sometimes called Combinatorial
Quantization. From a more mathematical perspective this is
the same as quantization of the algebra of functions on
the moduli space of flat connections on a Riemann surface.
In \cite{AGS1}, \cite{AGS2} we suggested the definition
of the moduli algebra which solves this problem and investigated
some elementary properties of this object.

In order to make
the paper more accessible to the reader, we collect in
Part 1  some motivations which
lead to the idea of  Combinatorial Quantization (following
\cite{FoRo}) and some results and constructions from
\cite{AGS1}, \cite{AGS2}. This provides a background for
the further consideration. We formulate
the main results in Section 4.

In Sections 5-8 we deal with the representation theory
of the moduli algebra. All unitary representations
are classified for moduli algebras corresponding to surfaces
of any genus with arbitrary number of marked points.
The results are in agreement with other methods of
quantization of the moduli space, like geometric quantization
\cite{ADW} and Conformal Field Theory \cite{Wit1}.

Section 9 includes  results on the action of the mapping
class group of the surface on the quantized moduli space.
As the mapping class group is acting on the moduli space
of flat connections before quantization, it is natural
to expect that this action survives in the quantum case.
This is confirmed by the experience of Conformal Field Theory
and we will indeed establish this result in the framework
of Combinatorial Quantization. Unitary generators of
the mapping class group action are described as elements
of the moduli algebra by explicit formulae. They are proved
to satisfy the relations for the (unique) central extension
of the mapping class group.

In this work we extensively use certain results and ideas
of the recent papers related to the program of Combinatorial
Quantization \cite{FoRo}, \cite{Bou}, \cite{AGS1}, \cite{AGS2},
\cite{BuRo}, \cite{AMR}, \cite{Tur}, \cite{AA}.

\part{Background Review}

\section{Chern-Simons Theory and the Moduli Space
of Flat Connections}
\setcounter{equation}{0}

The Chern-Simons theory is constructed by
the following data. We pick up some semi-simple Lie algebra
$\sg$, a coupling constant $k$ and a 3-manifold $M$.
The Chern-Simons action is a functional of a $\sg$-valued
gauge field $A$ on $M$ and has the form
\begin{equation} \label{CS}
CS(A)=\frac{k}{4\pi}Tr\int_M (A\wedge dA +\frac{2}{3}A^3).
\end{equation}
We do not describe truly exceptional properties of this action
and of the corresponding theory and refer the reader to the
original papers and numerous reviews. Let us only
mention that for $\sg$ being a compact Lie algebra,
$k$ is required to be integer in order to ensure the global
gauge invariance.

When the manifold $M$ has a structure of a direct product
of a 2-dimensional surface and a segment of a real line, the
Chern-Simons theory admits a Hamiltonian interpretation. As
in any topological theory, the Hamiltonian is equal to zero.
The action (\ref{CS}) induces a symplectic structure on the
space of connections \cite{AtBo}:
\begin{equation} \label{AtBo}
\Omega=\frac{k}{4\pi}Tr\int_{\Sigma} \delta_1A\wedge \delta_2A.
\end{equation}
This symplectic form is invariant with respect to
gauge transformations
\begin{equation} \label{gauge}
A^g=gAg^{-1}+dgg^{-1}.
\end{equation}
An easy check shows that  the moment map for the gauge
group action is proportional to the curvature
\begin{equation}
F=dA+A^2.
\end{equation}
The condition
\be \label{F=0}
F=0
\ee
emerges as a constraint from the Chern Simons action.
{}From this analysis we see that the phase space of the
Chern Simons model is a quotient of the space of flat
connections (\ref{F=0}) over the gauge group action
(\ref{gauge}). In this paper we often refer to this
space as moduli space which should always be understood
as a moduli space of flat connections on a Riemann surface.

\subsection{Combinatorial description}

Let us introduce a more efficient finite-dimensional
description of the moduli space which will be of
use throughout the paper. At the same moment we
introduce moduli spaces on surfaces with marked
points.  Let $G$ be a semi-simple
connected simply-connected Lie group corresponding
to the Lie algebra $\sg$ and $\Sigma_{g,m}$
be a surface of genus $g$ with $m$ marked points.
Assign a conjugacy class ${\cal C}_\nu\in G$ to the
$\nu$'s marked point. Denote $\pi=\pi_{g,m}$ the fundamental
group of the surface $\Sigma_{g,m}$. The group $\pi$ may be
generated by $2g+m$ invertible generators $a_i, b_i, i=1,\dots ,g$ and
$l_\nu, \nu=1,\dots , m$ subject to the relation
\be        \label{abm}
[b_g, a_g^{-1}]\dots [b_1, a_1^{-1}]l_m\dots l_1=id\ \ . \ \ (in\ \pi)
\ee
Here $[x, y]$ stays for the group commutator $xyx^{-1}y^{-1}$.
\begin{defn}
The moduli space of flat connections
${\bf m}_{g,m}^{\{ {\cal C}_\nu\} }$
on the Riemann surface of genus $g$ with $m$ marked points
is defined as
\be
{\bf m}_{g,m}^{\{ {\cal C}_\nu\} }=
\{ \rho \in \Hom(\pi, G), \rho(l_\nu)\in {\cal C}_\nu \} / G\ \ .
\ee
Here the group $G$ acts on the space of representations of $\pi$
by conjugations
\be
\rho^g(x)=g^{-1}\rho(x)g \ \ .
\ee
\end{defn}

In order to make
contact with the definition which involves flat connections we
represent generators of $\pi$ as circles on $\Sigma_{g,m}$
intersecting at the base point on the surface. Then a
flat connection $A$ induces a representation of the fundamental group
via
\begin{eqnarray}
M_\nu & = & \rho(l_\nu)\ =\ \Hol(A, l_\nu)\ \ ; \nonumber \\
A_i & = & \rho(a_i)\ =\ \Hol(A, a_i)\ \ ; \label{MAB} \\
B_i & = & \rho(b_i)\ =\ \Hol(A, b_i)\ \ . \nonumber
\end{eqnarray}
Here $\Hol(A, x)$ is a holonomy of the connection $A$ along
the cycle $x$. Connection $A$ being flat, $\rho$ satisfies
the defining relation
\be \label{BAM}
   M = [B_g, A_g^{-1}]\dots [B_1, A_1^{-1}] M_m\dots M_1=id
   \ \ .\ \  (in\  G)
\ee

The case of a surface with marked points requires some further
comment. The marked points are deleted from the surface and the
condition of flatness does not hold there. In general we permit
the curvature to develop $\delta$-function  singularities at
the marked points:
\be \label{Fsum}
F=\sum_\nu T_\nu \delta^{(2)}(z-z_\nu)\ \ .
\ee
Here $T_\nu$'s are elements of the algebra $\sg$.
However, if one permits arbitrary $T_\nu$'s,
the space of all connections fails to be  symplectic.
It regains this property if we fix conjugacy
classes of $T_\nu$ in $\sg$. This condition may be naturally
derived from the Chern-Simons action functional in the same fashion
as one obtains the flatness condition for a surface without marked
points \cite{EM}.

Assume that we have fixed some conjugacy classes for $T_\nu$'s
in the modified flatness condition  (\ref{Fsum}).
Then via equations (\ref{MAB}) one can define a subset of
representations of the fundamental group of a surface with
marked points. It is easy to check that the condition
(\ref{Fsum}) implies that $\rho(l_\nu)$ belong to certain
conjugacy classes in the group $G$ as it is stated in the definition of
the moduli space. In fact, one can get the conjugacy class
of $\rho(l_\nu)$ by applying the exponential map $exp: \sg
\rightarrow G$ to the conjugacy class of $T_\nu$.
In this way we obtain a family of symplectic spaces labeled
by two nonnegative integers (genus and the number of marked
points) and by conjugacy classes attached to the marked points.

Notations (\ref{MAB}) prove to be quite useful when constructing
functions on the moduli space. Take any function $f$ on the
the direct product of $2g+m$ copies of the group $G$
which is  invariant
with respect to  simultaneous conjugations of the arguments
\be \label{inv}
  f(g^{-1}g_1g,\dots , g^{-1}g_{(2g+m)}g)=
    f(g_1,\dots , g_{(2g+m)})\ \ .
\ee
Define a function $f(\rho)$ on the space of representations
via
\be \label{fro}
f(\rho)=f(M_1,\dots , M_m, A_1, B_1, \dots , A_g, B_g).
\ee
It is easy to see that this function descends to the quotient
space. In fact, any analytic function on the moduli space may
be obtained in this way.

\subsection{Poisson structure of the moduli space}

As the quantization program which we are going to develop
is close to deformation quantization, we describe the
Poisson bracket on the moduli space defined by Atiyah-Bott
symplectic structure. This Poisson bracket was
described in \cite{Gol}. However, we use another description
of the same object which was especially designed for the
needs of deformation quantization \cite{FoRo}.

Let us introduce some useful notations. We define left
and right invariant differential operators on a Lie group
with values in the dual space to its Lie algebra:
\begin{eqnarray}
<\nabla_L f (g), X> & = & \frac{d}{dt} f(e^{-tX}g)|_{t=0}\ \ ;
        \nonumber  \\[1mm]
<\nabla_R f (g), X> & = & \frac{d}{dt} f(ge^{tX})|_{t=0}\ \ .
\end{eqnarray}

As we work with semi-simple algebras, one can introduce a
nondegenerate Killing
form $K\in \sg \otimes \sg$. Let us consider one more
element   $r\in \sg \otimes \sg$. It is called
classical $r$-matrix if it satisfies the classical Yang-Baxter
equation in $U(\sg)^{\otimes 3}$:
\begin{equation}
[r_{12}, r_{13}]+[r_{12},r_{23}]+[r_{13},r_{23}]=0\ \ .
\end{equation}
Here we denote $r_{12}=r\otimes 1$, $r_{13}$ and $r_{23}$
are constructed in a similar way. For any solution $r$ of the classical
Yang-Baxter equation one can construct another solution $r'$
which is obtained by permutation of two copies of $\sg$.
Let us prepare their symmetric and antisymmetric combinations
\be
r_s=\frac{1}{2}(r+r')\ \ ,\  \ r_a=\frac{1}{2}(r-r')\ \ .
\ee
In general, neither $r_a$ nor $r_s$  satisfy the classical
Yang-Baxter equation.

In order to describe the Poisson bracket on the moduli space
we use its description via representations of the fundamental
group from the previous subsection. As the space of representations
is embedded into $G^{2g+m}$, it is convenient to order the
$2(2g+m)$ covariant differential operators in the following way
\begin{eqnarray}
\nabla_{2\nu-1} = \nabla_R^{M_\nu} & , &
\nabla_{2\nu}=\nabla_L^{M_\nu}\ \ \mbox{ for } \ \nu =1,\dots, m\ ;
\nonumber  \\[1mm]
\nabla_{m+4i-3}=\nabla_R^{A_i} & , &
\nabla_{m+4i-1}=\nabla_L^{A_i} \ \ \mbox{ for } \ \  i=1,\dots , g\ ;
\\[1mm]
\nabla_{m+4i-2}=\nabla_R^{B_i} & , &
\nabla_{m+4i-1}=\nabla_L^{B_i}\ \ \mbox{ for }\  i=1,\dots , g\ .
\nonumber
\end{eqnarray}
One can rewrite the condition for a function $f$ on $G^{2g+m}$
to be invariant in the sense of (\ref{inv}) as
\be
\sum_{s=1}^{2(2g+m)} \nabla_s f=0\ \ .
\ee
With these notations we present a description of the Poisson
structure on the moduli space.

\begin{theo} {\em (Fock-Rosly \cite{FoRo})}  \label{FoRoth}
Let $r\in \sg \otimes \sg$ be a solution of the classical
Yang-Baxter equation. Assume that its symmetric part coincides
with the Killing bilinear $r_s=K$.
Introduce a Poisson bracket on $G^{2g+m}$ by the following
formula
\be \label{FoRo}
\{ f , h\}=\frac{1}{2} \sum_i <r, \nabla_i f\wedge \nabla_i h>
+\sum_{i<j} <r, \nabla_i f\wedge \nabla_j h>.
\ee
This Poisson bracket restricts to the space of functions
which are invariant with respect to simultaneous conjugations. Being
mapped to the moduli space by means of equation (\ref{fro}), the bracket
(\ref{FoRo}) coincides with the canonical Poisson bracket defined
by the Atiyah-Bott symplectic structure.
\end{theo}

The Poisson bracket (\ref{FoRo}) may be easily evaluated
for simplest functions on $G^{2g+m}$. Examples of such
functions are given by matrix elements of $M_\nu, A_i$ and $B_i$
evaluated in some irreducible representations. We denote
a matrix which represents some $X\in G$ in the representation $\tau^I$
by $X^{I}$. Applying a pair of representations $\tau^I$ and $\tau^J$
to the classical $r$-matrix we produce $r^{IJ}=(\tau^I\otimes \tau^J)(r)$.
It is convenient to introduce tensor notations $X^1=X\otimes 1,
X^2=1\otimes X$. Let us give some examples of Poisson
brackets for particular matrix elements:
\begin{eqnarray} \label{examp}
\{ \M{1}{I}_\nu , \M{2}{J}_\mu \} & = &
r^{IJ}\M{1}{I}_\nu \M{2}{J}_\mu - \M{1}{I}_\nu r^{IJ} \M{2}{J}_\mu
\nonumber \\
&  & \hspace*{1.5cm}
- \M{2}{J}_\mu r^{IJ} \M{1}{I}_\nu + \M{1}{I}_\nu \M{2}{J}_\mu r^{IJ}
\ \ for \ \nu<\mu\ \  ;
\nonumber \\[3mm]
\{ \M{1}{I}_\nu , \M{2}{J}_\nu \} & = &
r^{IJ} \M{1}{I}_\nu \M{2}{J}_\nu - \M{1}{I}_\nu r^{IJ} \M{2}{J}_\nu
\nonumber \\
& & \hspace*{1.5cm}
- \M{2}{J}_\nu (r')^{IJ} \M{1}{I}_\nu + \M{1}{I}_\nu
   \M{2}{J}_\nu (r')^{IJ} \ \ ;
\\[3mm]
\{ \Ae{1}{I}_i , \Be{2}{J}_i \} & = &
r^{IJ}\Ae{1}{I}_i \Be{2}{J}_i - \Ae{1}{I}_i r^{IJ} \Be{2}{J}_i
- \Be{2}{J}_i (r')^{IJ} \Ae{1}{I}_i + \Ae{1}{I}_i \Be{2}{J}_i r^{IJ}\  .
\nonumber
\end{eqnarray}
As one can see, a Poisson bracket of two matrix elements
is quadratic in matrix elements of the same representations.
Classical $r$-matrices $r^{IJ}$ and $(r')^{IJ}$ play the
role of structure constants in this Poisson bracket algebra.
This point is the most important observation of \cite{FoRo}
as it makes it possible to proceed with quantization of
the algebra of functions on the moduli space.

Matrix elements of $M_\nu, A_i, B_i$ which we considered so far
do not define functions on the moduli space as they are not
invariant with respect to conjugations. Let us remark that
the Poisson bracket (\ref{FoRo}) is not conjugation invariant.
Nevertheless it may be consistently restricted to the set
of conjugation invariant functions which means that a Poisson
bracket of two such functions is again conjugation invariant.
Conjugation invariant functions are produced as linear
combinations of elements
$$ tr^J \left( C_1 [I_1, \dots, I_{2g+m}|J] M^{I_1}_1
  \dots M^{I_m}_m A_1^{I_{m+1}}\dots B_g^{I_{2g+m}}
   C_2 [I_1, \dots, I_{2g+m}|J]^*\right) $$
for arbitrary sets of labels $\{I_\nu\} , J$ and two intertwiners
$C_1,C_2 : V^{I_1} \o \dots \o V^{I_{2g+m}} \mapsto V^J$.

As we have seen in the previous subsection, any invariant
function on $G^{2g+m}$ defines a function on the moduli space.
In fact, we restrict an invariant function to the subset
defined by the conditions (cp. formula \ref{BAM} for notations)
\begin{equation} \label{MMC}
M_\nu\in {\cal C}_\nu\ \ , \ \ M=1\ \ .
\end{equation}
By the theorem of Fock and Rosly this restriction is consistent
with the structure of the bracket. This means that constraints
(\ref{MMC}) are preserved by Hamiltonian flows produced by
invariant functions. As it is stated in Theorem \ref{FoRoth}, the
resulting bracket on the moduli space coincides with the canonical
one.

\section{Quantization of the Moduli Space}

As we mentioned above, the structure constants of Poisson algebra
(\ref{FoRo}) are defined by a couple of classical $r$-matrices
$r$ and $r'$. The key point in the quantization procedure is the
fact that we know a family of solutions of the quantum Yang-Baxter equation
\begin{equation} \label{YB}
R_{12}R_{13}R_{23}=R_{23}R_{13}R_{12}
\end{equation}
labeled by a  parameter $h$ so that
\begin{equation} \label{cor}
R^{IJ}|_{h\rightarrow 0}= I+hr+{\cal O}(h^2)\ \ .
\end{equation}
We assume that $R^{IJ}$ exist for any $I$ and $J$ and that they
are represented by matrices in the same spaces as $r^{IJ}$.
Our intention is to replace the Poisson algebra with structure
constants $r, r'$ by an associative algebra with structure
constants $R, R'$. A natural
framework for this construction is provided by the theory of
quasi-triangular Hopf algebras.

\subsection{Quasi-triangular Hopf algebras}
\setcounter{equation}{0}

The objects we wish to investigate later are associated
with a quantum symmetry algebra $\G$. More precisely,
$\G$ is a ribbon Hopf-*-algebra, i.e. a $*$-algebra
$\G$ with co-unit $\e: \G \mapsto {\bf C} $, co-product
$\D: \G \mapsto \G \o \G$, antipode $\S: \G \mapsto \G$,
$R$-matrix $R \in \G \o \G $ and the ribbon element $v
\in \G $. We do not want to spell out all the standard
axioms these objects have to satisfy in order to give
a ribbon Hopf algebra (for a complete definition
see e.g. \cite{ReTu}). Let us stress, however,  that we
deal with structures for which the co-product $\D$ is
consistent with the action
$$ (\xi \o \eta)^* = \eta^* \o \xi^* \ \ \mbox{ for all}
   \ \ \eta,\xi \in \G  $$
of the *-operation $\ast$ on elements in the tensor product
$\G \o \G$. This case is of particular interest, since
it appears for the quantized universal enveloping  algebras
$U_q(\sg)$ when the complex parameter $q$ has values on the
unit circle \cite{MSIII}.

Given the standard expansion of $R \in \G \o \G$,
$R = \sum r^1_{\s} \o r^2_{\s}$, one constructs the elements
\be
  u = \sum \S(r^2_{\s})   r^1_{\s}\ \ . \label{u}
\ee
Among the properties of $u$ (cp. e.g. \cite{ReTu}) one finds
that the product $u \S(u)$ is in the center of $\G$. The
ribbon element $v$ is a central square root of
$u \S (u)$ which obeys the following relations
\ba
v^2 = u \S(u) \ \ & , & \ \ \S (v) = v \ \ , \ \ \
\e(v) =1 \ \ ,\label{v}      \\[2mm]
v^* = v^{-1}\ \ & , & \ \
\D(v)  =   (R'R)^{-1} (v \o v) \label{eigRR}\ \ .
\ea
The elements $u$ and $v$ can be combined to furnish a
grouplike unitary element $g = u^{-1} v \in \G$.
Examples of ribbon-Hopf-*-algebras are given by the
enveloping algebras of all simple Lie algebras
\cite{ReTu}.

In the following we make several additional assumptions
about the ribbon Hopf algebra $\G$. To begin with, we
will assume that $\G$ is {\em semisimple}. More restrictions
will be imposed in subsection 3.4 .

For every equivalence class $[J]$ of irreducible *-representations
of $\G$ we pick a particular representative $\t^J$ with carrier
space $V^J$. The tensor product $\t^I \bo \t^J$ of two representations
$\t^I,\t^J$ of the semisimple algebra $\G$ can be decomposed
into irreducibles $\t^K$. This decomposition determines
the Clebsch-Gordon maps $C^a[IJ|K]: V^I \o V^J \mapsto V^K$,
\be          \label{CGint}
C^a [IJ|K] (\t^I \bo \t^J) (\xi) = \t^K(\xi) C^a[IJ|K]\ \ .
\ee
The same representations $\t^K$ in general appears with some
multiplicity $N^{IJ}_K$. The superscript $a= 1, \dots, N^{IJ}_K$
keeps track of these subrepresentations. It is common to call
the numbers $N^{IJ}_K$ {\em fusion rules}. Normalization of these
Clebsch Gordon maps is connected with an extra assumption.
Notice that the ribbon element $v$ is central so that the
evaluation with irreducible representations $\t^I$
gives complex numbers $v_I = \t^I(v)$. We suppose that
there exists a set of square roots $\k_I, \k_I^2 = v_I, $
such that
\be
C^a[IJ|K]  (R')^{IJ} C^b[IJ|L]^*  =
\delta_{a,b} \delta_{K,L} \frac{\k_I \k_J }{\k_K} \ \ .
\label{pos}
\ee
Here $R' = \sum r^2_{\s} \o r^1_{\s}$ and $(R')^{IJ} =
(\t^I \o \t^J) (R')$. The adjoint of the Clebsch Gordon map
is meant with respect to the standard scalar product on
$V^I \o V^J$ induced by the scalar products on $V^I,V^J$.
Let us analyze this relation in more detail. As a consequence of
intertwining properties of the Clebsch Gordon maps and the $R$-element,
$\t^K(\xi)$ commutes with the left hand side of the equation.
So by Schurs' lemma, it is equal to the identity $e^K$ times some
complex factor $\omega_{ab} (IJ|K)$. After appropriate normalization,
$\omega_{ab}(IJ|K) = \d_{a,b} \omega(IJ|K)$ with a complex phase
$\omega(IJ|K)$. Next we exploit the $*$-operation
and relation (\ref{eigRR}) to find $\omega_{ab}(IJ|K)^2 =
v_I v_J/v_K$. This means that (\ref{pos}) can be ensured up to
a possible sign $\pm$. Here we assume that this sign is always
$+$. This assumption was crucial for the positivity in \cite{AGS1}.
It is met by the quantized  universal enveloping algebras of
all simple Lie algebras because they are obtained as
deformations of  Hopf-algebras which clearly
satisfy (\ref{pos}). As a consequence of the normalization
equation (\ref{pos}) and the equation $\D(e) = e \o e$ we
obtain the following completeness
\be  \label{complete}
    \sum_{K,a}  \frac{\k_{K}}{\k_I \k_J}
   (R')^{I J} C^a[I J |K ]^* C^a[I J |
     K]  = e^I \o e^J    \ \ .
\ee

We wish to combine the phases $\k_I$ into one element $\k$ in
the center of $\G$, i.e. by definition, $\k$ will denote
a central element
\be \k \in \G \ \ \ \mbox{ with } \ \ \t^J(\k) = \k_J
\label{k} \ \ . \ee
Such an element does exist and is unique. It has the property
$\k^* = \k^{-1}$.

The antipode $\S$ of $\G$ furnishes   a conjugation
in the set of equivalence classes of irreducible representations.
We use $[\bar J]$ to denote the class conjugate to $[J]$.
Some important properties of the fusion rules $N^{IJ}_K$
can be formulated with the help of this conjugation.
Among them are the relations
\be N^{K \bar K}_0 = 1 \ \ , \ \ N^{IJ}_K = N^{JI}_K =
    N^{J \bar K}_{\bar I}\ \ .  \ee
The numbers $v_I$ are symmetric under conjugation,
i.e. $v_K = v_{\bar K}$.

The {\em $q$-trace} of an element $X \in \End(V^K)$ is
defined by
\be
   tr^K_q (X) =  tr^K(X \t^K(g)) \label{sqtrace}
\ee
where $g= u^{-1} v$ is the grouplike element introduced
above. Let us also mention that the $q$-trace of the
identity map $e^K \in \End(V^K)$ computes the ``quantum
dimension'' $d_K$ of the representation $\t^K$ \cite{ReTu},
i.e.
\be
 d_K \equiv tr^K_q (e^K) \ \ . \label{qdim}
\ee
The numbers $d_K$ satisfy the equalities $d_I d_J =
\sum N^{IJ}_K d_K$ and $d_K = d_{\bar K}$. In general,
the quantum dimensions $d_I$ differ from the dimensions
$\d_I$ of the representation spaces $V^I$.

\subsection{The graph algebra $\L_{g,m}$}

Equipped with the technique and notations of the preceeding
subsection we introduce a quantized version of the Poisson
algebra (\ref{FoRo}) on $G^{2g+m}$. This is the first step
towards quantization of the moduli space.

\begin{defn} {\em (Graph algebra $\L_{g,m}$)}   \label{Agn}
The {\em graph-algebra $\L_{g,m}$} is an associative algebra
generated by matrix elements of $M_\nu^I, A^I_i, B^I_i
\in End(V^I) \o \L_{g,m}, \nu = 1, \dots, m, i = 1, \dots,g$.
The superscript $I$ runs through
the set of irreducible representations of a quantum
symmetry algebra  $\G$ with $R$-element $R$. Elements
in $\L_{g,m}$ are subject to the following relations
\ba
  \M{1}{I}_\nu R^{IJ}\M{2}{J}_\nu
      &=& \sum C^a[IJ|K]^* M_\nu^K C^a[IJ|K]\ \  ,
      \label{defAgn} \\[1mm]
  \Ae{1}{I}_i R^{IJ}\Ae{2}{J}_i
      &=& \sum C^a[IJ|K]^* A_i^K C^a[IJ|K]\ \ , \\[1mm]
  \Be{1}{I}_i R^{IJ}\Be{2}{J}_i
      &=& \sum C^a[IJ|K]^* B_i^K C^a[IJ|K]\ \ ,
      \label{defAgn2}\\[1mm]
  (R^{-1})^{IJ} \Ae{1}{I}_i R^{IJ} \Be{2}{J}_i
      &=& \Be{2}{J}_i (R')^{IJ} \Ae{1}{I}_i R^{IJ} \ \ ,
                                    \label{defAgn3-} \\[1mm]
  (R^{-1})^{IJ} \M{1}{I}_\nu R^{IJ} \M{2}{J}_\mu
      &=& \M{2}{J}_\mu  (R^{-1})^{IJ} \M{1}{I}_\nu R^{IJ}
      \ \ \mbox{for} \ \ \nu<\mu \ \ . \label{defAgn3}\\[1mm]
  (R^{-1})^{IJ} \M{1}{I}_\nu R^{IJ} \Ae{2}{J}_j
      &=&  \Ae{2}{J}_j (R^{-1})^{IJ} \M{1}{I}_\nu R^{IJ}
       \ \ \mbox{ for all } \ \ \nu,j \ \ , \label{defAgnnew1}\\[1mm]
  (R^{-1})^{IJ} \M{1}{I}_\nu R^{IJ}  \Be{2}{J}_j
      &=&  \Be{2}{J}_j (R^{-1})^{IJ} \M{1}{I}_\nu R^{IJ}
       \ \ \mbox{ for all } \ \ \nu,j \ \ , \label{defAgnnew2}\\[1mm]
  (R^{-1})^{IJ} \Ae{1}{I}_i R^{IJ} \Ae{2}{J}_j
      &=&  \Ae{2}{J}_j  (R^{-1})^{IJ} \Ae{1}{I}_i R^{IJ}
       \ \ \mbox{ for } \ \ i < j \ \ , \label{defAgnnew3}\\[1mm]
  (R^{-1})^{IJ} \Be{1}{I}_i R^{IJ} \Be{2}{J}_j
      &=&  \Be{2}{J}_j  (R^{-1})^{IJ} \Be{1}{I}_i R^{IJ}
       \ \ \mbox{ for } \ \ i < j \ \ , \label{defAgnnew3a}\\[1mm]
  (R^{-1})^{IJ} \Ae{1}{I}_i R^{IJ} \Be{2}{J}_j
      &=&  \Be{2}{J}_j (R^{-1})^{IJ} \Ae{1}{I}_i R^{IJ}
       \ \ \mbox{ for } \ \ i < j \ \ , \label{defAgnnew4} \\[2mm]
  (R^{-1})^{IJ} \Be{1}{I}_i R^{IJ} \Ae{2}{J}_j
      &=&  \Ae{2}{J}_j (R^{-1})^{IJ} \Be{1}{I}_i R^{IJ}
       \ \ \mbox{ for } \ \ i < j \ \ . \label{defAgnnew4a}
\ea
\end{defn}

Now we motivate and discuss some pieces of this long definition.

First of all, the most part of defining relations are quantized
counterparts of Poisson brackets (\ref{FoRo}). Let us pick as
example equation
\be \label{AeBe}
(R^{-1})^{IJ} \Ae{1}{I}_i R^{IJ} \Be{2}{J}_i
      = \Be{2}{J}_i (R')^{IJ} \Ae{1}{I}_i R^{IJ}\ \ .
\ee
In order to recover the Poisson algebra one has to use the
standard Ansatz of quantum mechanics (axiomatized by the deformation
quantization theory).

\begin{defn}
Let $X$ and $Y$ be elements of an associative algebra defined
over formal power series in $h$.
Assume that this algebra becomes abelian for $h=0$ and
\be \label{Poi}
XY-YX=hZ_1+h^2Z_2+\dots \ \ .
\ee
Then the Poisson bracket of $X$ and $Y$ is defined as
\be
\{ X, Y\} = Z_1\ \ .
\ee
\end{defn}

Applying this definition to equation (\ref{AeBe}) we expand
quantum $R$-matrices into the power series in $h$. As a coefficient
at $h^0$ we
discover the commutator of $\Ae{1}{I}_i$ and $\Be{2}{J}_i$ which
fits to the l.h.s. of (\ref{Poi}). To evaluate the first order
in $h$ we use formula (\ref{cor}) and recover the Poisson bracket
as presented in the r.h.s. of the last equation in (\ref{examp}).

What seems to be missing in this picture is counterparts of
mutual Poisson brackets of matrix elements of the same holonomy.
Let us temporarily reintroduce these missing relations:
\ba \label{AA}
 (R^{-1})^{IJ} \Ae{1}{I}_i R^{IJ} \Ae{2}{J}_i
      &=& \Ae{2}{J}_i (R')^{IJ} \Ae{1}{I}_i (R^{'-1})^{IJ} \ \ ,
\nonumber  \\[1mm]
 (R^{-1})^{IJ} \Be{1}{I}_i R^{IJ} \Be{2}{J}_i
      &=& \Be{2}{J}_i (R')^{IJ} \Be{1}{I}_i (R^{'-1})^{IJ} \ \ , \\[1mm]
 (R^{-1})^{IJ} \M{1}{I}_{\nu} R^{IJ} \M{2}{J}_{\nu}
      &=& \M{2}{J}_{\nu} (R')^{IJ} \M{1}{I}_{\nu} (R^{'-1})^{IJ} \ \ .
\nonumber
\ea
Together with other quadratic relations of the last definition,
equations (\ref{AA}) provide a quantization of the Poisson bracket
(\ref{FoRo}). Let us mention that in fact formula (\ref{FoRo})
encodes the same number of equations as (\ref{defAgn3-}-\ref{defAgnnew4a})
together with (\ref{AA}). Quantum exchange relations look somewhat more
complicated only for notational reasons.

The key observation which one can make looking at the quadratic
relations (\ref{defAgn3-}-\ref{defAgnnew4a}), (\ref{AA}) is the presence
of the Hopf algebra symmetry in the graph algebra.
More explicitly, let $\xi$ be an element of $\G$ and
\be
\D (\xi)=\sum \xi_{\s}^1\otimes \xi_{\s}^2.
\ee
The action of $\xi$ on the generators of $\L_{g, m}$ is defined as follows:
\ba \label{act}
 \xi (M_{\nu}^I) & = &
   \sum_{\s} \t^I(\S(\xi_{\s}^1)) M_{\nu}^I \t^I(\xi_{\s}^2)
 \ \ ,
\nonumber  \\[1mm]
 \xi (A_i^I) & = & \sum_{\s} \t^I(\S(\xi_{\s}^1)) A_i^I \t^I(\xi_{\s}^2)
 \ \ , \\[1mm]
 \xi (B_i^I) & = & \sum_{\s} \t^I(\S(\xi_{\s}^1)) B_i^I \t^I(\xi_{\s}^2)
 \ \ .
\nonumber
\ea
One can continue this action to the whole of $\L_{g, m}$ using the
property of generalized derivations \cite{Sud}
\be   \label{gender}
\xi (XY) =\sum_{\s} \xi^1_{\s} (X) \xi^2_{\s} (Y)\ \ .
\ee
Here $X$ and $Y$ are arbitrary elements of $\L_{g, m}$.
Formulas (\ref{act}) and (\ref{gender}) provide a  proper generalization
for the quantized case of simultaneous conjugations (\ref{inv}).
 Quadratic exchange relations which define the algebra
$\L_{g, m}$ are invariant with respect to the action of
the quantum symmetry in the following sense. For each of them
\be
\xi (l.h.s)= \xi (r.h.s)
\ee
for any $\xi \in \G$.

There is  no surprise that we discover the quantum invariance
in the quadratic exchange relations. Indeed, they are defined
by the same set of $R$-matrices as the  quasi-triangular Hopf
symmetry algebra.

In fact, we can make the quantum description closer to the
classical one if we trade the action of $\G$ for the coaction
of the dual Hopf algebra $\G^{\ast}$. The latter is generated
by the matrix elements of $g^I\in \End(V^I)\otimes \G^{\ast}$.
Among the others they satisfy quadratic exchange relations
\be
R^{IJ} g^I g^J=g^J g^I R^{IJ}.
\ee
These are famous defining relations for the algebra
of functions on a quantum group.

It is easy to check that the mapping ${\bf g}: \L_{g, m}
\rightarrow \G^{\ast}\otimes \L_{g, m}$ defined as
\ba \label{qcon}
M_{\nu}^I & \rightarrow & (g^I)^{-1} M_{\nu}^I g^I
 \ \ ,
\nonumber  \\[1mm]
A_i^I & \rightarrow  & (g^I)^{-1} A_i^I g^I
 \ \ , \\[1mm]
B_i^I & \rightarrow & (g^I)^{-1} B_i^I g^I
\nonumber
\ea
preserves exchange relations. The formalism (\ref{qcon}) is more
transparent. In particular, it has been extensively used in
\cite{BuRo} for description of lattice gauge models with quantum
gauge group. However, in view of the generalizations for  the quasi-Hopf
algebras, we stick to the more formal definition of the quantum invariance
which involves the action of $\G$
rather than the coaction of $\G^{\ast}$. Let us stress that these
two approaches
are completely equivalent.

Now we can explain the origin of the first three relations in the
definition of the graph algebra. Before quantization holonomies
$M_\nu, A_i$ and $B_i$ take values in the group $G$. This implies
that matrix elements of the same holonomy in different representations
are not algebraically independent. More explicitly, a product of two
matrix elements of some holonomy matrix in two different representations
may be decomposed into a sum of certain matrix elements by means of
Clebsch-Gordon maps:
\be \label{CGcl}
\M{1}{I}_\nu \M{2}{J}_\nu
      = \sum C_0^a[IJ|K]^* M_\nu^K C_0^a[IJ|K]\ \ .
\ee
Here the Clebsch-Gordon maps with the subscript $0$ refer
to the undeformed Lie group $G$.

Relations of the type (\ref{CGcl}) do not hold in the quantized
algebra. They would contradict e.g. the quadratic exchange relations
(\ref{AA}). The first three relations in the definition of the graph
algebra provide a proper substitute for (\ref{CGcl}). They are chosen
to be consistent with the
symmetry Hopf algebra action (\ref{act}), (\ref{gender}) on the
graph algebra. Equations (\ref{AA}) follow from the multiplication
laws (\ref{defAgn}-\ref{defAgn2}). That is why we do not include
(\ref{AA}) into the basic definition.

\subsection{Integration and $*$-properties of graph algebras}

This subsection is devoted to the $*$-operation and the integration
measure on the graph algebras. These two objects are useful technical
tools in the analysis of the representation theory. Also, the
$*$-operation is important for the physical interpretation as
an algebra of observables of a quantum system is always a $*$-algebra.

Let us recall that we assume consistency of the co-product for the
Hopf symmetry algebra with the special kind of $*$-operation which reverses
the order in the tensor product (see subsection 3.1).
As a consequence
of this choice the algebra $\L_{g, m}$ {\em is not} equipped with
a natural $*$-operation. However, we can save the situation using
the following trick.

The algebra $\G$ acting on $\L_{g, m}$, we can define a semi-direct
product $\S_{g, m}$  of these two objects.
It is generated by the elements $\xi
\in \G$ and by the elements of $\L_{g, m}$.
In order to describe the commutation relations of the symmetry
generators and the matrix elements of quantized holonomies
it is convenient to introduce the generating matrices
$\mu^I(\xi)\in End(V^I)\otimes \G$ for the symmetry algebra:
\be
\mu^I(\xi)= (\t^I \otimes id) \D (\xi)\ \ .
\ee
The cross relation between $\xi$ and quantum holonomies
in $\S_{g, m}$ look like
\ba
  \mu^J(\xi) M_\nu^J &=& M_\nu^J \mu^J(\xi)\ \ ,
  \label{defAgn4a} \\[2mm]
  \mu^J(\xi) A_i^J = A_i^J \mu^J(\xi) \ \ & , & \ \
  \mu^J(\xi) B_i^J = B_i^J \mu^J(\xi) \ \ . \label{defAgn4}
\ea

The semi-direct product $\S_{g, m}$ is already a $*$-algebra.
The $*$-operation on $\G$ coincides with the genuine $*$-operation
of the symmetry algebra. It is continued to the quantum holonomies
as \cite{AGS1}:
\ba  \label{ast}
  (M_\nu^I)^* &=& \s_\k(R^I (M_\nu^{I})^{-1} (R^{-1})^I)\ \ , \nn \\[1mm]
  (A_i^I)^*   &=& \s_\k( R^I (A_i^{I})^{-1}   (R^{-1})^I)\ \ ,  \\[1mm]
  (B_i^I)^*   &=& \s_\k( R^I (B_i^{I})^{-1} (R^{-1})^I )\ \ . \nn
\ea
Here we introduced
$(M_\nu^{I})^{-1}, (A_i^{I})^{-1}, (B_i^{I})^{-1} \in
\End(V^I) \o \S_{g,m}$ being the unique solutions of $M_\nu^I
(M_\nu^{I})^{-1} = e^I = (M_\nu^{I})^{-1} M_\nu^I$,\
$A_i^I  (A_i^{I})^{-1} = e^I = (A_i^{I})^{-1} A_i^I$ and $B_i^I
(B_i^{I})^{-1} = e^I = (B_i^{I})^{-1} B_i^I$.
The symbol $\s_\k$ stays for  the automorphism of  $\S_{g, m}$
obtained by conjugation with the unitary
element $\k \in \G$,
$$\s_\k (F) = \k^{-1} F \k$$
for any $F$.

Another important object which may be introduced for $\S_{g, m}$
is an invariant integration measure. We define a linear
functional $\om: \S_{g, m}\rightarrow {\bf C}$ as
\be  \label{integ}
\om(\ M_1^{I_1} \dots B_g^{I_{m+2g}}
\xi \ )=\e(\xi) \prod_{s=1}^{m+2g} \d_{I_s, 0}\ \ .
\ee
Formula (\ref{integ}) is defined on the set of monomials which span
the algebra $\S_{g, m}$. So, $\om$ is continued to the whole algebra
as a linear functional. The integral (\ref{integ}) may be restricted
to the algebra $\L_{g, m}$.
There  it furnishes the quantum analog of the
multidimensional Haar measure on $G^{m+2g}$.
In particular, one can formulate the invariance
of $\om$ as
\be
\om(\xi(X))=\e(\xi) \ \om(X)
\ee
for any $\xi \in \G$ and $X \in \L_{g,m}$.

Usually an interplay between the $*$-operation and the integration
is the positivity property of the Hermitian scalar product defined by the
integration functional
$$ (X, Y)=\om(X^*Y)\ \ . $$
However, this property never holds on $\S_{g, m}$ as $\om$ always
has a big kernel in $\G$. On the other hand, one can not formulate
positivity on $\L_{g, m}$ as this is not a $*$-algebra.
So, instead of the usual positivity we formulate the following
substitute.
\begin{theo} {\em (Positivity \cite{AGS1})}
\label{positivity}
Assume that the quantum dimensions of all irreducible representations
of the Hopf symmetry algebra $\G$ are strictly positive
\be
d_I > 0 \ for \ every \ I \ \ .
\ee
Then the restriction of the integral $\om$ to $\L_{g,m}$ is
positive in the following sense
\begin{eqnarray}  \label{posit}
\om(X^*X) & \geq & 0 \mbox{ for all }  X\in \L_{g, m}
\nonumber \\[1mm]
\mbox{ and } & & \om(X^*X)=0 \Rightarrow X=0\ \ .
\end{eqnarray}
\end{theo}
One can find a proof of this theorem in \cite{AGS1}.

Let us remark that the element $X^*$ does not belong to
$\L_{g, m}$. So, we need the bigger algebra $\S_{g, m}$ in order
to make sense of  (\ref{posit}). In fact, this kind of positivity
on the level of graph algebras will be sufficient to provide
positivity of the integration measure on the quantized moduli
algebras.

\subsection{The moduli algebra}

The construction of the moduli space starting from
$G^{2g+m}$ involves two additional steps. First, one has to
restrict the consideration to the subspace of invariant functions.
Next, one imposes the conditions (\ref{MMC}). We should find quantum
counterparts of both operations.

It is straightforward to generalize the first step.

\begin{defn} {\em (Algebras $\A_{g,m}$)}  \label{defAgm}
$\A_{g,m}$ is defined as a subalgebra of elements $A \in \L_{g, m}$
which are invariant with respect to the natural action of $\G$, i.e.
$ \xi (A) = A \e (\xi)$ for all $A \in \A_{g,m}$ and $\xi \in \G$.
\end{defn}

Since elements $\xi \in \G$ act trivially on $\A_{g,m}$, the
semi-direct product of $\G$ and $\A_{g,m}$ coincides with the
usual Cartesian product $\G \ti \A_{g,m}$ and hence the $*$-operation
on $\S_{g,m}$ furnishes a $*$-operation on $\A_{g,m}$.
Then the modified positivity which we defined on $\L_{g, m}$
ensures the usual positivity property of $\om$ on $\A_{g, m}$.
Let us also remark that elements of $\A_{g,m}$ are linear combinations
of expressions of the form:
$$ tr_q^J \left( C^q_1 [I_1, \dots, I_{2g+m}|J] M^{I_1}_1
  \dots M^{I_m}_m A_1^{I_{m+1}}\dots B_g^{I_{2g+m}}
   C^q_2 [I_1, \dots, I_{2g+m}|J]^*\right)\ \ . $$
Here $tr_q$ is the $q$-trace, $C^q_1, C^q_2$ are intertwining
operators for the Hopf algebra action. For generic values of
$q$ one can establish an isomorphism of the linear spaces
$\A_{g,m}$ and $\A^0_{g,m}$. The latter is  the space of conjugation
invariant
analytic functions on $G^{2g+m}$. Obviously, this isomorphism
can not be lifted to the level of algebraic structures as the
space of functions is abelian whereas $\A_{g,m}$ is not.

To perform the second step of the reduction to the moduli
space we single out a particular
graph algebra  corresponding to one marked point on a
Riemann surface.

\begin{defn} {\em (Loop algebra $\L$)}
The loop algebra $\L$ is an associative algebra isomorphic to
the graph algebra $\L_{0,1}$. It is generated by matrix elements
of the monodromies $M^I\in End(V^I)\otimes \L$.
\end{defn}

Closely related to $\L$ is the abelian
{\em fusion ({\rm or} Verlinde) algebra}.

\begin{defn} {\em (Fusion (or Verlinde) algebra)}
Let $\G$ be  semi-simple (quasi-) Hopf $*$-algebra.
Its fusion (or Verlinde) algebra $\V$ is an abelian $*$-algebra
spanned by the set of generators $c^I$. Here $I$ runs through
the set of all irreducible representations of $\G$. The multiplication
law and the $*$-operation in $\V$ are defined as follows:
\be \label{cVer}
    c^I c^J = \sum N^{IJ}_K c^K \ \ \mbox{and} \ \
    (c^I)^* = c^{\bar I}\ \ .
\ee
\end{defn}

One can employ the $q$-traces $tr^I_q$ to construct
central elements $c^I \in \L$
from the monodromies $M^I$
\be \label{c}
    c^I = \k_I tr^I_q( M^I ) \ \ .
\ee
In \cite{AGS2} we have demonstrated that the elements (\ref{c})
generate the Verlinde algebra.

It is important that under certain conditions (for details see
Section 5) the representations of the loop algebra $\L$
 and of the corresponding
Verlinde algebra $\V$ may be labeled by the same set of labels as the
representations of the symmetry algebra $\G$. Here we describe
the representations of the Verlinde algebra by explicit formulas.

Let us introduce a (possibly infinite) matrix
\be
s^{IJ}=(tr^I_q \o tr^J_q) (R'R)
\ee
with rows and columns labeled by the representations
of $\G$. Equations
\be \label{nuJ}
\vartheta^J(c^I)= \frac{s^{IJ}}{d^J}
\ee
define the set of representations of the Verlinde algebra.
We shall see that in a quite general situation  this set
of representations is complete.

It is convenient to introduce a special notation for the
relations
\be \label{PhiJ}
\Phi^J=\{ c^I=\vartheta^J(c^I)\}
\ee
which restricts the generators $c^I$ to a certain representation.
Imposing these central relations
we may get ideals in both $\V$ and $\L$.

Returning to arbitrary values of $g$ and $m$ we introduce
$m+1$ embeddings of $\L$ into $\L_{g,m}$ defined
by
\begin{eqnarray} \label{Mi}
e_\nu(M) & =& M_\nu \  \mbox{ for } \  \nu =1,\dots , m\ \ ;
\nonumber \\[1mm]
e_0(M) & = & [B_g, A_g^{-1}]\dots [B_1, A_1^{-1}] M_m\dots M_1\ \ .
\end{eqnarray}
If necessary, these embeddings may be lifted to the
corresponding semi-direct products with the symmetry algebra.

The embeddings (\ref{Mi}) provide a set of elements in  $\A_{g,m}$
\be
c^I_\nu=\k_I tr^I_q( M^I_\nu)\ , \ \nu =0,\dots , m\ \ .
\ee
It was shown in \cite{AGS2} that all of them belong
to the centre of $\A_{g,m}$. In particular, they commute
with each other. Now we are ready to define the moduli algebra.

\begin{defn} {\em (Moduli algebra)}  \label{defmodalg1}
Let $\G$ be a quasi-triangular Hopf symmetry algebra, $\Sigma_{g,m}$
be a closed oriented 2-dimensional surface of genus $g$ with $m$
marked points and $I_1,\dots , I_m$ be a set of $m$ irreducible
representations of $\G$ assigned to the marked points. The
moduli algebra ${\cal M}_{g,m}^{\{ I_\nu\}}$ is defined
by these data as a quotient of the invariant algebra $\A_{g,m}$,
\be
{\cal M}_{g,m}^{\{ I_\nu\}}=
\A_{g,m}/ \{ \Phi^0(M_0), \Phi^{I_\nu}(M_\nu), \nu=1,\dots , m\}.
\ee
Here $0$ labels the trivial representation of $\G$. The
moduli algebra $\MC^{I_\nu}_{g,m}$ inherits the *-operation
and the positive integration functional $\om$ from $\A_{g,m}$.
\end{defn}

In fact, relations $\Phi^I$ are proper quantum counterparts
of fixing the eigenvalues of the corresponding quantum holonomy.
In particular, the set of relations $\Phi^0$ is equivalent
to $M^I=\k_I^{-1}e^I$ for any $I$. The scalar factor $\k_I$
gives a `quantum correction' to the classical flatness
condition $M^I=e^I$.

\subsection{The finite-dimensional case}

As we already mentioned, we are mostly concerned with the case
of $q$ being a root of unity. Then  one
can not view the moduli algebra as a deformation of
the algebra of functions on the moduli space. Instead, one can reverse
the logic  and reinterpret the definition
of the graph algebra. The latter is defined by choosing a Riemann
surface with marked points, a ribbon Hopf $*$-algebra and
a set of  representations of this algebra, one representation for
each marked point. The main ingredient in these data is the
ribbon Hopf algebra. Instead of looking what happens to the
the moduli algebra at roots of unity we look at the symmetry
Hopf algebra, define it at roots of unity and thus
induce a new definition of the moduli algebra.
Here we follow this strategy. However, the relation between
the moduli algebras at generic $q$ and at roots of unity still
needs to be clarified.

The first important observation is that at roots of unity
quantum universal enveloping algebras have a big centre \cite{KC}.
A natural object is a quotient
over certain central relations which is already finite
dimensional \cite{Lus}. This is a motivation to consider
symmetry Hopf algebra with only finite number of
irreducible representations.

To ensure the positivity of $\om$ (see Subsection 3.3) we
require the quantum dimensions of all irreducible representations
to be strictly positive
\be
d_I >0 \ \ \mbox{ for all }  \ \ I\ \ .
\ee
In fact, this requirement may be satisfied only at roots of unity.
Let us remark that there are still indecomposable representations
with vanishing quantum dimensions. We will deal with them later
in this subsection.

As the number of representations is finite,
one can introduce a normalization constant
\be
\N \equiv     (\sum_K d_K^2)^{-1/2} < \infty \ .
\ee
It is useful to have a
matrix $S^{IJ}$ which differs from $s^{IJ}$ by this scalar
factor:
\be
     S_{IJ} \  \equiv \  \N (tr^I_q \o tr^J_q) (R'R) \ \ . \label{S}
\ee
We assume that {\em the matrix $S$ is invertible}.
A number of standard properties of $S$ can
be derived from the invertibility (and properties of the ribbon
Hopf-*-algebra). We list them here without further
discussion. Proofs can be found e.g. in
\cite{FrGa}.
\ba
    S_{IJ} = S_{JI} \ \ \ & , &
    \ \ \ S_{0J} = \N d_J \ \ , \nn \\[2mm]
    \sum_J S_{IJ} \overline{S_{KJ}} = \d_{IK} \ \ & , & \ \
    \sum_J S_{IJ} S_{JK} = C_{IK}\ \ ,  \label{Sprop}\\
    \sum_K N^{IJ}_K S_{KL} &=&  S_{JL} S_{IL} (\N d_L)^{-1} \nn
\ea
with $C_{IJ} = N^{IJ}_0$. For the relations in the second
line, the existence of an inverse of $S$ is obviously necessary.
Invertibility of $S$ is also among the defining
features of a modular Hopf-algebra in \cite{ReTu2}.
The last equation in the set (\ref{Sprop}) is usually referred
to `diagonalization' of fusion rules \cite{Ver}.

In the finite-dimensional situation on can consider
certain linear combinations $\chi^K$ of the $c^I$,
$$ \chi^K = \N d_K S_{K I} c^{\bar I} \ \ . $$
Here $\N = (\sum d_I^2)^{-1/2}$ and $S_{K I}$ are components
of the $S$-matrix (\ref{S}). The projectors  $\chi^K$
are central, orthogonal projectors, i.e. $(\chi^K)^*
= \chi^K$ and $\chi^K \chi^L = \d_{K,L} \chi^K$ ( see \cite{AGS2}).
They represent characteristic functions for the $1$-dimensional
ideals in the Verlinde algebra corresponding
to the representations considered in the previous subsection.

One serious technical problem which arises when we consider
quantum universal enveloping algebras at roots of unity
is the fact that they loose semi-simplicity. In particular,
finite-dimensional Hopf algebras under discussion are not
semi-simple. This problem may be cured in two different
ways. One way is to work with non semi-simple algebras.
Then the resulting moduli algebra is expected to have a
big ideal formed by functions which include matrix elements
of indecomposable representations. If one wants to introduce
the moduli algebra as a $*$-algebra this ideal should be factored
out. We refer to this operation as {\em truncation}.

Another way is to force the universal enveloping algebra to be
semi-simple \cite{MSIII}. Technically, one moves to the class
of quasi-Hopf algebras and  relaxes the axiom
\be
\Delta ( e) = e \otimes e\ \ .
\ee
This defines a class of weak quasi-Hopf algebras. Then
a non semi-simple Hopf algebra may be replaced by  a weak
quasi-Hopf algebra in the following way. One factors out
all the elements which vanish in all irreducible representations.
The resulting object is by definition semi-simple. All its
representations are completely reducible. In the case
of quantized universal enveloping algebras these are
so called physical representations related to Conformal Field Theory.
New semi-simple symmetry algebras are  also called
{\em truncated}.
They  do not  satisfy the axioms of the Hopf category. In particular,
the co-multiplication still may be defined but it is not co-associative.
We refer the reader to the original paper \cite{MSIII} for a more
detailed account. It is a conceivable conjecture that looking at the moduli
algebras corresponding to truncated quantum symmetries is equivalent
to dealing with truncated moduli algebras defined by a non semi-simple
quantum symmetry.

Admitting weak quasi-Hopf algebras we have to prove that moduli
algebras may be defined by these data. This has been done in
\cite{AGS1}, \cite{AGS2}. The qualitative difference is that the graph
algebra now is only quasi-associative. Indeed, it is designed
as an algebra of the objects covariant with respect to the symmetry
action. As the co-multiplication of the symmetry algebra is
quasi-co-associative, the multiplication of tensors is forced to be
quasi-associative. However, the moduli algebra is always  an
associative algebra.

All our results  are valid for truncated weak quasi-Hopf
symmetry algebras. This includes the most interesting cases of $U_q(\sg)$
at roots of unity which correspond to quantization of the Chern-Simons
theory for integer values of $k$ and for the compact Lie group $G$.
However, for pedagogical reasons we work throughout the paper
with the unrealistic case of a semi-simple Hopf algebra with finite
number of irreducible
representations. A more sophisticated version of the same
calculations goes through for the case of weak quasi-Hopf algebras.
In the end of the paper we comment on  the most important changes in
the consideration.

\section{Summary of Results}

\subsection{Representation theory of the moduli algebra}

The moduli algebra $\MC_{g,m}^{\{ I_\nu\}}$ corresponding
to the truncated universal enveloping algebra $U_q(\sg)$ is
supposed to coincide with the algebra of observables of the
Chern-Simons theory. Naturally, its $*$-representations may
be considered as candidates for the role of the Hilbert
space in this model.

Assuming that we construct the moduli algebra starting from a
semi-simple quasi-triangular (quasi)-Hopf algebra with
only finite number of representations we arrive at
the following result.

\begin{theo} {\em (Representations of the moduli algebra)}
For any set of representations $I_1,\dots , I_m$ assigned
to the marked points there exists a unique irreducible
*-representation of the moduli algebra ${\cal M}_{g,m}^{\{ I_\nu
\}}$ which acts in the space
\begin{eqnarray} \label{Mrep}
W_g^0(I_1,\dots , I_m) & =& \Inv (V^{I_1}\otimes \dots V^{I_m}\otimes
\Re^{\otimes g})
\\[1mm]
\mbox{ where } & & \Re =\oplus_I V^I\otimes V^{\bar{I}}\ \ .
\nonumber
\end{eqnarray}
Here {\it Inv} stays for invariant subspace with respect to the
natural action of the symmetry algebra.
\end{theo}

It is remarkable that the moduli algebra has a unique representation.
Thus, it is identified with the full matrix algebra with the natural
$*$-operation and may be regarded as an algebra of observables
of some quantum mechanical system with a finite-dimensional space
of states. Apparently, this is the case of the Chern-Simons theory
in the Hamiltonian formulation.

The space $\Re$ which enters the definition of $W_g^0(I_1,\dots , I_m)$
is an analog of the regular representation of  a finite or a
compact Lie group.

We shall show explicit formulas for the action
of $\MC_{g, m}^{\{ I_\nu\}}$ in the representation
(\ref{Mrep}). Let us only remark that there is an important
difference between the moduli algebras of zero and nonzero
genus. The first ones involve only $R$-matrices in the expressions
for the matrix elements of their representations. However, when the
genus of the surface is nontrivial, the knowledge of Clebsch-Gordon
maps is required.

\subsection{The action of the mapping class group on the
moduli algebra}

There is an important structure on  the moduli space
which we did not touch before. The moduli space
of flat connections on a surface $\Sigma_{g, m}$
carries the action
of the pure mapping class group $PM(g, m)$ of the
surface. The pure mapping class group is a subgroup
of the mapping class group $M(g, m)$ which preserves
the order of marked points.
$PM(g, m)$ acts by automorphisms
of the fundamental group of the surface. This action
lifts to $Hom(\pi ,G)$ as
\be
\rho^{\eta}(x)=\rho(x^{\eta})  \ \ \mbox{ for all }\ \
   x \in \pi_1 (\Sigma_{g,m}) \ \ .
\ee
Here $\rho$ is a representation of $\pi$ in $G$,
$x$ is an element of the fundamental group and $\eta$
is an element of the mapping class group. As
$\eta$ defines an automorphism of $\pi$, this action
descends to the moduli space.

In the course of quantization symmetries are usually
very important as they give us guidelines which features
of the classical theory should be preserved by quantization.
The action of the mapping class group on the moduli space of flat
connections preserves the symplectic structure. So, we should
expect that this action lifts to the quantized algebra of functions
on the moduli space.  This is indeed the case.
We need some more notations to describe this action.

The pure mapping class group may be generated by so called
Dehn twists. A Dehn twist corresponds to a circle on
a Riemann surface. The mapping class group transformation includes
cutting the surface along this fixed circle, relative
rotation of the boundaries of the cut by the angle of $2\pi$ and gluing
the sides of the cut back together. Thus one defines
a smooth mapping of the surface into itself which does
not belong to the connected component of the identical mapping.

As a first step we define a Verlinde subalgebra in the moduli
algebra for each circle on a Riemann surface. A circle
on a surface (or, more exactly, its homotopy type) defines
a conjugacy class in the fundamental group. Let us take
any element $x$ of this conjugacy class and represent it
as a word in the generators of the standard basis $l_\nu,
a_i, b_i$. As an example let us pick up an element
\be
x=l_1 b_1 a_1\ \ .
\ee
Now we define a bunch of quantum holonomy matrices
via
\be
X^I=\k_I^{-2} M^I_1 B^I_1 A_1^I\ \ .
\ee
Here $I$ is as usual a representation of the symmetry algebra.
In the standard way we define the generators of the Verlinde
algebra
\be
c^I(x)=\k_I tr_q X^I
\ee
corresponding to the circle $x$.
One can easily repeat this procedure for an arbitrary
circle on the surface.
The choice of the
representative in the conjugacy class appears to be irrelevant
for the definition
of $c^I$ due to the properties of the $q$-trace.
We refer to Section 9 for the general rule of counting
extra $\k^I$ factors.

The last ingredient which we need is a particular
element in the Verlinde algebra defined as
\be
\hh (x)=\sum_I v_I^{-1} \chi^I(x)\ \ .
\ee
Here projectors $\chi^I$ are constructed of  $c^I$ as in the previous
section.

We collect
the main results concerning the action of the mapping
class group in the following theorem.

\begin{theo}
The moduli algebra ${\cal M}_{g, m}^{\{ I_\nu\}}$ corresponding to a
surface of genus $g$ with $m$ marked points carries a
natural action of the pure  mapping class
group $PM(g, m)$ which preserves the order of the marked points.
$PM(g,m)$ acts by inner automorphisms of the moduli algebra.
The unitary generator of a particular Dehn twist defined by a circle
$x$ on the surface  is given by $\hh (x)$. This gives rise to a
projective representation of the mapping class group which
is unitary equivalent to the one described in \cite{ReTu}.
\end{theo}

\subsection{Comparison to other approaches}

We can compare the Combinatorial quantization approach
of this paper to two other quantization schemes.

Let us recall that the Hilbert space of the Chern-Simons theory
has been identified with the space of conformal blocks in the WZW
Conformal Field Theory corresponding to the same group $G$ and
with the same value of the coupling constant $k$ \cite{Wit1}.
The guess about the relation of the CS and WZW systems has been
explained in the framework of Geometric Quantization \cite{ADW}.
Technically speaking the space of conformal blocks may be
characterized as a space of solutions of certain linear differential
equations. In the simplest case of the 2-dimensional surface
being a sphere with marked points, these equations were suggested
in \cite{KZ} and usually called Knizhnik-Zamolodchikov
(KZ) equations.

One can notice that the dimension of the representation
space (\ref{Mrep}) coincides with the Verlinde formula for the
dimension of the space of conformal blocks. The natural isomorphism
between these two spaces emerges when one considers certain
asymptotics of the KZ-equations \cite{Dri2}. This is one of the
aspects of the relation between quantum groups and Kac-Moody
algebras discovered in \cite{KaLu}. The particular case of roots
of unity has been worked out in \cite{Fin}. We expect that
the general results of \cite{KaLu,Fin} imply the identification
of the geometric quantization and quantum group pictures for the
space of conformal blocks. However,
the current status of this
construction is not clear to us.

Another approach which we follow quite closely is deformation
quantization. The $r$-matrix presentation of the Poisson brackets
on the moduli space is especially designed to make deformation
quantization easy. It is natural to conjecture that for generic
values of the deformation parameter $q$ the moduli algebra
${\cal M}_{g,m}^{\{ I_\nu\}}$ provides a deformation
quantization of the algebra of functions on the moduli space.

As was recently discovered \cite{Fed}, the deformation quantization
is uniquely defined by choosing a symplectic connection on the
phase space. It is an intriguing question what kind of connection
on the moduli space chooses the Combinatorial quantization.

\part{Representation Theory of the Moduli Algebra}

Our basic strategy in this part on the representation theory
of graph- and moduli algebras is to study very simple building
blocks of the graph algebra first and then to put the pieces
together for the full theory to emerge. The simple building
blocks are the {\em loop algebra} $\L = \L_{0,1}$ (section 5)
and the {\em $AB$- (or handle-) algebra} $\T= \L_{1,0}$ (section 7).
A bunch of loop algebras $\L $ may be combined into a
{\em multi-loop algebra} $\L_m = \L_{0,m}$ (section 6).
{}From this one obtains the moduli algebras associated with
a punctured sphere. Section 8 concludes the representation
theory of moduli algebras with a complete description for
arbitrary genus $g$. The relation of this theory with
representations of the mapping class group is explained
in Section 9.

\section{The Loop Algebra and its Representations}
\setcounter{equation}{0}

This section contains an extensive treatment of the {\em loop
algebra} $\L$ which will be the most fundamental building block
in all the subsequent discussion. The main aim is to
develop a complete representation theory for $\L$. In
passing we note some results on Gauss decompositions
in the third subsection.

\subsection{The loop algebra $\L$}

The {\em loop algebra} $\L$ already appeared in Section 3.2
and 3.4 as a special example of graph algebras, namely $\L_{0,1}$.
It is generated by
matrix elements of the monodromies $M^I \equiv M^I_1 \in
\End(V^I) \o \L$. In this simple case of a graph algebra,
monodromies only have to obey {\em functoriality},
\be    \label{defM}
 \M{1}{I} R^{IJ} \M{2}{J}  =  \sum C^a[IJ|K]^* M^K C^a[IJ|K]
 \ \ .
\ee
Such relations were discussed at length in Section 3.2. For
the covariance properties of monodromies $M^I$ and the action
of the $*$-operation on $\S_{0,1} \supset \L$,  the reader is
referred to Sections 3.2, 3.3.

To motivate the results we are about to see in this section,
let us pick up some traces that were laid in the first part.
It was already noticed in Section 3.2 that {\em functoriality on the
loop} (eq. (\ref{defM})) determines the following exchange
relations for the monodromy
\be \label{Mex}
   (R')^{IJ} \M{1}{I} R^{IJ} \M{2}{J}  =
    \M{2}{J} R^{IJ} \M{1}{I} (R')^{IJ}  \ \ .
\ee
Relations of this form were found to describe the quantum
enveloping algebras of simple Lie algebras \cite{ReSTS}, which
are our main examples for the quantum symmetry $\G$. Thus --
following the ideology of \cite{FRT} -- one expects that the
deeper investigation to be carried out below will reveal
an isomorphism between $\L $ and $\G$ (at least up to some
subtleties). This expectation gains further support from our
discussion in Sections 3.4 and 3.5 where we found that the
following elements
$$
    \chi^K = \N d_K S_{K\bar I} c^I
           = \N d_K \k_I S_{K \bar I} tr^I_q(M^I)
$$
form a set of orthogonal projectors in the center of the loop
algebra. Since these {\em characteristic projectors} $\chi^K$
are labeled by the same index as the irreducible representations
of $\G$, a close correspondence between representations of $\G$
and $\L$ is quite plausible.

These remarks, however, have to be taken with a little grain
of salt. For the isomorphism of $\G$ and $\L$ to hold, we
have to restrict ourselves to the unrealistic case of a
finite-dimensional algebra $\G$ without truncation. If the tensor
product of representations of the symmetry algebra $\G$ is
truncated, the linear dimension of $\L$ is strictly
smaller than the dimension of $\G$, thus making an isomorphism
of the two spaces impossible. Still a close relation between
their representation theories will persist so that the lesson
we learn here in the non-truncated case is not in vain. We
will return to this discussion only in Section 10 when we
explain the adjustments which are required to treat truncated
structures.

\subsection{Representation theory of the loop algebra $\L $}

Now let us turn to our the main result in this section.
We wish to consider the representation theory of $\L$. Since
we suspect a close relation between $\G$ and $\L$, it should
be possible to construct representations of $\L$ in the same
spaces $V^{I}$ in which we already have  representations
$\t^I$ of the symmetry algebra $\G$. The following theorem
gives a concrete formulation  of this idea.

\begin{theo} {\em (Representations of the loop algebra)}
\label{thMrep}
The loop algebra $\L$ has a series of representations $D^I$
realized in the representation spaces $V^{I}$ of the underlying
quasitriangular Hopf algebra $\G$. In such a representation,
the generators of $\L$ can be expressed as
$$
D^I(M^{J})  =  (\k_J)^{-1} (R'R)^{JI}\ \ .
$$
$D^I$ extends to a *-representation of the semi-direct product
$\S_{0,1} \equiv \L \sd \G$ by means of the formula
$$
D^I(\xi)    =  \t^I(\xi) \ \ \mbox{ for all }\
             \  \xi \in \G\ \ .
$$
Compatibility with the $*$-operation on $\S_{0,1}$ means in
particular that $ D^I(M^*) = (D^I(M))^*$ for all $M \in \L$.
\end{theo}

\noindent
{\sc Proof:} To prove consistency with the multiplication rule
(\ref{defM}) let us evaluate the l.h.s. of (\ref{defM}) in the
representation $D^L$
\ba
 D^L( \M{1}{I} R^{IJ} \M{2}{J} )
   &=& (\k_I \k_J )^{-1}
        (R'_{13} R_{13} R_{12} R'_{23} R_{23})^{IJL}  \nn \\[1mm]
   &=& (\k_I \k_J )^{-1}
        (R'_{13} R'_{23} R_{12} R_{13} R_{23})^{IJL}  \nn \\[1mm]
   &=& (\k_I \k_J )^{-1}
        (R_{12} R'_{23} R'_{13} R_{13} R_{23})^{IJL}  \nn \\[1mm]
   &=& (\k_I \k_J )^{-1}
        (R_{12} (\D \o id)(R'R))^{IJL}\ \ . \nn
\ea
We used the Yang Baxter equation for $R$ twice and inserted
quasi-triangularity in the last line. Now one proceeds with the
help of equation (\ref{complete}).
\ba
   &=& (\k_K)^{-1} \sum  C^a[IJ|K]^* C^a[IJ|K]
         ((\D \o id)(R'R))^{IJL}\nn \\[1mm]
   &=& \sum (\k_K)^{-1}   C^a[IJ|K]^*
         (R'R)^{KL} C^a[IJ|K] \nn \\[1mm]
   &=& D^L \left( \sum C^a[IJ|K]^* M^K C^a[IJ|K] \right)
       \nn  \ \ .
\ea
To see that $D^I$ can be extended to $\S_{0,1}$
we have to check consistency of $D^I$ with the covariance
properties (\ref{defAgn4a}) of monodromies. Acting with $D^I$
on eq. (\ref{defAgn4a}) we obtain
$$ (\t^J \o \t^I)(\D (\xi)) \k_J^{-1} (R'R)^{JI} =
   \k_J^{-1} (R'R)^{JI} (\t^J \o \t^I)(\D (\xi)) $$
which holds due to the standard intertwining property of $R$.
So let us finally show that $D^I$ is a $*$-representation.
Because $(M^J)^{-1}$ is represented by $\k_J
((R'R)^{-1})^{JI}$ one gets
\ba
D^I((M^J)^*) & = & \k_J \k_I^{-1} R^{JI} ( R^{-1} (R')^{-1} )^{JI}
                  (R^{-1})^{JI} \k_I \nn \\
             & = & \k_J ((R')^{-1} R^{-1})^{JI} = (\k_J^{-1}
                   (R'R)^{JI})^* \nn \\
             & = & (D^I(M^J))^* \ \ . \nn
\ea
This concludes the proof of the proposition.

Having constructed a series of representations of $\L$, it is
instructive to look at the central elements $c^J$ (given by eq.
(\ref{c})) and evaluate them in these representations. Using
the definition (\ref{S}) of the matrix $S_{JI}$ we can express
the value of the central element $c^{J}$ in the
representation $D^{I}$ as follows
$$
     D^I(c^J)  =\frac{S_{JI}}{ \N d_I}\ \ .
$$
Evaluation of the elements $\chi^K = \N d_K S_{K\bar J} c^J$
shows that they are characteristic projectors for our
representations,
\be      \label{charJ}
  D^I(\chi^K) = \frac{d_K}{d_I} S_{K J}
    S_{\bar J I} =  \delta_{K,I}\ \ .
\ee
Actually the $\chi^K$ provide a complete set of minimal
central projectors in $\L$:

\begin{lemma} \label{irredlemma}
The set of representations $\{ D^I\} $ is faithful on
$\L$, i.e. for any nonzero element $X \in \L$ there is at
least one label $I$ such that $D^I(X)$
is nonzero. In absence of truncation this implies that
the representations $D^I$ are irreducible.
\end{lemma}

\noindent
{\sc Proof:} Let $\d_I$ denote the dimension of $V^I$.
In order to prove the irreducibility we have to investigate the
space $\De^I \subset \End(V^I)$ obtained as the image of $\L$
under $D^I$. The assertion of the lemma holds, if $\De^I$
has the dimension $\d_I^2$, i.e. if $\De^I$ is the full
matrix algebra on $V^I$. It will be fundamental to notice
that the space $\End(V^I)$ carries a representation $ad^I$
of the symmetry Hopf-algebra $\G$,
$$
    ad^I(\xi) b^I \equiv  \sum \t^I(\xi^2_\s) b^I \t^I(\S
   (\xi^1_\s)) \ \ \mbox{ for all } \ \ b^I \in \End(V^I) \ \ .
$$
With respect to this action, the space $\End(V^I)$  decomposes
into a direct sum of subspaces one for every irreducible
representation in the decomposition of the tensor product
$$
   (\t^{\bar I} \bo \t^I)
   \cong \bigoplus N^{I\bar I}_K \t^K\ \ .
$$
The corresponding irreducible subspaces $S^I(K,a)$ of $\End(V^I)$
are labeled by pairs $K,a$ such that $C^a[I \bar I | K] \neq 0$.
Now let us fit the space $\De^I$ of representation matrices of
$\L$ into this picture. Elements in $\De^I$ are obtained from
the components of
$$ D^I(M^J) = \k_J^{-1} (R'R)^{JI}  \ \ ,  $$
where $J$ runs through all possible labels. If the standard
intertwining relation $\D(\xi) R'R = R'R \D(\xi)$ is written
in the somewhat non-standard form
$$
   \sum (e \o \xi^2_\s) R'R (e  \o \S(\xi^1_\s)) =
   \sum (\S(\xi^1_\s) \o e) R'R (\xi^2_\s \o e) \ \ ,
$$
one concludes that $\De^I \subset \End(V^I)$ is invariant under
the action $ad^I$ of the algebra $\G$. So the before mentioned
decomposition of $\End(V^I)$ induces a decomposition of the
subspace $\De^I \subset \End(V^I)$.

Elements in the irreducible subspaces $S^I (K,a) \cap \De^I$ of
$\End(V^I)$ are constructed from  representations matrices of the
loop algebra $\L$ as
\be \label{kaform}
    D^I ( C[I \bar I |0] M^J (R')^{I \bar I}
    C^a[I \bar I|K]^*)\ \ .
\ee
We wish to show that for every pair $K,a$ such that $C^a[I
\bar I| K] \neq 0$ at least one element of the form (\ref{kaform})
is nonzero. Equivalently, we have to find one set of labels
$J,L,b$ so that the maps $E_{ab}^{IJ} (K|L): V^K \mapsto V^L$
defined by
$$
    E_{ab}^{IJ}(K|L) \equiv C^b [ \bar J J |L] D^I( C[I \bar I|0]
    M^J (R')^{I \bar I}   C^a[I \bar I |K]^*) (R')^{ \bar J J }
    C[ \bar J J |0]^*
$$
do not vanish. A simple calculation shows
$$
  \t^L(\xi) E_{ab}^{IJ}(K|L) = E_{ab}^{IJ}(K|L)  \t^K(\xi) \ \ .
$$
It follows with the help of Schurs' lemma that $E_{ab}^{IJ} =
\d_{K,L} e_{ab}^{IJ}(K)$. The complex number $e_{ab}^{IJ}(K)$
determines whether the map $E_{ab}^{IJ}(K|L)$ vanishes or not.

In conclusion we find that the space $\De^I$ is equal to $\End
(V^I)$, if and only if for every pair $K,a$ such that $C^a[I
\bar I|K] \neq 0$ there is a pair $J,b$ such that $e_{ab}^{IJ}
(K)$ is nonzero. In fact, if the stated condition is satisfied,
it guarantees the existence of at least one element in $S^I(K,a)
\cap \De^I$. Considering that $\G$ acts on both $S^I(K,a)$ and
$\De^I$, we can conclude $S^I (K,a) \cap  \De^I = S^I (K,a)$ so
that $\De^I = \End (V^I)$ follows from $\bigoplus_{K,a} S^I (K,a) =
\End (V^I) $.

A long but straightforward calculation indeed shows that
$$ \sum_{b,J} e_{ab}^{IJ}(K) (e_{ab}^{IJ}(K))^* v_J d_J^2
   \neq 0\ \ . $$
This means that sufficiently many numbers $ e_{ab}^{IJ}(K)$ are
nonzero and establishes the irreducibility of $D^I$.

The whole set of representations $D^I$ of $\L$ maps the
loop algebra $\L$ to the algebra $\bigoplus_I \End (V^I)$.
Since $\L$ is spanned by the matrix elements of
monodromies $M^J \in \End(V^J) \o \L$ it has the dimension
$\sum_J \d_J^2$. The latter coincides with the dimension
of the space $\bigoplus_I \End (V^I)$ of representation
matrices so that the set $\{D^I\}$ is faithful.

\subsection{Gauss decomposition and loop-symmetry isomorphism}

The previous lemma has two prominent consequences which we would
like to mention in passing, even though they will not be used later.
They can be regarded as a refinement of the technique developed
by Faddeev, Reshetikhin and Takhtajan in \cite{FRT}.
First we can exploit Lemma \ref{irredlemma} to introduce
the quantum Gauss decomposition for the matrix generators of the
loop algebra. All representations $D^I$ of $\L$ are irreducible
so that there are elements $M_\pm^J \in \End(V^J) \o \L$ with
the property
$$  D^I(M_+^J) = (R')^{JI} \ \ , \ \
    D^I(M_-^J) = (R^{-1})^{JI} \ \ . $$
The elements $M^J_\pm$ are clearly invertible, i.e. there exist
elements $(M_{\pm}^J)^{-1}$ such that $ M^J_{\pm}(M^J_{\pm})^{-1}
= e^J =(M^J_{\pm})^{-1}M^J_{\pm}$. If we evaluate the product
$(\k_J)^{-1} M_+^J (M_-^J)^{-1} $ in the representation $D^I$,
it is found to agree with the representation matrix of $M^J$,
$$
  D^I(\k_J^{-1} M_+^J (M_-^J)^{-1})  = \k_J^{-1} (R'R)^{JI}
    = D^I(M^J)  \ \ .
$$
{}From faithfulness of the representation theory (Lemma
\ref{irredlemma}) we conclude that the matrix generators
satisfy $\k_J^{-1} M_+^J (M_-^J)^{-1} = M^J$.

\begin{coro} {\em (Gauss decomposition)} \label{Gauss}
In absence of truncation there exist elements $M^J_\pm
\in \End(V^J) \o \L$ such that
\be
     M^J= \k_J^{-1} M^J_{+}(M^J_{-})^{-1},
\ee
where $(M^J_{-})^{-1} \in \End(V^J) \o \L$ is the inverse
of $M_-^J$.
\end{coro}

A second consequence of Lemma \ref{irredlemma} was already
motivated by our discussion in the first subsection.

\begin{coro}  {\em (Loop-symmetry isomorphism)} \label{lsisom}
In absence of truncation the symmetry
Hopf-algebra $\G$ is isomorphic (as a Hopf-*-algebra) to
the algebra $\L$ generated by matrix elements of $M^I\in
\End(V^I) \o \L$ subject to the relations
$$
\M{1}{I} R^{IJ}\M{2}{J}=\sum C^a[IJ|K]^* M^K C^a[IJ|K] \ \ ,
$$
and supplied with the following co-product, co-unit,
antipode and *-operation
\ba
 \D (M^I)
   & = & \k_I^{-1} M^I_+ \tilde M^I (M_-^I)^{-1} \ \ , \nn \\[1mm]
 \e (M^I) & = & \k_I^{-1} e^I  \ \ , \nn \\[1mm]
 \S (M_+^I) = (M_+^I)^{-1} \ \ & , & \ \ \S(M_-^I) = (M_-^I)^{-1}
  \ \ , \nn \\[1mm]
 (M^I)^* & = & \s_\k ( (M^I_-)^{-1} (M^I)^{-1} M^I_- )  \ \ . \nn
\ea
Here $M_\pm^I$ are the Gauss components of $M^I$ and $\s_\k$ is
conjugation with $\k$ regarded as an element in $\L$. On the
right hand side of the first relation, the factors are elements
in $\End(V^I) \o \L \o \L$. $M^I_+, (M^I_-)^{-1}$ are
supposed to have trivial entry in the last component while
$\tilde M^I \equiv \sum m_\s \o e \o M_\s$ with $e$ being the unit
in $\L$ and $M^I = \sum m_\s \o M_\s \in \End(V^I) \o \L$.
\end{coro}

\noindent
{\sc Proof:} To understand the formulas for the action of
the co-product, co-unit etc., one has to evaluate them in the
representations $D^I$. On representation spaces they turn
into standard relations for the $R$-matrix, co-unit $\e$
etc. of the quantum symmetry $\G$. Faithfulness of the
representation theory allows to transport the results back
to the level of algebras $\L, \G$.

\section{The Moduli Algebra on a Sphere with Marked Points}
\setcounter{equation}{0}

We will now extend our analysis to some kind of
{\em multi-loop algebra}. In the terminology of Part 1, the
main actor of this section is the graph algebra $\L_{0,m}$
which is assigned to an $m$-punctured sphere. This leads us
to our simplest examples of moduli algebras. We will be able
to give a complete description of their representations
theory. The section also contains a proof of the first
{\em pinching theorem} that was used in \cite{AGS2} to
normalize the Chern-Simons measure $\omega_{CS}$.

\subsection{The graph algebra $\L_m$}

A full definition of the graph algebra $\L_m = \L_{0,m}$ was
given in Section 3.2. It is generated by
matrix elements of a family of monodromies
$M^I_\nu, \nu = 1, \dots, m$. In addition to the usual
functoriality, these monodromies are subject to the following
exchange relations,
\be
 (R^{-1})^{IJ} \M{1}{I}_{\nu}  R^{IJ} \M{2}{J}_{\mu}
  = \M{2}{J}_{\mu}   (R^{-1})^{IJ} \M{1}{I}_{\nu}  R^{IJ}
   \ \ \mbox{ if } \ \  \nu < \mu \ \ . \label{defSm2}
\ee
As one can see, the graph algebra $\L_m$ contains a bunch of
loop algebras. For fixed subscript $\nu$, the matrix elements
of $M_\nu^I$ generate a subalgebra
$\L(M_\nu)$ of $\L_m$ which is isomorphic to the loop algebra
$\L$. It was suggested in Section 3.4 to construct the elements
\be
 c^{I}_{\nu}= \k_I tr^I_{q}(M^{I}_{\nu})\ \ .
\ee
For fixed index $\nu$, they certainly generate a fusion algebra
and commute with every element in $\L(M_\nu)$. Moreover,
using the equation
$tr_q^I (  M^I_\nu ) =  tr^I_q \left[(R^{-1})^{IJ} M^I_\nu R^{IJ}
\right]$ we infer from eq. (\ref{defSm2}) that
$$
   c^I_\nu M^J_\mu =
    M^J_\mu c^I_\nu \ \ \mbox{ for } \ \    \nu < \mu
$$
and the same for $\mu > \nu$. This means that the {\em $c^I_\nu$
provide a large family of central elements } in the graph
algebra $\L_m$. Needless to say that we can pass to characteristic
projectors $\chi^K_\nu$ with the same transformation that was
used in Section 3.5. From these remarks we certainly expect
representations of the graph algebra $\L_m$ to be labeled by
tuples $ (I_1, \dots, I_m)$. This will be confirmed in the next
subsection.

\subsection{Representation theory of the graph algebra $\L_m$}

Let us turn to the representation theory of $\L_m$. Our strategy
is to construct the representations in tensor products
\be \label{Imn}
 \Im (I_1,\dots , I_m) = V^{I_1}\o \dots \o V^{I_m}
\ee
of the representations spaces of the underlying Hopf algebra $\G$.
To state the formulas we introduce the notation $D^{I_\nu}_\nu$
for the representation $D^{I_\nu}$ of the $\nu^{th}$ copy $\L(M_\nu)
\subset \L_m$ of the loop algebra on the $\nu^{th}$ factor in the
tensor product (\ref{Imn}), i.e. for every $X \in \L(M_\nu)$,
$$
D^{I_\nu}_\nu (X) = id_{I_1} \o \dots \o id_{I_{\nu-1}} \o
      D^{I_\nu} (X) \o id_{I_{\nu+1}} \o \dots \o id_{I_m}\ \ .
$$
Given a set of labels $I_1, \dots, I_m$, it is convenient to
use $\imath_\nu, \nu = 1, \dots m,$ as a shorthand for the
representations $\imath_1 = \e$ and
$$
 \imath_\nu (\xi)  =  (\tau^{I_1} \bo
   \dots \bo \t^{I_{\nu-1}})(\xi) \o id_{I_{\nu}} \o \dots
   \o id_{I_m} \ \ (\nu \geq 2)
$$
of the Hopf algebra $\G$ on the tensor product $V^{I_1} \o \dots
\o V^{I_m}$. $(\t^{I_1} \bo \dots \bo \t^{I_{\nu-1}})$ is the
ordinary tensor product of representations $\t^{I_1}$ through
$\t^{I_{\nu-1}}$ of $\G$. The representation $\imath=\imath_{m+1}$
coincides with the natural action of $\G$ in the space
$\Im (I_1, \dots, I_m)$. With these notations we are prepared
to define representations of $\L_m$.

\begin{theo} {\em (Representations of the algebra $\L_m$)}
\label{Smrep}
The algebra $\L_m$ has a series of representations
$D^{I_1,\dots, I_m}$ realized in the tensor product (\ref{Imn})
of representation spaces of the underlying quasi-triangular
Hopf algebra $\G$. In such a representation the generators of
$\L_m$ can be expressed as
$$
D^{I_1, \dots, I_m} (M_\nu^{J})  =  (\t^J \o \imath_\nu)(R')
   D_\nu^{I_\nu}(M_\nu^{J}) (\t^J \o \imath_\nu)((R')^{-1}) \ \ .
$$
These representations extend to the semi-direct
product $\S_m \equiv \L_m \sd \G$ with the help of
$$
D^{I_1, \dots, I_m} (\xi)    =   \imath (\xi) \ \
\mbox{ for all } \ \  \xi \in \G \ \ .
$$
The set of representations $\{ D^{I_1, \dots ,I_m}\}$ is
faithful on $\L_m$. In absence of truncation this implies that
the representations $D^{I_1, \dots, I_m}$ are irreducible.
\end{theo}

A more explicit formula for the action of $M^I_\nu$ involves the
elements $K_n \in \G^{\o (n+1)}$ with  $ K_1 = e \o e$ and
$$
 K_n = (K_{n-1}\o e) R'_{1n}  \ \ \mbox{ for } \ \ n \geq 2\ \ .
$$
Using the definition of $\imath_\nu$,  $D^{I_\nu}_\nu$ and
quasi-triangularity one may derive that
$$ D^{I_1, \dots,I_m} (M_\nu^{J})  =  (\k_J)^{-1}
   (K_{\nu} R'_{1(\nu+1)} R_{1(\nu+1)} K_{\nu}^{-1} \o
   e^{(m-\nu)})^{JI_1\dots I_m}    \ \ ,
$$
where $e^{(n)}$ is the unit element in $\G^{\o n}$.

We are certainly interested in the $*$-properties of the
representations $D^{I_1, \dots, I_m}$. They are described
by the following proposition.

\begin{prop} {\em (Scalar product for $\Im(I_1, \dots, I_m)$) }
\label{Imsp}
Suppose that $\t,\t'$ are two *-representations
of the Hopf algebra $\G$ on Hilbert spaces $V,V'$.
With $(.,.)$ denoting the standard scalar product
on $V \o V'$, the formula
\ba
      \langle v_1,v_2 \rangle  & \equiv &
      ( v_1, (\t \o \t')( R \eta) v_2 ) \ \ ,  \nn \\
      & \mbox{ with } & \eta \equiv \D(\k)
       ( \k^{-1} \o \k^{-1})            \nn
\ea
defines a (positive definite) scalar product for elements
$v_1, v_2  \in V \o V'$. The tensor product $(\t \bo \t' )$ is
a *-representation of $\G$ with respect to $\langle .,. \rangle$.
Iteration gives a scalar product $\langle .,. \rangle$ on
$\Im(V_1, \dots, V_m)$. With respect to this scalar product,
the $D^{I_1, \dots, I_m}$ are $*$-representations.
\end{prop}

Scalar products of this type have been proposed by Durhuus et al.
\cite{DJN}. They are  motivated by the fact that tensor products
$\t \bo \t'$ of $*$-representations are incompatible with the
$*$-operation when we use the scalar product $(.,.)$
on $V\o V'$ to give sense to $((\t \bo \t')(\xi))^*$. To avoid
confusion we would like to stress that the two different scalar
products $(.,.)$ and $\langle .,.\rangle$ furnish two different
notions of ``adjoint'' for linear maps $X : V \o V' \mapsto
V \o V'$. The adjoint with respect to $(.,.)$ is denoted by
$*$ and was used e.g. in $C^a[IJ|K]^*$.  For the new adjoint
provided by $\langle .,.\rangle $ we will occasionally employ
the symbol $\dagger $.

\noindent
{\sc Proof: } The proof of Proposition \ref{Imsp} and the main
statements in Theorem \ref{Smrep} is rather standard and we can
leave it as an exercise. Let us, however, give some arguments
which establish the irreducibility statement in Theorem
\ref{Smrep}.

For this purpose, we count the number of linear independent
elements in the image of the representation $D^{I_1, \dots ,I_m}$
and show that it is given by $\prod_\nu \d^2_{I_\nu} = \dim
(\Im (I_1, \dots, I_m))^2$. We do this by induction over $m$.
The case
$m=1$ has been dealt with in Lemma \ref{irredlemma}. So suppose
that the representation $D^{I_1,\dots , I_m}$ of the set $\L_m$
is irreducible. The representations
$D^{I_1, \dots, I_m, I_{m+1}}$ can be restricted to the
first $m$ monodromies. This restriction coincides with
$D^{I_1, \dots, I_m}$ acting on the first $m$ factors of
the tensor product $V^{I_1} \o \dots \o V^{I_{m+1}}$. So by
assumption we already found $\d_1= \prod_{\nu = 1}^m \d_{I_\nu}^2$
linear independent maps. It is worth noticing that all the maps
obtained by representing the first $m$ monodromies act
trivially on the last factor $V^{I_{m+1}}$. With this in
mind, let us turn to the $(m+1)^{th}$ monodromy $M_{m+1}^J$.
A short calculation shows that they are represented by
$$  D^{I_1, \dots, I_{m+1}}(M^J_{m+1}) =
   (\t^J \o \imath \o \t^{I_{m+1}}) \left( R^{-1}_{23} R'_{13}
         R_{13} R_{23} \right) \ \  $$
with $\imath = (\t^{I_1} \bo \t^{I_2} \bo \dots \bo \t^{I_m})$.
By Lemma \ref{irredlemma} such representation matrices account for
$\d_2 = \d^2_{I_{m+1}}$ linear independent maps. The latter act in a
very special way on the representation space. In fact, every
such map is of the form
$$ (\imath \o \t^{I_{m+1}}) \left( R^{-1} (e \o m) R \right) $$
with some $m \in \G$.

Let us finally check that products of basis elements in the two
discussed sets of representation matrices are linear independent.
Thereby we will establish the existence of $\d_1 \d_2 =
\prod_{\nu=1}^{m+1} (\d_{I_\nu})^2$ linear independent maps in the
image of the representation $D^{I_1, \dots, I_m}$ and hence the
irreducibility. So consider a basis $m^{(1)}_\s \o e^{I_{m+1}}, \s
= 1, \dots, \d_1$, in the space of representation matrices for  the
first $m$ monodromies and  similarly $\d_2$ linear independent
maps of the form $ (\imath \o \t^{I_{m+1}}) \left( R^{-1}
(e \o m^{(2)}_\a) R \right) , \a= 1, \dots, \d_2$ which come from
representing the monodromies $M^I_{m+1}$. With these two sets of
basis elements, linear relations between products are of the form
$$
\sum_{\s,\a} \lambda_{\s,\a} (m^{(1)}_\s \o e^{I_{m+1}})
(\imath \o \t^{I_{m+1}}) \left( R^{-1}
   (e \o m^{(2)}_\a) R \right) =0 \ \ .
$$
Let us introduce the expansion $R^{-1} = \sum s^1_\s \o s^2_\s$.
We multiply the equation with $(\imath \o \t^{I_{m+1}})(R^{-1}
(\S(s^1_\s) \o e)$ from the right and with
$(\imath \o \t^{I_{m+1}})(e \o s^2_\s)$ from the left and sum over
$\s$. This results in
$$
\sum_{\s,\a} \lambda_{\s,\a} (m^{(1)}_\s \o e^{I_{m+1}})
(\imath \o \t^{I_{m+1}}) (e \o m^{(2)}_\a ) =0 \ \ .
$$
Consequently, the complex coefficients $\lambda_{\s,\a}$ have to
vanish. This concludes the proof.

For further considerations we need a specific subalgebra in the
graph algebra $\L_2$.

\begin{defn} {\em (Diagonal subalgebra $\L_2^d $)} The diagonal
subalgebra $\L_2^d \subset \L_2$ is defined as
$$ \L_2^d = \L_2 \sum_K \chi_1^K \chi_2^{\bar K}\ \ . $$
\end{defn}

Irreducible representations of the algebra $\L_2^d$ are labeled
by the index $I$ which runs through the set of irreducible
representations of $\G$. These representations are realized in
the spaces $V^I \o V^{\bar I}$.

\subsection{The moduli algebra $\lM_m^{\{K_\nu\}}$}

Before we define the moduli algebra, we want to consider the
*-algebra $\A_m$ of invariants within the space $\L_m
$ (cp. Definition \ref{defAgm}).  Since $\A_m$ is a subalgebra
of $\L_m$, the
representations $D^{I_1, \dots, I_m}$ can be restricted and furnish
representations of $\A_m$ on the representation space $\Im(I_1,
\dots I_m)$. We denote the restricted representations by the same
letter $D^{I_1, \dots, I_m}$. As a representation of $\A_m$,
$D^{I_1, \dots, I_m }$
are reducible, or -- in other words -- $D^{I_1, \dots,I_m}(\A_m)$ has a
nontrivial commutant. Obviously, the latter contains all maps on
$\Im (I_1, \dots , I_m)$  which represent elements in $\G$ so that
invariant subspaces for the representation of $\A_m$ are determined
by the decomposition
\be \label{dec}
      \Im(I_1, \dots , I_m) =
      \bigoplus_J V^J\otimes W^J(I_1,\dots , I_m)\ \ .
\ee
Here the sum runs over irreducible representations of $\G$ and
$W^J(I_1,\dots , I_m)$ are multiplicity spaces. The decomposition
(\ref{dec}) is always possible for a semisimple symmetry algebra
$\G$. Because of Proposition \ref{Imsp}, the action of $\G$ on
$\Im (I_1, \dots, I_m)$ is consistent with the adjoint $\dagger$.
This implies  that formula (\ref{dec}) is compatible with
the scalar product $\langle .,. \rangle$ on $\Im(I_1, \dots, I_m)$
and consequently the multiplicity spaces $W^J(I_1, \dots I_m)$
come equipped with a canonical scalar product.

Let us argue that the restriction of $D^{I_1, \dots, I_m}$ to
$\A_m$ is irreducible on the spaces $W^J(I_1, \dots, I_m)$. We
noticed before that the set of representations $\{ D^{I_1,
\dots, I_m}\}$ of $\L_m$ is faithful.
When we restrict $D^{I_1, \dots, I_m}$ to the subalgebra
$\A_m \subset \L_m$, faithfulness survives. Now let $\d(D_m)$
denote the dimension of the space of representation
matrices of $\A_m$ on the direct sum
\be \label{Arepspace}
  \bigoplus_{I,I_1, \dots, I_m} W^I(I_1, \dots, I_m)\ \ .
\ee
Because of faithfulness of the representation theory of $\A_m$,
$\d(D_m)$ is equal to the dimension $\d(A_m)$ of the algebra
$\A_m$. The latter is easy to compute. Recall that elements
in $\A_m$ are obtained as linear combinations of
$$ tr^J_q\left( C_1 [I_1, \dots, I_m|J] M^{I_1}_1
  \dots M^{I_m}_m    C_2 [I_1, \dots, I_m|J]^*\right) $$
for arbitrary sets of labels $\{I_\nu, J\}$ and two intertwiners
$C_1,C_2 : V^{I_1} \o \dots \o V^{I_m} \mapsto V^J$. So we find
\ba
    \d (D_m) = \d (\A_m) & = &  \sum_{I_1, \dots, I_m,J}
   (\sum_{J_1, \dots, J_{m-2}} \sum N^{I_1 I_2}_{J_1} \dots
   N^{J_{m-2} I_m}_J)^2\ \   \nn \\[2mm]
   & = &
   \mbox{{\it dim\/}} \left( \bigoplus_{J,\{I_\nu\}}
   \End(W^J( I_1, \dots, I_m) \right) \ \ .  \nn
\ea
This result for $\d (D_m)$ shows that every map on the space
$W^J(I_1, \dots, I_m)$ appears in the  image of $\A_m$ under the
representation $D^{I_1, \dots, I_m}$. As a conclusion one should
keep in mind that each space {\em $W^J(I_1, \dots, I_m)$ carries
an irreducible *-representation of $\A_m$}. Since we are
interested in the moduli algebra, we do not want to formulate
this as a proposition.

Let us now pass from $\A_m$ to the moduli algebras. According to
our discussion in Section 3.4 a moduli algebra is prepared
by implementing $m+1$ additional (flatness-) relations. In the
finite dimensional case considered here, we may describe the
resulting object with the help of characteristic projectors
$\chi^I$. More explicitly, we will need the $m$ characteristic
projectors $\chi_\nu^{K_\nu}$ that were assigned to the loop
algebras $\L(M_\nu)$ in the first subsection. As we have seen there,
these characteristic projectors are central in $\L_m$ and hence
also in $\A_m \subset \L_m$. In addition to the $\chi^{K_\nu}_\nu$,
we will employ one more set of elements $\chi^K_0 \in \A_m$ that
is assigned to the total product $r=l_m \dots l_1$ of loops $l_m$
through $l_1$, i.e.
\ba
   \chi^K_0 & \equiv & \N d_K S_{K \bar I} c^I_0 = \N d_K \k_K
    S_{K \bar I} tr^I_q(M^I(r)) \ \label{chi0}\\[1mm]
    \mbox{ with } & & M^I(r) = \k_I^{m-1} M^I_m \dots M^I_1\ \ .
    \nn
\ea
By choice of the factor $\k_I^{m-1}$ in front of $M^I(r)$, the
corresponding elements $c^I_0$ generate a fusion algebra. While
the elements $\chi^K_0$ are not central in $\L_m$, it was shown
in \cite{AGS2} that they are {\em central} in $\A_m$. This will
be confirmed in the discussion below.

Now we are prepared to restate Definition \ref{defmodalg1} of
the {\em moduli algebra} $\lM_m^{\{K_\nu\}}$ in the finite
dimensional case.

\begin{defn} {\em (Moduli algebra)}
The moduli algebra $\lM^{\{K_\nu\}}_m$ of a sphere with
$m$ punctures marked by $K_\nu, \nu = 1, \dots, m,$ is the
$*$-algebra
\be
   \lM^{\{K_\nu\}}_m \equiv \chi^0_0 \prod_{\nu = 1}^m
    \chi^{K_\nu}_\nu     \A_m \ \ .
\ee
Here $\chi^{K_\nu}_\nu, \chi^0_0 \in \A_m$ are the central
projectors introduced in the text preceding this definition.
\end{defn}

To determine the representation theory of the moduli algebra we
have to evaluate the characteristic projectors within the
representation $D^{I_1, \dots, I_m}$. Let us start with
$\chi^{K_\nu}_\nu ,\nu = 1, \dots, m$. It is easy to see that
\be              \label{charInu}
     D^{I_1, \dots, I_m} (\chi^{K_\nu}_\nu) =
    \d_{I_\nu, K_\nu} e^{I_1} \o \dots \o e^{I_m} \ \ .
\ee
The reasoning is similar to the derivation of eq. (\ref{charJ}) and
uses that traces are invariant under conjugation with $(\t^J \o
\imath_\nu)(R')$ . The formula (\ref{charInu}) implies that the
*-representation $D^{I_1, \dots, I_m}$ of the moduli algebra
$\lM^{\{K_\nu\}}_m$ on the spaces $W^J (I_1, \dots , I_m)$ is
nonzero, if and only if $I_\nu = K_\nu$ for all $\nu = 1,
\dots, m$.

Evaluation of $\chi^0_0$ in the representations $D^{K_1, \dots,K_m}$
is more difficult. From formula (\ref{chi0}) we see that $\chi^0_0
= \N^2 d_L c^L_0$ with $c^L_0 = \k_I tr^L_q(M^L(r))$. When
$c^L_0$ is evaluated with $D^{K_1, \dots K_m}$ using our explicit
formulas after Theorem \ref{Smrep} we encounter the following
expression in the argument or the $q$-trace $tr^L_q$
\ba
   & &
       K_m R'_{1(m+1)} R_{1(m+1)} K_m^{-1}(K_{m-1}
       R'_{1m} R'_{1m} K_{m-1}^{-1} \o e)
      \dots R_{12}  \nn \\[1mm]
    &=&
       K_m R'_{1(m+1)} R_{1(m+1)} (R'_{1m})^{-1}
       R'_{1m} R'_{1m} (R'_{1(m-1)})^{-1}
      \dots R_{12} \nn \\[1mm]
    &=&  (id \o \D^{(m-1)})(R'R)       \nn
\ea
In the derivation the definition of $K_\nu \in \G^{\o_\nu}$
was inserted. The result for $D^{K_1, \dots, K_m}$ evaluated
on $c^L_0$ is
$$
    D^{K_1, \dots, K_m} (c^L_0)
    = tr^L_q ((\t^L \o \imath_m)(R'R))\ \ .
$$
On the multiplicity spaces $W^J(K_1, \dots, K_m)$ this gets
represented by
$$
    D^{K_1, \dots, K_m} (c^L_0)|_{W^J(K_1, \dots, K_m)}
     =  tr^L_q ((\t^L \o \t^J)(R'R)) \ \ .
$$
The same argument that resulted in the formula (\ref{charJ}) can
finally be employed to conclude
$$
   D^{K_1, \dots, K_m} (\chi^0_0) |_{W^J(K_1, \dots, K_m)}
   =    \d_{J,0} id_{W^0(K_1, \dots, K_m)} \ \ .
$$
So we end up with only one nonzero representation of the
moduli algebra $\lM_m^{\{K_\nu\}}$ on the space $W^0(K_1, \dots,
K_m )$. Because of faithfulness of the representations theory,
other representations cannot exist. We may summarize these
findings in the following theorem.

\begin{theo} {\em (Representations of the moduli algebra; genus 0)}
\label{repmod1}
For any set  $K_1,.$ $\dots, K_m $ labeling $m$ points on a
Riemann surface of genus $0$, there is a unique irreducible
*-representation of the corresponding moduli algebra $\lM^{
\{K_\nu\}}_{m}$ on the space $W^0(K_1, ...$ $., K_m)$ (defined
through the decomposition (\ref{dec})). This representation can
be obtained explicitly by restricting the representation
$D^{K_1, \dots, K_m}$ of $\L_m$ to the moduli algebra
$\lM^{\{K_\nu\}}_m$.
\end{theo}

\subsection{The first pinching theorem}

The different algebras $\lM_m^{\{K_\nu\}}$ are related with each
other. In fact one can construct various inclusions which are
parametrized by the choice of a circle on the punctured surface.
Geometrically, the inclusions corresponds to a pinching of
the surface. We will explain this only for one particular
pinching-circle, but the idea is more general.

To state and prove the {\em first pinching theorem} we need to
introduce some new notations. Suppose that ${\cal X}$ is a
subalgebra of an algebra ${\cal Y}$. The (relative) commutant
of ${\cal X} \subset {\cal Y}$ will be denoted by ${\cal C}
({\cal X}, {\cal Y})$.

Consider the moduli algebra $\lM_ m^{ K_1, \dots , K_m}$
corresponding  to a sphere with $m$ marked points. Pick up a cycle
$l=l_{m} l_{m-1}$ (i.e. the product of the two elementary loops
$l_{m}$ and $l_{m-1}$) and construct the fusion algebra ${\cal V}(l)$
assigned to $l$, i.e. the algebra generated by
$$ c^I(l) = tr^I_q( M^I_{m} M^I_{m-1}) \ \ . $$
We wish to investigate the commutant of  ${\cal V}(l)$ in
$ \lM_m^{K_1,  \dots , K_m}$. The result is given
by the following theorem.

\begin{theo} {\em (First pinching theorem)}
\label{1pintch}
The commutant of ${\cal V}(l)$ in  $\lM_{ m}^{K_1, \dots , K_m}$
splits into the direct sum of products of moduli algebras corresponding to
$m-1$ and $3$ marked points
\be \label{ThD}
{\cal C}({\cal V}(l), \lM_{ m}^{ K_1, \dots , K_m}) \cong
\bigoplus_K \lM_{m-1}^{K_1, \dots, K_{m-2}, K}\o \lM_{ 3}
^{\bar K, K_{m-1}, K_m}.
\ee
Here the sum runs over all classes of irreducible representations
of the symmetry Hopf algebra $\G$.
\end{theo}

In Topological Field Theory the evaluation of the commutant should be
interpreted as a fusion of two marked points into one. One
can imagine that we create a long neck which separates these two points
from the rest of the surface. When we cut the neck, the surface splits
into two pieces. The ``main part'' carries the rest of marked points
and a new one created by the cut. The other piece has only three
punctures, two of them are those that we wish to fuse and the new
one appears because of the cut. Iteration of this procedure
results in a product of 3-punctured spheres.

\noindent
{\sc Proof:} The proof consists of two parts. First we construct
an embedding $\Phi$ of the algebra on the right hand side of relation
(\ref{ThD}) into the moduli algebra $\lM_{m}^{K_1, \dots, K_m}$.
Then we show that the commutant of the image of $\Phi$ is isomorphic
to the fusion algebra  ${\cal V}(l)$. By semisimplicity, this is
equivalent to the statement in the theorem.

In preparation let us observe an isomorphism between moduli
algebras $\lM_p^{\{I_\nu\}}$ and the algebras
\be  \label{algmm1}
    \chi^{I_p}_0 \prod_{\nu = 1}^{p-1}
    \chi^{I_\nu}_\nu    \A_{p-1}
\ee
which are obtained from the graph algebras $\S_{p-1}$ with the
help of characteristic projectors (for notations cp. last
subsection). In fact, our discussion of the representation
theory of moduli algebras can be applied to show that algebras
(\ref{algmm1}) possess a unique irreducible representation on
the multiplicity space $ W^{I_p} (I_1, \dots, I_{p-1})$. For
dimensional reasons, this
space is isomorphic to the representation space $W^0(I_1,
\dots, I_p)$ of $\lM_p^{\{ I_\nu\}}$ so that an isomorphism
of algebras follows from Theorem \ref{repmod1}. Alternatively,
this isomorphism may be obtained from the results in \cite{AGS2} on
the independence of moduli algebras from the choice of the graph.

Let us consider the graph algebras $\L_{m-2}$ and $\L_2$.
Generators for $\L_{m-2}$ will be denoted by $\tilde M^I_\nu,
\nu = 1, \dots, m-2,$ and for $\L_2$ we use $\hat M^I_i, i =
1,2$. Analogous conventions apply for other elements. In
particular we will need our standard characteristic
projectors $\tilde \chi^{K_\nu}_\nu, \nu= 0, \dots, m-2,$
in $\A_{m-2} \subset \L_{m-2}$ as well as $\hat
\chi^{K_{m-2+i}}_i, i = 0,1,2,$ in $\A_2 \subset \L_2$.

There is an obvious embedding  $\phi: \L_{m-2} \o \L_2
\mapsto \L_m$ defined by
$$
   \phi( \tilde M^I_\nu) =  M^I_\nu \ \ \ \ ,\ \ \ \
   \phi( \hat M^I_i) =  M^I_{m-2+i}
$$
for all $\nu = 1, \dots, m-2$ and $i = 1,2$.
The definition of $\phi$ implies that
$$
    \phi(\tilde \chi^{K_\nu}_\nu) = \chi^{K_\nu}_\nu \ \ \ \ ,
    \ \ \ \ \phi(\hat \chi^{K_{m-2+i}} _i ) =
     \chi^{K_{m-2+i}}_{m-2+i}
$$
where $\nu$ ranges over $1, \dots, m-2$ and $i = 1,2$.
Consequently, the map $\phi$ induces an embedding
\be \label{embedd1}
   \phi:   \bigoplus_K \lM_{m-1}^{K_1, \dots, K_{m-2}, K}\o
    \lM_{ 3}^{\bar K, K_{m-1}, K_m} \mapsto
    \A_m \prod_{\nu = 1}^m \chi^{K_\nu}_\nu \ \ .
\ee
In writing this we also used the isomorphism observed in the
second paragraph of this proof.

The result (\ref{embedd1}) is not yet strong enough. In fact we
have embedded the direct sum of moduli algebras on the left hand
side into an algebra that is much larger than the moduli algebra
$\lM_m^{K_1, \dots, K_m}$. Notice that the latter contains a
factor $\chi^0_0$ in its definition which does not appear
on the right hand side of (\ref{embedd1}). So it remains to
understand why the full matrix algebras
 $$
   M(K) \equiv \lM_{m-1}^{K_1, \dots , K_{m-2}, K}\o
    \lM_{ 3}^{\bar K, K_{m-1}, K_m}
$$
are embedded into the direct summand $\lM_m^{K_1, \dots, K_m}
\subset \A_m \prod_{\nu = 1}^m \chi^{K_\nu}_\nu$. The projection
into the moduli algebra $\MC_m^{K_1, \dots, K_m}$ is furnished
by the element $\chi^0_0$. It suffices to show that the
unit element in $\phi(M(K))$ is projected to a nontrivial
element of the moduli algebra. If we denote the
product of loops $l_{m-2} \dots l_1$ by $\tilde l$ and use
$\hat l = l_m l_{m-1}$, we can write the unit element in
$\phi(M(K))$ as $ \chi^K(\tilde l) \chi^{\bar K} (\hat l)$.
So we wish to demonstrate that $ E \equiv
\chi^0_0 \chi^K(\tilde l) \chi^{\bar K} (\hat l) \neq 0$.
We do this by showing that $\omega(E) \neq 0$. With the help
of  Lemma 3 in  \cite{AGS2} one obtains indeed
\ba
   \omega( E) & = & \omega (\chi^0_0 \chi^K(\tilde l) \chi^{\bar K}
    (\hat l))  \nn \\[2mm]
      & = & \sum_L \omega (\chi^0_0 \chi_0^L)
       =  \omega((\chi^0_0)^* \chi^0_0) > 0 \ \ .    \nn
\ea
We used that $\chi^0_0 \chi^L_0 = \delta_{L,0} \chi_0^0 =
(\chi^0_0)^* \chi^0_0$ and positivity of the functional
$\omega$ (cp. Theorem \ref{positivity}). According to our prior
remarks we can now conclude that $\Phi \equiv \chi^0_0
\phi $ defines the desired embedding.

Elements in the fusion algebra ${\cal V}(l)$ over the circle
$l = \hat l$ obviously commute with the image of $\Phi $.
We will show now that the moduli algebra $\lM_m^{\{K_\nu\}}$
contains no other elements with this property. Once more,
a comparison of dimensions is useful. Indeed, square roots
of the dimensions of the matrix blocks $\lM_{m-1}^{K_1,
\dots, K_{m-2}, K}\o \lM_{3}^{\bar K, K_{m-1}, K_m}$
add up to the square root of the dimension of the moduli
algebra $\lM_m^{\{K_\nu\}}$. This implies that every
block is embedded with multiplicity one into the moduli
algebra and hence the commutant of the image of $\Phi$
cannot contain more than the algebra of its minimal central
projectors. The latter coincides with the fusion algebra
${\cal V}(l)$. This concludes the proof of the first
pinching theorem.

\section{The Handle Algebra and its Representations}
\setcounter{equation}{0}

Before we get to discuss surfaces of higher genus, we have to
introduce one more elementary building block: the {\em $AB$- (or
handle-) algebra}. It will turn out to be associated to
a handle of the surface, much as a loop algebra came with
every puncture. The $AB$-algebra is a graph algebra -- namely
the algebra $\L_{1,0}$ -- assigned to a pair of links winding
around a handle. These links correspond to the $a-$ and $b-$
cycles (hence the name $AB$-algebra). The main topic is again
the representation theory for the $AB$-algebra in the
second subsection. At the end of the section, we discuss
some properties of the algebra which will be used in
the next two sections.

\subsection{The $AB$-algebra}

Following our tradition, let us recall that the $AB$-algebra
$ \T = \L_{1,0}$ is generated by
matrix elements of two monodromies $A^I = A_1^I, B^I =
B^I_1 \in \End(V^I) \o \T$. As usual, the monodromies satisfy
functoriality. A more characteristic feature are the exchange
relations between $A^I$ and $B^J$,
\be
 (R^{-1})^{IJ} \Ae{1}{I} R^{IJ} \Be{2}{J}   =
  \Be{2}{J} (R')^{IJ} \Ae{1}{I} R^{IJ} \ \ .
  \label{defAB3}
\ee
For a complete description of $\T = \L_{1,0}$ including covariance
properties and the $*$-operation on $\S_{1,0} = \T \sd \G$ we refer
the reader to Sections 3.2, 3.3.

We recognize two copies of the loop algebra inside the $AB$-algebra
which are associated with the monodromies $A^I$ and $B^J$. As
before one may construct elements $c^I$ from these monodromies.
Unlike in the previous section, they turn out not to be central
in the $AB$-algebra. In fact we will see shortly that $\T$
has a trivial center.

The algebra $\T$ is isomorphic to some well-known object.
Namely, there are many ways to identify it with the quantized algebra
of functions on the Heisenberg double corresponding to the symmetry
Hopf algebra $\G$. The Heisenberg double is a Poisson-Lie counterpart of
the cotangent bundle to a Lie group. The corresponding quantized
algebra of functions is a generalization of the algebra of finite
order differential operators on a Lie group. The isomorphism between
the algebra $\T$ and the quantized algebra of functions on the
Heisenberg double is not canonical. One of the reasons for that is
a wide group of automorphisms of the $AB$-algebra $\T$. They will
be discussed in the last subsection.

\subsection{Representation theory of the $AB$-algebra}

The representation theory of the $AB$-algebra is closely related to the
representation theory of the underlying Hopf algebra as in the case
of the loop algebra. Here we will find precisely one representation
which acts in the space of the regular representation of the Hopf
algebra $\G$.
\be\label{reg}
 \Re =\bigoplus_{I} \End( V^I )     \ \ .
\ee
Here $\End(V^I)$ is regarded as a complex vector space. We will
introduce a scalar product on $\Re$ below.

\begin{theo} {\em (Representation $\pi$ of the $AB$-algebra)}
\label{ABrep}
The $AB$-algebra $\T$ has one representation $\pi$ realized on
the representation space $\Re$. In this representation there
exists a cyclic vector $\vac$ such that
\be
      \pi(B^I) \vac = \vac (\k_I)^{-1}  \label{Bxivac}
\ee
This property determines the representation $\pi$ uniquely. $\pi$
can be extended to a representation of the semi-direct product
$\T \sd \G$ such that the vector $\vac$ is invariant under the
action of elements $\xi \in \G$, i.e.
\be   \pi(\xi) \vac = \vac \e(\xi) \ \ . \label{xivac} \ee
The representation $\pi$ is faithful on $\T$, in absence
of truncation it is irreducible.
\end{theo}

A representation space of $\T$ with properties as specified
in the theorem is generated from the ``ground state''
$\vac $ by application of operators $A^I_{cd}$. So it is
obviously isomorphic to $\Re$. Moreover, the exchange relations
of monodromies $B^J$ and $A^I$ and the transformation properties
of $A^I$ together with relations (\ref{Bxivac}) determine the
action of matrix elements in $B^J$ and of elements $\xi \in \G$
on arbitrary states in $\Re$.
\ba
   \pi (\Be{2}{J})  (R')^{IJ}\pi(\Ae{1}{I}) \vac  & = &
    ((R'R)^{-1})^{IJ} (R')^{IJ} \pi (\Ae{1}{I}) \vac \k_J^{-1}
    \nn \\[1mm]
   \pi (\mu^I(\xi)) \pi (A^I) \vac & = & \pi (A^I) \vac \e(\xi)
    \nn
\ea
with $\mu^I(\xi) = (\t^I \o id)(\D(\xi)) \in \End(V^I) \o \G$.
The action of $A^J$ can be derived by using functoriality of
monodromies $A^K$
$$  \pi (\Ae{1}{J}) R^{JI} \pi(\Ae{2}{I}) \vac
    = \sum C^a[JI|K]^* \pi(A^K) \vac C^a[JI|K] \ \ . $$
Faithfulness and irreducibility of $\pi$ will be demonstrated
at the end of the next subsection.

{}From the structure of $\pi$ we obtain an important consequence in
connection with Theorem \ref{positivity} in the background review.
Let $\L (A) \subset \T$ be the loop algebra generated by matrix
elements of $A^I$. Analogous to eq. (\ref{integ}) we define a
functional $\omega_A : \L (A) \sd \G \mapsto {\bf C}$ by
\be
     \omega_A ( A^I \xi) = \delta_{I,0} \e(\xi)\ \ .
     \label{omega}
\ee
The positivity result in Theorem \ref{positivity} furnishes the
following proposition.

\begin{prop} {\em (Scalar product for $\Re$)}
Suppose that $a_1, a_2$ are two elements in $\Re$ and
that they are obtained from the ground state $\vac$ by application
of $A_1, A_2 \in \L(A)$, i.e. $\pi(A_i)\vac = a_i$ for $i = 1,2$. Then
the formula
\be          \label{Resp}
   \langle a_1 , a_2 \rangle \equiv \omega_A((A_1)^* A_2 ),
\ee
with $\omega_A $ given through eq. (\ref{omega}), defines a (positive
definite) scalar product on $\Re$, iff all quantum dimensions $d_I$
are positive. With respect to this scalar product, the map $\pi:
\S_{1,0} = \T \sd \G \mapsto \End(\Re )$ defined in Theorem
\ref{ABrep} is a  $*$-representation.
\end{prop}

\subsection{Embedding the loop algebra into the $AB$-algebra}

The $AB$-algebra is closely related to the loop algebra
considered above. Obviously, $M^I \mapsto A^I$ as well as
$M^I \mapsto B^I$ define embeddings of the loop algebra into
$AB$-algebra. Here we give a more sophisticated embedding
which will be of special importance for us. None of
the following results has conceptual importance. Nevertheless
the subsection serves a twofold purpose: it provides some
more background material needed in the proofs of the next
section and prepares for Section 9 as well. The calculations
done here are typical for the discussion of the mapping class
group in our approach.

\begin{lemma} {\em (Automorphisms of the $AB$-algebra)}
The maps $i,j$ defined by
\ba \label{Aut}
i (A^I) =  \k_I^{-1} B^I A^I & , & i(B^J) = B^J\ \ ;  \nn \\
j (B^I) =  \k_I^{-1} B^I A^I & , & j(A^J) = A^J \ \  \nn
\ea
extend to *-automorphisms of the $AB$-algebra $\T$.
\end{lemma}

\noindent
{\sc Proof:} It suffices to give the proof for $i$.
To begin with, let us determine the multiplication rules of the
product $\k_I^{-1} B^I A^I$.
\ba
\k_I^{-1} \Be{1}{I} \Ae{1}{I} R^{IJ}  \k_J^{-1} \Be{2}{J} \Ae{2}{J}
  &=& (\k_I\k_J)^{-1} \Be{1}{I} R^{IJ} \Be{2}{J} (R')^{IJ}
          \Ae{1}{I} R^{IJ} \Ae{2}{J}\nn \\[1mm]
  &=& (\k_I \k_J)^{-1}  \sum C^a[IJ|K]^*  B^K C^a[IJ|K]
         \cdot \nn  \\
     & &  \hspace*{1cm} \cdot \
                (R')^{IJ}   C^b[IJ|L]^* A^L C^b[IJ|L]   \nn \\[1mm]
  &=&  \sum  \k_K^{-1}
        C^a[IJ|K]^*  B^K A^K C^b[IJ|K] \ \ .  \nn
\ea
This coincides with the multiplication rules for $A^I$. The
exchange relations for $\k_I^{-1} B^I A^I$ with $B^J$ are
\ba
\k_I^{-1} (R^{-1})^{IJ} \Be{1}{I} \Ae{1}{I}
  R^{IJ} \Be{2}{J}
& = & \k_I^{-1}
 (R^{-1})^{IJ} \Be{1}{I} R^{IJ} \Be{2}{J}
 (R')^{IJ} \Ae{1}{I} R^{IJ}                \nn \\[1mm]
& = & \k_I^{-1}
  \Be{2}{J} (R')^{IJ} \Be{1}{I}
     \Ae{1}{I} R^{IJ} \ \ . \nn
\ea
In the process of this calculation we inserted a relation of the
type (\ref{Mex}) which follows from the operator products of $B^I$.
Exchange relations for $(\k_J)^{-1} B^J A^J$ and $A^I$ are derived
in the same way. They coincide with the exchange relations of $A^I$
and $B^J$. Consistency with the covariance relations is obvious.
So it remains to discuss the properties of $i$ with respect to the
$*$-operation.
\ba
i ( (A^I)^* ) & = &  i(\s_\k(R^I (A^I)^{-1} (R^{-1})^I)) \nn\\[1mm]
 & = &  \s_\k (\k_I R^I (A^I)^{-1} (B^I)^{-1}
        (R^{-1})^I) =  (i ( A^I))^* \ \ .   \nn
\ea
This proves the lemma.

Now let us come to the main theme of this subsection, namely
the embedding of the loop algebra into the $AB$-algebra.

\begin{lemma} Let elements $(A^I)^{-1}, (B^I)^{-1} \in \End(V^I)
\o \T$ be defined through the equations $(A^I)^{-1} A^I = e^I
= A^I (A^I)^{-1}$ and $(B^I)^{-1} B^I = e^I = B^I (B^I)^{-1}$
as before. Then the map $\hvrho$,
\be \label{ABAB}
    \hvrho(M^I)  = \k^{3}_I B^I (A^I)^{-1} (B^I)^{-1} A^I
\ee
extends to an embedding of the loop algebra $\L$ into the
$AB$-algebra $\T$. The embedding extends to semidirect
products with the symmetry algebra $\G$ and respects the
action of $*$.
\end{lemma}

\noindent
{\sc Proof:} We have to determine the multiplication rules
for the product on the r.h.s of equation (\ref{ABAB}). In the
proof one first inserts the exchange relation for A and B in the
middle. Then we can use that $(B^J)^{-1} (A^J)^{-1}$ has the
same exchange relations with $A^I\ [B^I] $ as $(B^J)^{-1} \
[(A^J)^{-1}]$ (preceding lemma).
\ba
& &
\k^{3}_I \Be{1}{I} (\Ae{1}{I})^{-1} (\Be{1}{I})^{-1}
\Ae{1}{I} R^{IJ} \k^{3}_J \Be{2}{J} (\Ae{2}{J})^{-1}
(\Be{2}{J})^{-1} \Ae{2}{J}
 \nn \\[1mm]
&=&  (\k_I \k_J)^{3}
\Be{1}{I} (\Ae{1}{I})^{-1} (\Be{1}{I})^{-1} R^{IJ}
\Be{2}{J} (R')^{IJ} \Ae{1}{I} R^{IJ} (\Ae{2}{J})^{-1}
 (\Be{2}{J})^{-1} \Ae{2}{J}
\nn \\[1mm]
&=&  (\k_I \k_J)^{3}
\Be{1}{I} R^{IJ} \Be{2}{J} (R^{-1})^{IJ}(\Ae{1}{I})^{-1}
(\Be{1}{I})^{-1} ((R')^{-1})^{IJ} \cdot \nn \\[1mm]
& & \hspace*{2cm} \cdot (\Ae{2}{J})^{-1} (\Be{2}{J})^{-1}
(R^{-1})^{IJ} \Ae{1}{I} R^{IJ} \Ae{2}{J} \ \ . \nn
\ea
Finally we use the multiplication rules for $A^I, B^I$ and
$\k_I (A^I)^{-1} (B^I)^{-1}$. With the normalization (\ref{pos})
this gives
\ba
&=&  (\k_I \k_J)^2
   \sum C^a[IJ|K]^* B^K \frac{\k_K}{\k_I \k_J} \k_K  (A^K)^{-1}
   (B^K)^{-1}  \frac{\k_K}{\k_I \k_J} A^K C^a[IJ|K]   \nn  \\[1mm]
&=&\sum C^a[IJ|K]^* (\k_K)^3 B^K   (A^K)^{-1} (B^K)^{-1}
       A^K C^a[IJ|K]  \ \ . \nn
\ea
To show that $\hvrho :\L \mapsto \T$ respects the action of $*$ is
straightforward. It uses the fact that $((A^I)^{-1})^* = \s_\k(R^I
A^I (R^{-1})^{I})$ and $((B^I)^{-1})^* = \s_\k(R^I B^I (R^{-1})^I)$.

\noindent
{\bf Remark:} Let us remark that the embedding given by (\ref{ABAB})
is invariant with respect to automorphisms described in equation
(\ref{Aut}), i.e.
\be
i(\hvrho(M)) = \hvrho(M) \ \ \ , \ \ \ j(\hvrho(M)) = \hvrho(M)
\ee
for all $M \in \L $. The proof is straightforward. Indeed,
\ba  i(\k^3_I B^I (A^I)^{-1} (B^I)^{-1} A^I)
   & = & \k^3_I B^I (\k^I (A^I)^{-1} (B^I)^{-1}) (B^I)^{-1}
     (k_I^{-1} B^I A^I)  \nn \\[1mm]
   & = &   \k^3_I B^I (A^I)^{-1} (B^I)^{-1} A^I \nn
\ea
and similarly for $j$ instead of $i$.

For later application we wish to evaluate the image of $\hvrho$
in the representation $\pi$.

\begin{lemma}  \label{Minpi}
The quantum monodromies $\k^3_I B^I (A^I)^{-1}
(B^I)^{-1} A^I$ which appear as the image of $M^I$ under
the embedding $\hvrho: \L \mapsto \T$, are represented by
\be  \label{pirhoM}
  \pi(\k^3_I B^{I} (A^{I})^{-1} (B^{I})^{-1} A^{I}) =
  \pi (\k_I^{-1} (R'R)^I ) \ \ .
\ee
Here $\pi$ is regarded as a representation of the
semi-direct product $ \T \sd \G$ and
$(R'R)^I \equiv (\t^I \o id )(R'R) \in \End(V^I) \o \G$.
\end{lemma}

\noindent
{\sc Proof:} The computation is done  in several steps
using the eq. (\ref{Bxivac}) for $\pi$. To begin with, let
us show that
$$ \pi((B^I)^{-1} A^I) \vac = \pi(A^I )\vac
   \k_I^{-3}\ \ .$$
In fact, when a variant of relation (\ref{defAB3}) is applied to
the ground state $\vac$ and eq. (\ref{Bxivac}) is inserted
one obtains
\be \pi((\Be{2}{J})^{-1} (R^{-1})^{IJ} \Ae{1}{I}) \vac
   = (R')^{IJ} \pi(\Ae{1}{I}) \vac \k_I \ \ .
\label{form1}\ee
Suppose that we expand the inverse of the $R$-matrix
according to $R^{-1} = \sum s^1_\s \o s^2_\s$. Now
multiply the above equation from the left with
$\t^I(\S(s^1_\s))$ and from the right by $\t^J(s^2_\s)$
and sum over $\s$. Then taking the product of the two
components (for $I = J $) results in
$$ \pi((B^I)^{-1} A^I) \vac = \t^{I} (r^1_\t \s^2_\s \S(s^1_\s)
     r^2_\t) \pi(A^I) \vac  \k_I\ \ .$$
Finally, a short computation reveals that $ r^1_\t s^2_\s
\S(s^1_\s) r^2_\t  = u^{-1} \S(u^{-1}) = v^{-2} = \k^{-4}$.
This gives the anticipated formula. As a corollary we
note that
$$  \pi(\hvrho(M^I))\vac = \vac \k_I^{-1} \ \ . $$
To proceed with the evaluation of $\hvrho(M^I)$ on more general
states, one needs the exchange relation
$$ \hvrho(\M{1}{I}) R^{IJ} \Ae{2}{J} (R^{-1})^{IJ} =
   ({R'}^{-1})^{IJ} \Ae{2}{J} (R')^{IJ} \hvrho(\M{1}{I})\ \ .
$$
Its derivation is left as an exercise. From this we may conclude
$$ \pi(\hvrho(\M{1}{I}) R^{IJ} \Ae{2}{J} )\vac  =
   ({R'}^{-1})^{IJ} \pi(\Ae{2}{J}) \vac \k_I^{-1} (R'R)^{IJ} \ \ .
$$
We keep the answer in mind while we notice another formula which
can be derived from eq. (\ref{xivac}) and quasi-triangularity of
the $R$-element.
$$
    \pi(\k_I^{-1} (R'R)^I R^{IJ} \Ae{2}{J} )\vac =
    ({R'}^{-1})^{IJ} \pi(\Ae{2}{J})
    \vac \k_I^{-1} (R'R)^{IJ}\ \ .
$$

At this point we are prepared to finally prove the irreducibility
of $\pi$ in Theorem \ref{ABrep}.

\noindent
{\sc Proof of irreducibility of $\pi$:} We will show directly
that the representation $\pi$ is irreducible. Faithfulness follows
from a counting argument. The idea of the proof is this: first we
decompose the representation space $\Re$ into certain subspaces
and show that a subalgebra of $\T$ acts irreducibly on these
subspaces. Then we employ this result to show that every vector
in $\Re$ is cyclic. Because of the first part of the proof, it
actually suffices to establish cyclicity for one vector in each
of the considered subspaces.

To begin with, let us note that the monodromies $B^I$ can be
used to project onto subspaces of $\Re$ which are isomorphic
to $\End(V^J)$. The elements
\ba \chi^J_B & \equiv & \N d_J S_{J \bar I} \k_I tr^I_q(B^I)
     \nn \\[1mm]
 \mbox{ satisfy } & & \pi(\chi^J_B A^I) \vac = \pi(A^I) \vac
    \d_{I,J}\ \ .        \label{chiBA}
\ea
This follows easily with the help of a formula similar to eq.
(\ref{form1}) and the argument after the proof of Theorem
\ref{thMrep}. The subspaces $\pi(\chi^J_B) \Re$ carry
irreducible representations of a certain subalgebra in $\T$.

\begin{lemma}  \label{LDembed}
The map $\upsilon: \L_2 \mapsto \T$ defined by
$$ \upsilon( M_1^I) = \k_I^2 (A^I)^{-1}  (B^I)^{-1} A^I \ \ \ ,
\ \ \ \upsilon (M^I_2) = B^I $$
restricts to an embedding of the diagonal subalgebra $\L_2^d
\subset \L_2$ into the handle algebra. The image of $\L_2^d$
under $\upsilon$ is represented irreducibly on the subspaces
$\pi( \chi^J_B)\Re \subset \Re$.
\end{lemma}

{\sc Proof:} It is straightforward to prove that $\upsilon$
gives a homomorphism. In order to understand that $\upsilon$
furnishes an embedding of $\L_2^d$ into $\T$ it is sufficient
study the representations of the image $\upsilon (\L_2^d)$.

We infer from
\ba
  \pi ((\Be{2}{J})^{-1} (R')^{IJ} \Ae{1}{I}) \vac & = & \k_I^{-1}
  (R'R)^{IJ} (R')^{IJ} \pi(\Ae{1}{I}) \vac \nn \\[1mm]
   \pi \left((\k_J^2 (\Ae{2}{J})^{-1}(\Be{2}{J})^{-1}\Ae{2}{J})
   (R')^{IJ} \Ae{1}{I}\right) \vac & = &
    (R')^{IJ} \pi(\Ae{1}{I}) \vac \k_J^{-1} (RR')^{IJ}\ \ . \nn
\ea
that the action of monodromies $B^J$ and $\k_J^2 (A^J)^{-1}
(B^J)^{-1} A^J$ restricts to the subspace $\End(V^I)$
on which $\chi^I_B$ projects. To gain a better understanding of
this action, let us map the spaces $ \End(V^I) \subset
\Re$ to $ V^I \o V^{\bar I}$ by means of
$$
     w^I \equiv \k_I^{-3} C[I \bar I|0] (R')^{I \bar I}
              a^I R^{I \bar I} \in V^I \o V^{\bar I}
$$
for every $a^I \in \End(V^I) \subset \Re$. It is easy to
check that the monodromies $B^J$ and $(A^J)^{-1}(B^J)^{-1}
A^J$ act on $w^I$ according to
\ba
  \pi(B^I) w^J & = &  w^J \k_I^{-1} (R'_{12} R'_{13} R_{13}
          (R_{12}')^{-1}) ^{I J \bar J} \ \ , \nn \\[1mm]
  \pi((A^I)^{-1} (B^I)^{-1} A^I) w^J & = &  w^J \k_I^{-1}
         (R'_{12} R_{12}) ^{I J \bar J} \ \ . \nn
\ea
When these formulas are compared with the expressions
in Theorem \ref{Smrep} we see that the action $B^J$ and
$(A^J)^{-1}(B^J)^{-1} A^J$ on $\End(V^I)$ is equivalent
to the action of the monodromies $M^J_1, M^J_2 \in \L_2$
on the space $V^I \o V^{\bar I}$. The latter was considered in the
previous section and is known to be irreducible by Proposition
\ref{Smrep}. This proves the lemma.

Now we can go back to discuss the irreducibility of $\pi$.
We still have to show that
every subspace $\End(V^I) = \pi(\chi^I_B) \Re$ contains at least one
cyclic vector $\Psi^I \in \pi(\chi^I_B) \Re$. Since $\vac$ is cyclic, it
suffices to find elements $T^I \in \T$ -- one for each
vector $\Psi^I$ -- with $\pi(T^I) \Psi^I = \vac$. Using the notation
$c^I_A = \k_I tr^I_q(A^I)$, we choose $\Psi^I \equiv \pi(c^I_A )\vac
\in \pi(\chi^I_B) \Re $. The standard relations $c^I_A c^J_A = \sum
N^{IJ}_K c^K_A$ furnish
\ba  \pi(T^I) \Psi & = & \pi(\chi^0_B c^{\bar I}_A c^I_A) \vac \nn \\[1mm]
    & = & \pi(\chi^0_B \sum N^{I \bar I}_K c^K_A) \vac
    = \pi(c^0_A) \vac = \vac \nn
\ea
for $T \equiv \chi^0_B c^{\bar I}_A \in \T$. The calculation
employs the fact the $c^K_A$ is constructed from $A^K$ so that
eq. (\ref{chiBA}) furnishes $\pi(\chi^0_B \chi^K_A )\vac =
\pi(\chi^0_A) \vac$. This concludes the proof.

\section{The Moduli Algebra for Higher Genera}
\setcounter{equation}{0}

In this section we consider the most general case of the moduli
algebra on a Riemann surface of arbitrary genus $g$ with an arbitrary
number $m$ of marked points. Our strategy remains similar to preceding
sections. We will represent the graph algebra in certain tensor product
of representations of the underlying Hopf algebra
(cf. all previous sections). Having a set of representations of
the graph algebra, we construct the representations of the moduli
algebra in the gauge invariant subspaces of the tensor products (cf.
Subsection  6.3).

\subsection{Representation theory of the graph algebra $\L_{g,m}$}

The graph algebra $\L_{g,m}$ was defined in Definition \ref{Agn}.
Notice that it consists of $m$ loop algebras (corresponding to
the  monodromies $M^I_\nu$) and $g$ $AB$-algebras. They are
pieced together in a rather standard fashion. It will turn out
that the discussion of the representation theory is completely
parallel to the corresponding arguments in Section 6.2.

So let us begin to construct representations of $\L_{g,m}$. As usual, we
will choose some appropriate representation space of the underlying Hopf
algebra $\G$ on which we then realize $\L_{g,m}$. Let us introduce
the symbol $\Im_g(I_1,\dots , I_m)$ to denote the representation
space
\be \label{Imgn}
 \Im_g(I_1,\dots , I_m) = V^{I_1}\o \dots V^{I_m} \o \Re^{\o g}\ \ .
\ee
With all the experience we gathered in the preceding subsections we
can guess a suitable scalar product for $\Im_g(I_1, \dots, I_m)$
right away. All the $m+g$ factors within the tensor product (\ref{Imgn})
come with a distinguished scalar product. For the carrier spaces $V^I$
of the representations $\t^J$ of the Hopf algebra $\G$ this was part
of the input. On $\Re$ we use the scalar product given through eq.
(\ref{Resp}), so that the natural action of $\G$ on $\Re$ becomes
a *-representation. In conclusion, the space $\Im_g (I_1, \dots, I_m)$
is a $m+g$-fold tensor product of Hilbert spaces, each of which
carries a *-representation of $\G$. Now we are in a position to
employ the Proposition \ref{Imsp} to construct a scalar product
$\langle . , . \rangle$ on $\Im_g (I_1, \dots, I_m)$.

$\Im_g(I_1,\dots , I_m)$ has the best chances to carry a representation
of $\L_{g,m}$ because we know already that its first tensor factors
$V^{I_1}\o \dots \o V^{I_m}$ carry a representation of $\L_m \subset
\L_{g,m}$. The other part in the tensor product (\ref{Imgn}) is a
$g$-fold tensor power of the regular representation. Each copy of
$\Re$ can carry a representation of $AB$-algebra (Theorem \ref{ABrep}).
In fact, if $\L_m$  and the $g$ copies of the $AB$-algebra would come
into $\L_{g,m}$ in the form of a Cartesian product, the construction
of the representation would have been obvious. Even though this is
not the case, we can use our previous experience to define
representations $\Lambda_g (I_1,\dots , I_n)$ of the graph algebra
$\L_{g,m}$ on the spaces $\Im_g (I_1,\dots , I_n)$. Let us introduce
the notation $\pi_i$ for the representation $\pi$ of the $i^{th}$ copy
$\T_i \subset \L_{g,m}$ of the $AB$-algebra implemented in the $i^{th}$
copy of $\Re$ in the tensor product  (\ref{Imgn}), i.e.  for every
$T \in \T_i$,
\be
\pi_i(T) =id_{I_1} \o \dots \o id_{\Re_{i-1}} \o \pi(T) \o
      id_{\Re_{i+1}} \o \dots \o id_{\Re_{g}} \ \ .
\ee
Another useful object is the following representation of the Hopf
algebra $\G$.
\ba
 \jmath_1 & = & \e \ \ \ \ \mbox{ and } \nn \\[1mm]
 \jmath_i (\xi) & = & (\tau^{I_1} \bo
   \dots \bo \pi_{i-1})(\xi) \o id_{\Re_{i}} \o \dots
   \o id_{\Re_g} \ \ (i \geq 2) \nn
\ea
Here $\bo$ denotes the tensor product of representations defined
through the co-product in the standard way. The representation
$\jmath=\jmath_{g+1}$ coincides with the natural action of $\G$ in
the space  $\Im_g(I_1,\dots , I_m)$. Now we are ready to formulate
a theorem.

\begin{theo}  {\em (Representations of the algebra $\L_{g,m}$)}
\label{Sgmrep}
A representation $\Lambda_g^{I_1, \dots, I_m}$ of $ \L_{g,m}$ acts
in the space $\Im_g(I_1,\dots, I_m)$ by means of the following set
of equations.
\ba
 \Lambda_g^{I_1, \dots, I_m} (M^J_\nu) &=&
  D^{I_1,\dots , I_n} (M^J_\nu) \o id_{\Re}^{\o g} \ \ , \nn \\[1mm]
 \Lambda_g^{I_1,\dots , I_m} (A_i^J) &=& (\t^J \o \jmath_i)(R')
  \pi_i(A^J_i) (\t^J \o \jmath_i)((R')^{-1})\ \ , \nn  \\[1mm]
 \Lambda_g^{I_1,\dots , I_m} (B_i^J) &=& (\t^J \o \jmath_i)(R')
  \pi_i(B^J_i) (\t^J \o \jmath_i)((R')^{-1})\ \ . \nn
\ea
$\Lambda_g^{I_1, \dots, I_m}$ extends to a representation of
the semi-direct product $\S_{g,m} = \L_{g,m} \sd \G$ by means of
$$  \Lambda_g ^{I_1,\dots ,I_m} (\xi)  =  \jmath (\xi)\ \ .  $$
When $\Im_g(I_1,.., I_m)$ is equipped with the scalar product
$\langle .,.\rangle$ as discussed after eq. (\ref{Imgn}), the
map $\Lambda_g^{I_1, \dots I_m}: \S_{g,m} \mapsto \End(\Im_g
(I_1, \dots,I_m))$ is a  *-representation.
The set of representations $\{ \Lambda_g ^{\{I_\nu\}}\}$ is faithful
on $\L_{g,m}$. In absence of truncation this implies that the every
representation $ \Lambda_g ^{\{I_\nu\}}$ is irreducible.
\end{theo}

\noindent
{\sc Proof: } On the basis of our prior work, there is nothing
left to be discussed here.
With the relations (\ref{defAgnnew1}) through (\ref{defAgnnew4})
being of the same type as our well known eq. (\ref{defAgn3}), the
proof is straightforward (i.e. essentially identical to the proof
of Theorem \ref{Smrep}).

\subsection{The moduli algebra $\lM^{\{K_\nu\}}_g$}

Having the set of representations $\Lambda_g^{I_1, \dots, I_m}$
of $\L_{g,m}$ at hand we can proceed as in Subsection 6.3. So we
notice that the space $\Im_g(I_1,\dots , I_m)$ carries a
representation of the gauge invariant subalgebra $\A_{g,m}$
within $\L_{g,m}$. Under the action of $\A_{g,m}$ the space
$\Im_g(I_1,\dots , I_m)$ splits into a sum of irreducible
representations. As is the case of $\L_m$, this decomposition reads
\be \label{decg}
  \Im_g (I_1, \dots, I_m) =
  \sum_J V^J \otimes W_g^J(I_1,\dots, I_m)\ \ ,
\ee
where $J$ runs through all equivalence classes of irreducible
representations of $\G$ and $W^J_g(I_1, \dots, I_m)$ are
multiplicity spaces.  We note that
the scalar product $\langle .,. \rangle $ on $\Im_q(I_1, \dots,I_m)$
restricts to $W^J_g (I_1,\dots , I_m)$ and that these spaces carry
an irreducible *-representation of $\A_{g,m}$. These representations
can be used to obtain representations of the moduli algebra  which we
discuss next.

Recall that we have $m$ marked points on the Riemann surface.
They are labeled by $K_1, \dots K_m$. Each of these points
is surrounded by one of the loops $l_\nu$. To implement the
corresponding flatness relations in the finite dimensional
case, we use the characteristic projectors $\chi^{K_\nu}_\nu$
which are constructed from monodromies $M^I_\nu$. From the
arguments in Section 6.1 one concludes that these projectors
are central in $\A_{g,m} \subset \L_{g,m}$. Still we need one
more projector that comes with the total product $r_g =
[b_g, a_g^{-1}] \dots [b_1, a_1^{-1}] l_m \dots l_1$. As in
section 3.5 we define
\ba
   \chi^K_0 & \equiv & \N d_K S_{K \bar I} c^I_0 = \N d_K \k_K
    S_{K \bar I} tr^I_q(M^I(r_g)) \ \label{chi0g}\\[1mm]
    \mbox{ with } & & M^I(r_g) = \k_I^{4g+m-1} [B_g^I, (A_g^I)^{-1}]
    \dots M^I_1\ \ .
    \nn
\ea
Here $[B^I_i, (A_i^I)^{-1}] = B^I_i (A^I_i)^{-1} (B^I_i)^{-1}
A^I_i$. The elements $\chi^K_0$ generate a fusion algebra and
are {\em central} in $\A_m$ (see \cite{AGS2}).

\begin{defn} {\em (Moduli algebra for arbitrary surface)}
The moduli algebra $\lM^{\{K_\nu\}}_{g,m}$ of a surface of genus g
with $m$ punctures marked by $K_\nu, \nu = 1, \dots, m,$ is the
$*$-algebra
\be
   \lM^{\{K_\nu\}}_{g,m} \equiv \chi^0_0 \prod_{\nu = 1}^m
    \chi^{K_\nu} _\nu    \A_{g,m} \ \ .
\ee
Here $\chi^{K_\nu}_\nu, \chi^0_0 \in \A_{g,m}$ are the central
projectors introduced in the text preceding this definition.
\end{defn}

Along the lines of the corresponding discussion in Section 6.3,
we can evaluate the projectors $\chi^K\in \A_{g,m}$ on the
representation spaces $W_g^J(I_1, \dots,I_m)$. For the
character $\chi^0_0$ one uses the fact that the product
$ \tilde M^I_i = \k^3_I B^I_i (A_i^I)^{-1} (B_i^I)^{-1}
A_i^I$ which is assigned to $[b_i,a^{-1}_i]$, embeds the loop
algebra so that the element assigned to $r_g$ looks as if
it would come from $m+g$ loops. Because of eq. (\ref{pirhoM}),
this holds also for the way the elements are represented on
$\Im_g(I_1, \dots I_m)$. Consequently, the evaluation of
characters on multiplicity spaces is identical to the
calculation we did in the Section 6.3 and we can state
the results right away.

\begin{theo} {\em (Representations of the moduli algebra, genus g)}
For any set  $K_1,$ $\dots$, $K_m $ labeling $m$ points on a
Riemann surface of genus $g$, there is a unique irreducible
*-representation of the corresponding moduli algebra $\lM^{\{
K_\nu\}}_{g,m}$ on the space $W_g^0(K_1, \dots, K_m)$ (defined
through the decomposition (\ref{decg})). This representation
can be obtained explicitly by restricting the representation
$\Lambda^{K_1, \dots K_m}_g$ of $\L_{g,m}$ to the moduli algebra
$\lM^{\{K_\nu\}}_{g,m} \subset \L_{g,m}$.
\end{theo}

\subsection{The second pinching theorem}

We would like to extend our description of relations between
moduli algebras that was initiated in Section 6.4. It is
obvious how Proposition \ref{1pintch} carries over to the
moduli algebras on higher genus surfaces. But now we have
another possibility: we can also shrink the surface along
a circle that wraps around some handle. To be specific,
we chose the fusion algebra ${\cal V}(b_1)$ constructed
from the elements $B_1^I$ and evaluate its commutant.
Recall that ${\cal V}(b_1)$ is generated by
$$     c^I  = \k_I tr^I_q(B_1^I) \ \ . $$

\begin{prop} {\em (Second  Pinching Theorem) } The commutant
of ${\cal V}(b_1)$ in the moduli algebra $\lM_{g, m}^{ K_1, \dots
, K_m}$ splits into the direct sum of moduli algebras of genus
$g-1$ with $m+2$ marked points,
\be \label{ThC}
{\cal C}({\cal V}(b_1), \lM_{g, m}^{ K_1, \dots , K_m}) \cong
\bigoplus_K  \lM_{g-1, m+2}^{K_1, \dots , K_m, K, \bar{K}}.
\ee
Here the sum runs over all classes of irreducible representations
of the symmetry Hopf algebra.
\end{prop}

In the language of Topological Field Theory, evaluation of the
commutant corresponds to shrinking the cycle $b_1$ so that we
get a surface of lower genus. It has two marked points at the
place where the handle is pinched -- one on either side of
the cut. Shrinking all the $b$-cycles one after another, one
produces spheres with $m+2g$ marked points.

\noindent
{\sc Proof:} To prove the theorem  we proceed as in Section 6.4.
Let us first show how to embed the direct sum on the right hand
side of (\ref{ThC}) into the moduli algebra. The set of
matrix generators $M^I_\nu\ (\nu=1, \dots, m)\ , \k_I^2
B^I_1 (A^I_1)^{-1} (B^I_1)^{-1}, A^I_1$ and $ A^I_i,B^I_i\  (i = 2,
\dots, g)$ can be regarded as an image of the standard generators
of $L_{g-1,m+2}$ under a map $\phi: \L_{g-1,m+2} \mapsto
\L_{g,m}$. $\phi$ is not an embedding. But the proof for
Theorem  \ref{ABrep} at the end of Subsection 7.3
shows that $\phi$ restricts to an embedding of the direct sum
$\bigoplus_K  \lM_{g-1, m+2}^{\{K_1, \dots , K_m, K, \bar{K}\}}$
into the moduli algebra $\lM_{g, m}^{\{ K_1, \dots , K_m\}}$
(cp. also Lemma \ref{LDembed}).
Counting dimensions we find that every summand is embedded
with multiplicity one so that the commutant of the image of
$\phi$ is exactly the fusion algebra. For details the reader is
referred to the proof of the Proposition \ref{1pintch}.
\section{Representations of Mapping Class Groups}
\setcounter{equation}{0}

The action of the mapping class group $M(g,m)$ on the fundamental
group of a surface $\Sigma_{g,m}$ of genus $g$ with $m$ punctures
induces an action on the graph algebras
$\L_{g,m}$, i.e. a homomorphism from $M(g,m)$ into the automorphism
group of $\L_{g,m}$. In general, only automorphisms corresponding
to the pure mapping class group $PM(g,m) \subset M(g,m)$ restrict
to automorphisms of the moduli algebras $\lM_{g,m}^{\{K_\nu\}}$. Since
the moduli algebras are simple, every automorphism is inner and
can be implemented by a unitary element. We will give a simple
prescription to construct such unitaries for a generating set of
elements in $PM(g,m)$. They will furnish projective representations
of $PM(g,m)$. The latter are equivalent to those found by
Reshetikhin and Turaev \cite{ReTu} (extended to the possible
presence of punctures)\footnote{As we discussed in \cite{AGS1,AGS2},
our theory generalizes to quasi-Hopf algebras. So the discussion
of mapping class groups covers the extension of the Reshetikhin-
Turaev construction found by Altschuler and Coste \cite{AlCo}.}.

\subsection{Action of the mapping class group on moduli algebras}

To begin with, let us recall that the quantum monodromies
$M_\nu^I, A^I_i, B^I_i \in \L_{g,m}$ are assigned to the
generators $l_\nu, a_i, b_i \in \pi_1(\Sigma_{g,m}\setminus D)$
\footnote{If we remove a disk $D$ from the surface, the fundamental
group is freely generated by $l_\nu, a_i,b_i$. We will later glue
the disk back and thereby introduce the relation \ref{abm}}.
We want to display this more clearly in the notations and
hence introduce
$$
   M^I(l_\nu) \equiv M^I_\nu \ \ , \ \
   M^I(a_i) \equiv A^I_i \ \ , \ \
   M^I(b_i) \equiv B^I_i \ \ .
$$
We can go even further and define $M^I(\C)$ for arbitrary
elements $\C \in \pi_1(\Sigma_{g,m}\setminus D)$. Suppose that
$\C = \C_1 \C_2 $ with two elements $\C_i, i = 1,2,$ in the
fundamental group of $\Sigma_{g,m}$. Then we set
\be
   M^I(\C)  \equiv   \k_I^{w(\C_1,\C_2)} M^I(\C_1)
   M^I(\C_2)  \label{MonC} \ \ .
\ee
Here $w(\C_1, \C_2) = \pm 1$ is determined from the
pair $(\C_1, \C_2)$ as follows. Suppose that we have
presentations of $\C_1$ and $\C_2$ in terms of generators
$\l_\nu^{\pm}, a_i^{\pm},b_i^{\pm}$. We assume that
these presentations are reduced so that neighboring elements
are never inverse to each other. From these presentations
we extract two generators $c_i \in \{\l_\nu^{\pm},
a_i^{\pm},b_i^{\pm}\}$ such that $\C_1 = \C_1' c_1$
and $\C_2 = c_2 \C_2'$. Without restriction we can assume
that $c_1 \neq c_2^{-1}$. The weight $w$ that we want to
describe satisfies $w(\C_1, \C_2) = w(c_1,c_2)$. To define
the latter we need to introduce two maps $t,s$ from the
set of generators $\{ l_\nu^{\pm}, a_I^{\pm}, b_i^\pm \}$
to the integers $1, \dots, 2m+4g$.
\ba
       t(l_\nu) = 2 \nu = s(l_\nu) + 1 \ \ & , & \ \
       t(a_i) = 2m + 4i = s(a_i) + 2 \ \ , \nn \\[1mm]
       t(b_i) = & 2m + 4i -1 & = s(b_i) + 2 \ \ \nn
\ea
and $t(c^-) = s(c)$ for every $c \in \{ l_\nu^{\pm},
a_I^{\pm}, b_i^\pm \}$. Now we can complete the description
of the weight $w$.
\be
       w(c_1,c_2) \equiv \left\{ \begin {array}{l}
       +1 \ \ \mbox{ if } \ \ t(c_1) < s(c_2) \\[1mm]
       -1 \ \ \mbox{ if } \ \ t(c_1) > s(c_2)
       \end{array} \right.
\ee
for all $c_1,c_2 \in  \{ l_\nu^{\pm}, a_I^{\pm}, b_i^\pm \}$.
Our choice of the weight factor in the definition (\ref{MonC})
is designed such that all monodromies $M^I(\C)$ satisfy
functoriality. We have seen particular examples of this in
subsection 7.3.

The mapping class group $M(g,m)$ of a surface $\Sigma_{g,m}$ is
defined as the group of diffeomorphisms of $\Sigma_{g,m}$
into itself modulo its identity component. Similarly, $M(g,m;B)$
is obtained from diffeomorphisms of $\Sigma_{g,m} \setminus D$
which leave the boundary $B = \partial D$ pointwise fixed.
Elements in $M(g,m)$ and $M(g,m;B)$ may interchange the punctures.
This furnishes the usual canonical homomorphism from
mapping class groups into the symmetric group. The kernel of
this homomorphism is called pure mapping class group. We will
denote it by $PM(g,m)$ and $PM(g,m;B)$.

It is well known that elements $\varrho$ in the  mapping class group
$M(g,m;B)$ of the $m$-punctured surface act on the fundamental group
$\pi_1(\Sigma_{g,m} \setminus D)$ as outer automorphisms. We will
not distinguish in notation between elements $\varrho$ in $M(g,m;B)$
and the corresponding elements $\varrho \in \Aut (\pi_1(\Sigma_{g,m}
\setminus D))$. The action of the mapping class group on the
fundamental group lifts to an action on the graph algebras by
means of the formula
$$
 \hat \varrho ( M^I(\C)) = M^I(\varrho(\C)) \ \ \mbox{ for all } \ \
   \C \in \pi_1(\Sigma_{g,m}\setminus D)\ \ .
$$
Ultimately, we are more interested in automorphisms of moduli
algebras. As a first step in the reduction from the graph algebra
$\L_{g,m}$ to the moduli algebras, one notices that the
automorphisms $\hat \varrho $ are consistent with the transformation
law under the action of $\xi \in \G$.
Consequently, the action of $\hat \varrho$ on $\L_{g,m}$ descends to
an action on $\A_{g,m}$ (for notations compare Section 8). To
proceed towards the moduli algebra, two more steps are necessary. First
we have to multiply $\A_{g,m}$ by the projectors $\chi^{K_\nu}_\nu$,
then we need to implement flatness for the circle $r_g = [b_g, a_g^{-1}]
\dots l_1$  with the help of $\chi^0_0$. It is intuitively clear that
only automorphisms assigned to elements in the pure mapping class
group $PM(g,m;B)$ survive the first step which corresponds to coloring
the punctures. To understand the effect of $\chi^0_0$, we recall
the relation (a proof can be found in \cite{AGS2})
\be
   \chi^0_0 M^L(r_g) = \chi^0_0 \k_L^{-1} e^L\ \ .
    \label{chimr}
\ee
It implements the defining relation $r_g =id$ for $\pi_1(\Sigma_{g,m})$
into the moduli algebra. One may check that all the automorphisms
constructed from elements in the (pure) mapping class group respect
this relation (i.e. $\chi^0_0$ in invariant under their action).
On the other hand, some of the elements in $PM(g,m;B)$ act trivially
on the moduli algebras because $M^L(r_g) \sim \k_L^{-1} e^L$. A more
detailed investigation shows that nontrivial automorphisms correspond
to elements in $PM(g,m)$. We formulate this as a proposition.

\begin{prop} {\em (Action of mapping class group)}
\label{MCGact} For every
$\varrho \in M(g,m;B)$ there exits an automorphism $\hat \varrho:
\L_{g,m} \mapsto \L_{g,m}$ of graph algebras,
$$ \hat \varrho ( M^I(\C)) = M^I(\varrho(\C)) \ \ \mbox{ for all } \ \
       \C \in \pi_1(\Sigma_{g,m}\setminus D)\ \ . $$
These automorphisms furnish an action of the mapping class group
$M(g,m;B)$ on the graph algebra $\L_{g,m}$. Automorphisms $\hat
\eta$ corresponding to elements $\eta \in PM(g,m;B)$ restrict to
the moduli algebras $\lM_{g,m}^{\{K_\nu\}}$ and give rise
to an action of the pure mapping class group  $PM(g,m)$
on moduli algebras.
\end{prop}

The proof of this proposition is technically not very difficult.
Given an explicit description of the action of $M(g,m;B)$ on the
fundamental group, it is essentially based on calculations similar
to those performed in Section 7.3. Guided by the explanations above,
the reader may try to verify the result. A complete proof will also
appear in \cite{ASh2}.

\subsection{Projective representations of mapping class groups}

In Proposition \ref{MCGact} we obtained an action of the
pure mapping class group on moduli algebras. Let us recall that
the moduli algebras $\lM^{\{K_\nu\}}_{g,m}$ are simple so that
all automorphisms are inner. This means that for every element
$\eta$ in  the pure mapping class group there is a unitary element
$\hh(\eta)$  in the moduli algebra so that
$$    \hh(\eta) \ A = \hat \eta (A) \hat h(\eta) \ \ . $$
Such elements $\hat h$ provide a projective representation of
the pure mapping class group. In fact, let $\eta_i, i= 1,2,$ be
two elements in $PM(g,m)$ and denote corresponding elements
in the moduli algebra by $\hat h_i= \hh (\eta_i)$. Suppose that the
product $\heta = \hat \eta_1 \hat \eta_2$ is implemented by $\hat h$.
Then
$$  \hat h^* \hat h_1 \hat h_2\  A = A \  \hat h^* \hat h_1
    \hat h_2 $$
for all elements $A \in \lM^{\{K_\nu\}}_{g,m}$. Since the moduli
algebra is simple, only scalars can commute with all elements
and hence $\hat h_1 \hat h_2 = \varpi \hat h$ with a complex
number $\varpi$.

In constructing this representation, the only problem is to
find explicit expressions for the elements $\hat h$, at least
for a generating set of elements $\eta \in PM(g,m)$. This is
surprisingly simple. Suppose that $\eta$ is a Dehn twist along
the circle $x(\eta)$. Let us regard $x(\eta)$ as
an element in the fundamental group so that
$M^I(x(\eta))$ is well defined. Then we set
\ba
    \hat h (\eta)& \equiv & \sum \th^{-1} \N d_I v_I
    \k_I tr_q^I(M^I (x(\eta))) \ \label{hdef} \\[1mm]
    \mbox{ where } & & \th = \sum \N v_I d_I^2\ \ \label{th}
\ea
and $\N =(\sum d_I^2 )^{-1/2}$. We wish to demonstrate that
$\hat h (x)$ is unitary. Setting $c_x^L \equiv \k_L tr^L_q
( M^L(x)) $ we obtain
\ba
\th \th^* \hh (x) \hh (x)^*
          & = & \sum \frac{v_I}{v_J} d_I d_J \N^2
                 c_x^I c_x^J \nn \\[1mm]
          & = & \sum \frac{v_I}{v_J} d_I d_J \N^2
                 N^{IJ}_K c^K_x \nn \\[1mm]
          & = & \sum  \N \frac{1}{v_K} d_I S_{I \bar K}
                c^K_x \nn \\[1mm]
          & = & \sum S_{0 I} S_{I \bar K} \frac{1}{v_K}
                c_x^K \nn \\[1mm]
          & = & \sum C_{0 \bar K} \frac{1}{v_K} c_x^K = 1
               \ \ .     \nn
\ea
Here we made use of the following expression for the matrix $S$
$$
S_{I \bar K} = \N \frac{v_I v_K}{v_J} N^{I \bar K} _{\bar J}
                   d_{\bar J} \ \ .
$$
Now we recall that $\vth (c_x ^I) = d_I $ defines a
representations $\vth$ of the fusion algebra $\V (x)$ over the
circle $x$. The definition (\ref{hdef}) shows that in particular
$\vth( \hh (x)) = 1 $. If we apply $\vth$ to the equality
$1 = \th \th^* \hh (x) \hh (x)^* $, we see that $1 = \vth(\th
\th^* \hh (x) \hh (x)^*) = \th \th^*$. This means that $\th$
is a phase and hence
\be
        \hh (x)  \hh(x )^* = \th ^* \th = 1
\ee
Notice that $\hat h (x) $ is an element in the fusion algebra
${\cal V} (x)$ over $x$. In particular
it can be expressed as a linear combination of the characters
$\chi^I(x)$. The outcome of a short calculation using
the same ideas as the above proof for the unitarity of $\hh
(x)$ is
\be
     \hat h(x)
     = \sum v^{-1}_I \chi^I (x) \ \ .
     \label{hchi}
\ee
Our result (\ref{hchi}) shows that
$\hat h(x)$ is simply the {\em inverse of the ribbon
element $v(x)$ over the circle $x$}.

\begin{theo} {\em (Representation of $PM(g,m)$) }
\label{repMCG}
Suppose that $\eta \in PM(g,m) $ is a Dehn twist along the
circle $x(\eta)$ on the surface $\Sigma_{g,m}$. Then the unitary
element $\hh (x(\eta))$ defined through eq. (\ref{hdef})
implements the action of $\eta \in PM(g,m)$ on the moduli
algebras, i.e.
 $$
   \hat h(x(\eta)) A = \hat \eta( A ) \hat h(x(\eta)) \ \
 $$
holds for all elements $A$ in the moduli algebras.
The map $\eta \mapsto \hat h(x(\eta))$ defines a unitary
projective representation of the pure mapping class group
$PM(g,m)$ (recall that Dehn twists $\eta$ generate $PM(g,m)$).
\end{theo}

In spite of its appearance, the preceeding Theorem is
relatively cumbersome to prove because a lot of cases have to
be investigated separately. To keep this work compact, we will
give the proof elsewhere \cite{ASh2}.

\subsection{Equivalence with Reshetikhin Turaev representation}

We turn now to
a more explicit description of the projective representation
which we have just constructed. The plan is to evaluate the
action of $\hh (x(\eta))$ on the states of the $CS$-theory. This
will then allow to compare our representation with the action
of mapping class groups on conformal blocks and the Reshetikhin-
Turaev representation.

Let us begin by stating the main theorem of this section.

\begin{theo} {\em (Equivalence with RT-representation)}
\label{RTequiv}
The projective representation $\eta \mapsto \hh (x(\eta)) \in \A_{g,m}$
of the pure mapping class group $PM(g,m)$ is unitarily equivalent
the representation found by Reshetikhin and Turaev in \cite{ReTu}.
\end{theo}

The rest of the section is devoted to the proof of this theorem.
For this comparison of our representation with the Reshetikhin-
Turaev representation it suffices to evaluate the elements $\hh
(x(\eta))$ for a generating set of Dehn twists. We use a set which
consists of Dehn twists $\eta_{\nu p}, 1 \leq \nu \leq m,
\nu < p \leq m+2g$ and  $ \a_i, \b_i, \d_i, e_i, i = 1,
\dots g,$ along the circles $x(\eta_{\nu p})$ and $x(\a_i),
x(\b_i), x(\d_i), x(\e_i)$. More precisely we define
\ba
  x(\eta_{\nu p}) = l_p l_\nu  \ \ & , & \ \ x(\a_1) = l_{m+1}
   \ \ , \nn \\[1mm]
  \mbox{ and }\   x(\a_j)  & = & l_{m+2j-1} l_{m+2j-2}
  \ \ \mbox{ for } \ \ j = 2, \dots, g\ \ , \nn \\[1mm]
  x(\b_i) = a_i\ \ &,&\ \ x(\d_i)=l_{2i-1}\ \ ,\nn\\[1mm]
  x( \e_i) = l_{m+2i-1} \dots l_{m+1}\ \ & & \mbox{ for }
  \ \  i = 1, \dots, g\ \ . \nn
\ea
where $l_\nu, \nu = 1, \dots, m$ and $a_i, i = 1, \dots, g,$
are defined as before and we introduced $l_p, m< p < m+2g$ by
$$ l_{m+2i-1} \equiv a^{-1}_i b^{-1}_i a_i \ \ ,
   \ \ l_{m+2i} \equiv b_i \ \  $$
for $i = 1, \dots, g$. The circles are shown in the Figure
1,2 .

\begin{picture}(300,210)(5,-45)  \it
\setlength{\unitlength}{.9pt}
\def\handle{\begin{picture}(70,70)(0,0) \put(20,0){\line(0,1){45}}
\put(35,0){\line(0,1){45}} \put(55,0){\line(0,1){45}}
\put(70,0){\line(0,1){45}} \put(45,45){\oval(20,20)[t]}
\put(45,45){\oval(50,50)[t]} \end{picture}}
\thicklines \linethickness{0.4mm}
\put(35,80){\oval(50,100)[l]} \put(345,80){\oval(50,100)[r]}
\put(35,30){\line(1,0){310}} \put(20,90){\handle}
\put(110,90){\handle} \put(200,90){\handle} \put(300,90)
{\circle*{3}} \put(325,90){\circle*{3}} \put(350,90){\circle*{3}}
\put(60,130){\line(1,0){10}} \put(150,130){\line(1,0){10}}
\put(240,130){\line(1,0){10}}
\put(95,130){\line(1,0){30}} \put(185,130){\line(1,0){30}}
\put(275,130){\line(1,0){70}} \put(102,90){$\dots$}
\put(192,90){$\dots$} \put(307,90){$\dots$} \put(332,90){$\dots$}
\put(327,82){$\nu$} \put(352,82){1}
\put(302,82){m}  \put(255,82){\tiny{m+1}}
\put(220,82){\tiny{m+2}} \put(137,82){\tiny{p}}
\put(38,82){\tiny{m+2g}} \put(325,100){\oval(20,20)[t]}
\put(150,100){\oval(20,20)[tr]} \put(125,100){\oval(20,20)[tl]}
\put(115,75){\line(0,1){25}} \put(160,75){\line(0,1){25}}
\put(315,75){\line(0,1){25}} \put(335,75){\line(0,1){25}}
\put(170,75){\oval(20,20)[bl]} \put(170,75){\oval(110,70)[bl]}
\put(305,75){\oval(20,20)[br]} \put(305,75){\oval(60,70)[br]}
\put(170,65){\line(1,0){135}} \put(170,40){\line(1,0){135}}
\rm
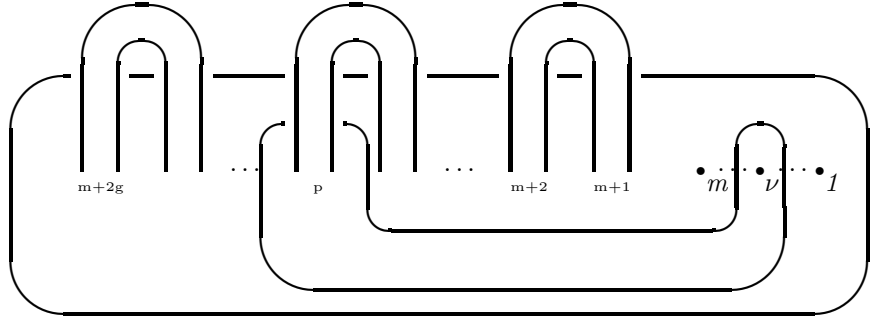
\put(25,5){\parbox[t]{290pt}{\rm \small {\bf Fig. 1:} Curves
$x(\eta_{\nu p}), 1 \leq \nu \leq m, \nu < p \leq m+2g$ on a surface
$\Sigma_{g,m}$ with $m$ marked points wrap around the $\nu^{th}$
puncture and the $p^{th}$ ``leg'' if $m < p \leq m +2g$. If
$p \leq m$, $x(\eta_{\nu p})$ encloses the $\nu^{th}$ and the $p^{th}$
puncture.}}
\end{picture}
\vspace*{-1cm}

\begin{picture}(300,190)(5,-35)  \it
\setlength{\unitlength}{.8pt}
\def\handle{\begin{picture}(70,70)(0,0) \put(20,0){\line(0,1){45}}
\put(35,0){\line(0,1){45}} \put(55,0){\line(0,1){45}}
\put(70,0){\line(0,1){45}} \put(45,45){\oval(20,20)[t]}
\put(45,45){\oval(50,50)[t]} \end{picture}}
\thicklines \linethickness{0.4mm}
\put(35,80){\oval(50,100)[l]}
\put(370,30){\line(1,0){5}}
\put(370,130){\line(1,0){5}}
%\put(365,80){\oval(50,100)[r]}
\put(35,30){\line(1,0){330}} \put(20,90){\handle}
\put(110,90){\handle}  \put(290,90){\handle}
\put(200,90){\line(0,1){45}} \put(215,90){\line(0,1){45}}
\put(255,90){\line(0,1){45}} \put(270,90){\line(0,1){45}}
\put(225,135){\oval(20,20)[tl]} \put(225,135){\oval(50,50)[tl]}
\put(245,135){\oval(20,20)[tr]} \put(245,135){\oval(50,50)[tr]}
\put(225,145){\line(1,0){20}} \put(225,160){\line(1,0){20}}
\put(60,130){\line(1,0){10}} \put(150,130){\line(1,0){10}}
\put(220,130){\line(1,0){30}} \put(330,130){\line(1,0){10}}
\put(95,130){\line(1,0){30}} \put(185,130){\line(1,0){10}}
\put(275,130){\line(1,0){30}} \put(102,90){$\dots$}
\put(282,90){$\dots$}
\put(208,57){\line(0,1){78}} \put(262,57){\line(0,1){78}}
\put(225,135){\oval(34,34)[tl]} \put(245,135){\oval(34,34)[tr]}
\put(225,57){\oval(34,34)[bl]} \put(245,57){\oval(34,34)[br]}
\put(225,152){\line(1,0){20}} \put(225,40){\line(1,0){20}}
\put(160,100){\oval(20,20)[tl]} \put(220,100){\oval(20,20)[tr]}
\put(160,80){\oval(20,20)[bl]} \put(220,80){\oval(20,20)[br]}
\put(185,110){\line(1,0){10}} \put(160,70){\line(1,0){60}}
\put(150,80){\line(0,1){20}} \put(230,80){\line(0,1){20}}
\put(250,100){\oval(20,20)[tl]} \put(250,80){\oval(20,20)[bl]}
\put(240,80){\line(0,1){20}} \put(370,80){\line(0,1){20}}
\put(360,100){\oval(20,20)[tr]} \put(360,80){\oval(20,20)[br]}
\put(250,70){\line(1,0){110}}  \put(275,110){\line(1,0){30}}
\put(330,110){\line(1,0){10}}
\put(70,100){\oval(20,20)[tl]} \put(95,100){\oval(20,20)[tr]}
\put(70,80){\oval(20,20)[bl]} \put(95,80){\oval(20,20)[br]}
\put(70,70){\line(1,0){25}}
\put(60,80){\line(0,1){20}} \put(105,80){\line(0,1){20}}
\put(40,157){g} \put(130,157){i}  \put(195,157){i-1}
\put(310,157){1}
\put(90,60){$x(\d_g)$} \put(170,60){$x(\a_i)$}
\put(340,60){$x(\e_{i-1})$} \put(265,40){$x(\b_{i-1})$}
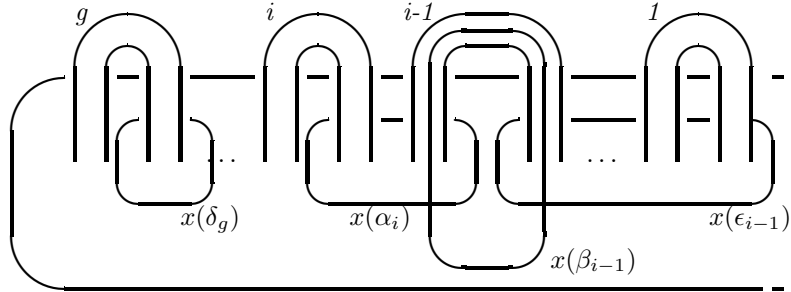
\put(20,5){\parbox[t]{290pt}{\rm \small {\bf Fig. 2:} Curves
$x(\a_i), x(\b_i), x(\d_i), x(\e_i), i = 1, \dots ,g,$ on a
surface $\Sigma_{g,m}$  along which the Dehn twists
$\a_i, \b_i, \d_i, \e_i$ are performed. If there are
punctures on the surface, they are not encircled by
the curves. }}
\rm \end{picture}

As in Lemma \ref{LDembed} we can map the
representation spaces $ \Im (K_1, \dots, K_m) \o \Re^{\o_g}$ to
$\bigoplus_{\{ J_i\}} \Im (K_1, \dots, K_m) \o V^{J_1} \o
V^{\bar J_1} \o \dots \o V^{J_g} \o V^{\bar J_g}$. This
space carries a representation $D_g^{K_1, \dots, K_m}$
of the graph algebra $\L_{g,m}$. The representation $D_g^{K_1,
\dots,K_m}$ is unitarily equivalent to
$\Lambda_g^{K_1, \dots,K_m}$.

Now we evaluate the basic generators of the pure mapping class
group in the representations $D_g^{K_1, \dots,K_m}$.

\begin{lemma}
\label{etaprep} In the representations $D_g^{\{K_\nu\}}$ on
$\bigoplus_{\{ J_i\}} \Im (K_1, \dots, K_m) \o V^{J_1} \o
V^{\bar J_1} \o \dots \o V^{J_g} \o V^{\bar J_g}$ the
elements $\hh_{\nu p} \in \L_{g,m}, 1 \leq \nu \leq m,
\nu < p \leq m+2g,$ are represented by
$$
    D_g^{\{K_\rho\}} (\hat h_{\nu p})
    =  \bigoplus_{\{J_i\}}
   \left( v_\nu^{-1} v_p^{-1}  Q_{\nu p}\right)
   ^{K_1 \dots K_m J_1  \bar J_1 \dots \bar J_g}
$$
and elements $\ha_i \hb_i, \hd_i, \he_i, i = 1, \dots , g,$ are
represented by
\ba
    D_g^{\{K_\rho\}} (\ha_i) & = &
    \bigoplus_{\{J_j\}} \left( v^{-1}_{m+2i-2} v^{-1}_{m+2i-1}
      Q_{(m+2i-2)(m+2i-1)} \right)^{K_1 \dots K_m
      J_1 \bar J_1 \dots J_g \bar J_g} \ \ , \nn \\[1mm]
    D_g^{\{K_\rho\}} (\hb_i) & = &  \th^{-1}
    \bigoplus_{\{J_j\}} \sum_K \N d_K  v_K v_{\bar J_i}
      C[J_i \bar J_i |0]  (R')^{J_i K} (R')^{\bar K \bar J_i}
      C[ K \bar K|0]^c \ \ ,
    \nn \\[1mm]
    D_g^{\{K_\rho\}} (\hd_i) & = & \bigoplus_{\{J_j\}}
    (v_{m+2i-1})^{K_1 \dots K_m  J_1 \bar J_1 \dots J_g
    \bar J_g}
    \ \ , \nn \\[1mm]
    D_g^{\{K_\rho\}} (\he_i) & = &
    \bigoplus_{\{J_j\}} \left( v^{-1}_{m+1}..v^{-1}_{m+2i-1}
    Q_{(m+1)(m+2i-1)}..Q_{(m+2i-2)(m+2i-1)} \right)
    ^{K_1 \dots \bar J_g} .\nn
\ea
Here $\hh_{\nu p} = \hh (x(\eta_{\nu p})), \ha_i = \hh (x(\a_i)),
etc$. To simplify the expressions we also used the element
$Q_{\nu \mu} \in \G^{\o_{m+2g}}$,
\ba Q_{\nu \mu}  & = &  R'_{\nu (\nu+1) }
    \dots R'_{\nu \mu} R_{\nu \mu}  (R'_{\nu (\mu- 1)})^{-1}
   \dots  (R'_{\nu (\nu+1)})^{-1} \in \G^{\o_m}  \ \ . \nn \\[1mm]
     \mbox{ and } \ & &  \ C[K \bar K|0]^c  \equiv
  (R')^{K \bar K} C[ K \bar K|0]^* \frac{d_K}{v_K}\ \ . \nn
\ea
Fig. 3,4 give a graphical presentation of the result.
\end{lemma}

\vspace*{2cm}
\begin{picture}(160,100)(0,0)
\thicklines \linethickness{0.4mm}
\put(00,30){\line(1,0){150}} \put(00,130){\line(1,0){150}}
\put(10,30){\line(0,1){100}} \put(140,30){\line(0,1){100}}
\put(40,30){\line(0,1){25}}  \put(60,30){\line(0,1){25}}
\put(100,30){\line(0,1){25}} \put(120,30){\line(0,1){20}}
\put(40,65){\line(0,1){65}}  \put(60,65){\line(0,1){30}}
\put(100,65){\line(0,1){30}} \put(120,110){\line(0,1){20}}
\put(60,105){\line(0,1){25}} \put(25,70){\line(0,1){20}}
\put(100,105){\line(0,1){25}}\put(45,100){\line(1,0){65}}
\put(35,60){\line(1,0){75}}  \put(35,90){\oval(20,20)[tl]}
\put(35,70){\oval(20,20)[bl]} \put(110,110){\oval(20,20)[br]}
\put(110,50){\oval(20,20)[tr]} \put(73,80){$\dots$}
\put(36,20){$p$} \put(117,20){$\nu$}
\put(40,42){\circle{5}} \put(120,42){\circle{5}}
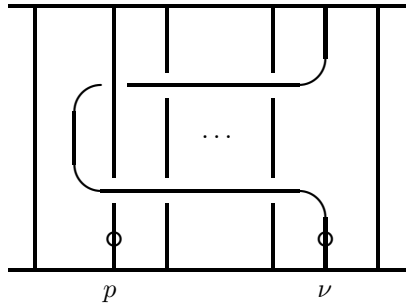
\put(180,20){\parbox[b]{5cm}{\rm \small {\bf Fig. 3:} Pictorial
presentation of Lemma \ref{etaprep}. The picture contains
m + 2g strands, one for every puncture and two for every handle.
It shows
the action of $\eta_{\nu p}, 1 \leq \nu \leq m, \nu \leq p \leq
m+2g$ on states. Over- and undercrossings correspond to factors
$R, R^{-1}$ and the open circle on a line $s$ means
multiplication with $v_{K_\s } ^{-1}$. }}
\end{picture}

\begin{picture}(320,380)(-5,-50)
\def\obr{ \begin{picture}(40,10)(0,0)
\put(-6,0){$\overbrace{\phantom{======}}$} \end{picture}}
\thicklines \linethickness{0.4mm}
\put(10,320){\obr} \put(90,320){\obr} \put(230,320){\obr}
\put(50,170){\obr}  \put(190,170){\obr} \put(270,170){\obr}
\put(29,333){$i$} \put(102,333){$i-1$}
\put(249,333){$i$} \put(69,183){$i$}
\put(209,183){$i$} \put(289,183){$1$}
\put(0,320){\line(1,0){140}} \put(0,240){\line(1,0){140}}
\put(0,170){\line(1,0){140}} \put(0,90){\line(1,0){140}}
\put(180,320){\line(1,0){140}}  \put(180,240){\line(1,0){140}}
\put(180,170){\line(1,0){140}} \put(180,90){\line(1,0){140}}
\put(70,210){$\ha_i$}   \put(250,210){$\hb_i \th$ }
\put(70,60){$\hd_i$}   \put(250,60){$\he_i $}
\put(10,90){\line(0,1){80}}  \put(50,90){\line(0,1){80}}
\put(90,90){\line(0,1){80}}  \put(130,90){\line(0,1){80}}
\put(90,130){\circle{5}}
\put(225,123){\oval(20,16)[l]}  \put(300,107){\oval(20,16)[tr]}
\put(260,138){\oval(20,16)[tr]}  \put(300,139){\oval(20,16)[br]}
\put(225,154){\oval(20,16)[l]}  \put(260,170){\oval(20,16)[br]}
\put(190,90){\line(0,1){80}}  \put(230,90){\line(0,1){20}}
\put(270,90){\line(0,1){20}}  \put(310,90){\line(0,1){17}}
\put(310,139){\line(0,1){31}} \put(230,120){\line(0,1){21}}
\put(230,151){\line(0,1){19}} \put(310,139){\line(0,1){31}}
\put(270,120){\line(0,1){18}} \put(225,115){\line(1,0){75}}
\put(235,131){\line(1,0){30}} \put(275,131){\line(1,0){25}}
\put(225,146){\line(1,0){35}}  \put(235,162){\line(1,0){25}}
\put(240,138){\circle*{2}}    \put(250,138){\circle*{2}}
\put(260,138){\circle*{2}}   \put(230,100){\circle{5}}
\put(270,100){\circle{5}}    \put(310,100){\circle{5}}
\put(50,250){\circle{5}}  \put(90,250){\circle{5}}
\put(50,240){\line(0,1){20}} \put(90,240){\line(0,1){15}}
\put(50,270){\line(0,1){50}} \put(90,305){\line(0,1){15}}
\put(45,265){\line(1,0){35}} \put(55,295){\line(1,0){25}}
\put(45,280){\oval(20,30)[l]} \put(80,255){\oval(20,20)[rt]}
\put(80,305){\oval(20,20)[rb]}
\put(10,240){\line(0,1){80}}  \put(130,240){\line(0,1){80}}
\put(190,240){\line(0,1){80}} \put(310,240){\line(0,1){80}}
\put(230,240){\line(0,1){30}} \put(270,240){\line(0,1){30}}
\put(230,290){\line(0,1){30}} \put(270,290){\line(0,1){30}}
\bezier{30}(250,290)(256,290)(264,282)
\bezier{30}(250,270)(245,270)(236,278)
\bezier{30}(269.5,270)(269.5,275)(268.5,276)
\bezier{30}(229.5,290)(229.5,285)(231,284)
\put(250,270){\oval(40,40)[tl]} \put(250,290){\oval(40,40)[br]}
\put(270,258){\circle*{5}}  \put(270,302){\circle*{5}}
\put(20,30){\parbox[t]{300pt}{\rm \small
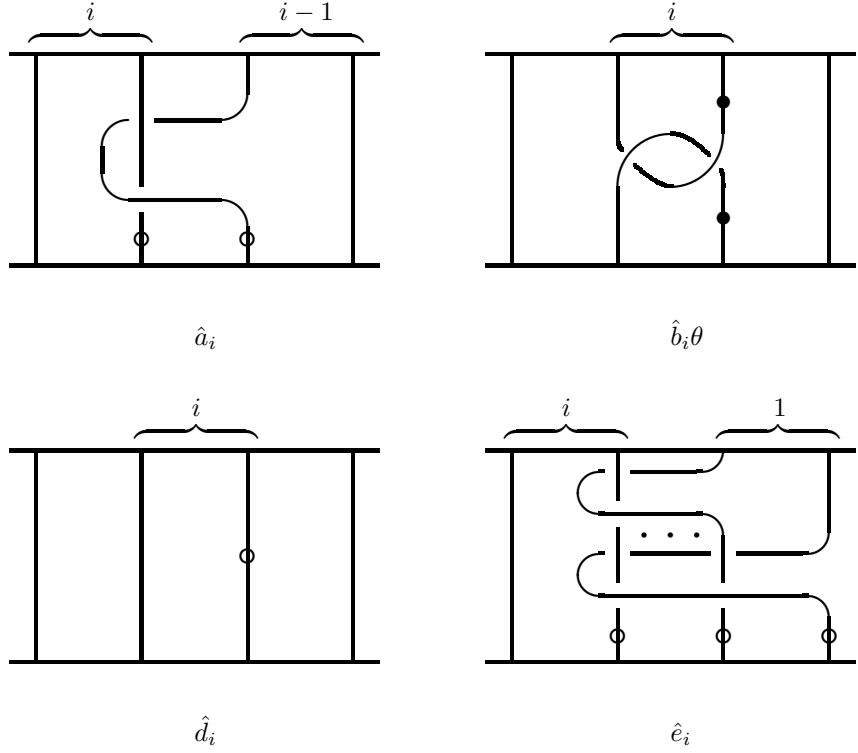
{\bf Fig. 4:} Pictorial presentation of the action of
$\ha_i, \dots, \he_i$ on states. Pairs of strands correspond
to the $g$ handles of the surface and describe a map on
$\Re$. As usual, over- and undercrossings have to be
interpreted as factors $R, R^{-1}$ and open (closed) circles
in a line stay for a factor $v^{-1}$ ($v$). Maxima and minima
mean insertion of normalized Clebsch-Gordon maps $C[\bar I I|0]$
and their conjugates. }}  \end{picture}

\noindent
{\sc Proof of the Lemma:} The simplest case is the action of
$\hat \eta_{\nu p}$
for $ p \leq m$. For notational reasons we give the proof only
in the example $\nu = 1, p = 3$. The general case can be
treated with the same ideas. With the definition (\ref{hdef})
of $\hh_{13}$ and the formulas in Theorem (\ref{Smrep}) for
the representation on monodromies $M^I_1$ and $M^I_3$ one
finds
\ba & &  D^{\{K_\rho\}} (\hat h_{1 3}) \nn \\[1mm]
    & = &
     \sum \N \th^{-1} d_I v_I tr^I_q
     \left[ ( R'_{12} R'_{13} R'_{14} R_{14} (R'_{13})^{-1}
       R_{12} )^{I K_1 K_2 K_3} \right] \nn \\[1mm]
    &=&
     \sum \N \th^{-1} d_I v_I tr^I_q
     \left[ ( R'_{12} R^{-1}_{34} R'_{14} R_{14} R_{34}
       R_{12} )^{I K_1 K_2  K_3 } \right] \nn \\[1mm]
    &=&
     \sum \N \th^{-1} d_I v_I tr^I_q
     \left[ ( R^{-1}_{34} [(id \o id \o \D) (R'_{13}
      R_{13})]_{1324} R_{34} )^{I K_1 K_2 K_3} \right] \ \ .
      \nn
\ea
For the equalities we used quasitriangularity of $R$ several
times and the subscript $_{1324}$ means that the second and third
component of the expression in brackets are exchanged. Now let
us recall that $\sum \N \th^{-1} d_I v_I tr^I_q((R'R)^I) = v^{-1}$
so that
\ba
    D^{\{K_\rho\}} (\hat h_{1 3}) & = &
  \left[ ( R^{-1}_{34} [(id \o id \o \D) (e \o e \o
  v^{-1})]_{1324}    R_{34} )^{I K_1 K_2 K_3} \right]
    \  \ \nn \\[1mm]
   & = &
   v^{-1}_{K_1} v^{-1}_{K_3}
   ( R^{-1}_{34}  R'_{24} R_{24} R_{34} )^{I K_1 K_2 K_3}
      \  \ \nn \\[1mm]
   & = &
   v^{-1}_{K_1} v^{-1}_{K_3}
   ( R'_{23}  R'_{24} R_{24} (R'_{23})^{-1})^{I K_1 K_2 K_3}
      \  \ \nn \ \ .
\ea
We inserted the formula (\ref{eigRR}) for the action of
$\D$ on the ribbon element $v$ and employed the Yang Baxter
equation. The last formula coincides with the statement
in the Lemma.

To compute the action of the $\hh_{\nu p}, p > m,$ on states
we use a map from the representation space $\Re$ of the $AB$-
algebra $\T$ to the space $\bigoplus_J V^J \o V^{\bar J}$
described by
\be \label{wJ}
w^J \equiv \k_J^{-3} C[J \bar J|0] (R')^{J \bar J}
              a^J R^{J \bar J} \in V^J \o V^{\bar J}
\ee
for every $a^J \in \End(V^J) \subset \Re$. It is easy to
check that the monodromies $B^I$ and $(A^I)^{-1}(B^I)^{-1}
A^I$ act on $w^J$ according to
\ba
  \pi(B^I) w^J & = &  w^J \k_I^{-1} (R'_{12} R'_{13} R_{13}
          (R_{12}')^{-1}) ^{I J \bar J} \ \ , \nn \\[1mm]
  \pi((A^I)^{-1} (B^I)^{-1} A^I) w^J & = &  w^J \k_I^{-1}
         (R'_{12} R_{12}) ^{I J \bar J} \ \ . \nn
\ea
These formulas should be compared with the expressions
in Theorem \ref{Smrep} for the representation of
monodromies $M_\nu^I$. The result can easily be
generalized to the case of $g$ handles. We find that
the monodromies $M^I (l_p), p > m, $ act on the summands in
$\bigoplus_{\{ J_i\}} \Im (K_1, \dots, K_m) \o
V^{J_1} \o V^{\bar J_1} \o \dots \o V^{J_g} \o V^{\bar J_g}$
according to
$$
   (\k_I)^{-1} (K_{p} R'_{1(p+1)} R_{1(p+1)} K_{p}^{-1} \o
   e^{(m-\nu)})^{I K_1 \dots K_m J_1 \bar J_1 \dots J_g
   \bar J_g}
$$
for all $p \leq m+2g$. The notations were introduced in
a remark after Theorem \ref{Smrep}. With this result, the
evaluation of $\hh_{\nu p}$ on states is reduced to the
previous case where $p \leq m$.

The evaluation of $\ha_i, \hd_i, \he_i$ is equally simple since
the curves $x(\a_i), x(\d_i), x( \e_i)$ are at most products of
elements $l_p, 1 \leq p \leq 2m+g$. Details of the straightforward
computation are left as an exercise. Only $\hb_i$ requires some
new calculation. Since $x(\b_i) = a_i$ we need to know how
$\hh (a_i)$ acts on states. Let us give the answer for the handle
algebra, i.e. $g =1$ and $\hb = \hh (a) = \hh (a_1)$. When $w^J$
is defined as in eq. (\ref{wJ}) a direct
calculation establishes
$$
  \pi(\hb) w^J   =   \th^{-1} \sum_K \N d_K v_K v_J
      w^K  C[J \bar J |0]  (R')^{J K} (R')^{\bar K \bar J}
      C[K \bar K|0]^c\ \ .
$$
The generalization of this formula to $g$ handles is
obvious. This concludes the proof of the Lemma.

\noindent
{\sc Proof of Theorem \ref{RTequiv}:} Let us recall how to obtain
the Reshetikhin Turaev representation of the mapping class group.
First one presents elements in the mapping class group by tangle
diagrams. For our generating set of Dehn twists $\eta_{\nu p},
\a_i,\b_i, \d_i, \e_i$ the corresponding diagrams can be found
e.g. in \cite{MaPo}. Then the elements in the diagram are
represented by objects associated with a Hopf algebra $\G$ in
the usual way (cp. also figure captions for Figure 3,4). This
provides us with a set of maps between representation spaces
of $\G$, one for each elements in the mapping class group.

The formulas for the action of $\hh_{\nu p}, \ha_i,\hb_i, \hd_i,
\he_i$ on states of the Chern Simons theory were described by
the pictures in Figure 3,4. These pictures agree with the
tangle diagrams which present the corresponding elements
$\eta_{\nu p}, \a_i,\b_i, \d_i, \e_i$ in the pure mapping class
group. In fact, only the picture for $\he_i$ in Fig. 4 differs
from the generators of the tangle algebra used in \cite{MaPo}.
However the two pictures are equivalent up to a Kirby move \cite{Kir}.
Because the representation constructed above is projective,
application of Kirby moves might change the result by a phase.
More precisely, whenever a manipulation with tangle diagrams involves
an $\cO_1$
move (i.e. creation or annihilation of a closed circle
\begin{picture}(9,9)(0,0) \put(5,5){\circle{8}} \end{picture})
we have to encounter a factor $\th^{\pm 1}$. The $\cO_2$ move
is respected by the representation and hence does not produce
additional factors. Any manipulation that leads from our picture
for $\he_i$ to the corresponding diagram used by Matveev and Polyak
involves one $\cO_1$ move and thus gives rise to a factor $\th$.
This concludes the proof of the theorem.

For the torus, the representation space $\Re$ can be
identified with the quantum symmetry algebra $\G$ itself.
Our representation $\pi$ of the mapping class group on
$\Re$ may then be compared with formulas of Lyubashenko
\cite{Lyb},  Majid \cite{Maj} and Kerler  \cite{Ker}
for the action of mapping class groups on Hopf algebras.

\section{Comments}
\setcounter{equation}{0}

\subsection{Truncation}

All the theory developed above was valid under the assumption
that the symmetry algebras $\G$ is semisimple. It is well
known that this requirement is not satisfied for the quantum
group algebras $U_q(\sg)$ when $q$ is a root of unity. To treat
this important case we proposed (cp. \cite{AGS1}) to use the
semisimple truncation of $U_q(\sg), q^p =1,$ which has been
constructed in \cite{MSIII}. In this truncation, semisimplicity
is gained in exchange for co-associativity, i.e. the truncated
$U_q^T(\sg)$ of \cite{MSIII} are only quasi--co-associative.
In addition, the co-product $\D$ of these truncated structures
is not unit preserving (i.e. $\D(e) \neq e \o e$). This leads
to a generalization of Drinfeld's axioms \cite{Dri2} and the
corresponding algebraic structures were called ``weak quasi-
Hopf-algebra'' in \cite{MSIII}.

Our discussions in \cite{AGS1} and especially in \cite{AGS2}
provide all background information needed in dealing
with nontrivial reassociators $\vp$, i.e. if the co-product
of $\G$ is not co-associative. Using the ``substitution
rules'' from \cite{AGS2}, the above discussion generalizes
to the quasi-Hopf case without any difficulties. Let us just recall
that all graph algebras become quasi-associative and only algebras
$\A_{g,m}$ and the moduli algebras stay associative. On the other
hand, the effect of the truncation $\D(e) \neq e \o e$ is much
more subtle and we would like to clarify this with the following
discussion. At some points throughout this paper we explicitly
used the equation $\d_I \d_J = \sum N^{IJ}_K \d_K$ for the
dimensions $\d_I$ of the irreducible representations $\t^I$.
This fails to hold in the presence of truncation and we have
to encounter the inequality $\d_I \d_J \geq \sum N^{IJ}_K
\d_K $ instead.

A first effect of such truncations is that they reduce the
number of linear independent matrix elements in the quantum
monodromies $M^I_\nu, A^I_i, etc. $ beyond the value of
$\d_I^2$. In fact one may derive the following relations
$$   \t^I(\S(e^1_\s)) M^I(\C) \t^I(e^2_\s) = M^I(\C)
  \mbox{ with } \ \ \D(e) = \sum e^1_\s \o e^2_\s \ \ . $$
Here $M^I(\C)$ is defined as in eq. (\ref{MonC}) and may in
particular be equal to $A^I_i, B^I_i$ or $M^I_\nu$.
Let us shortly sketch the proof under the simplifying assumption
that $\D$ is co-associative. Since the element $e\in \G$
is supposed to be the identity in the graph algebras, we
know that $M^I(\C) = e M^I(\C)$. From the covariance properties
of monodromies one concludes
$$     M^I (\C) = \t^I(\S(e^{11}_{\s\t})) M^I (\C)
       \t^I(e^{12}_{\s\t}) e^2_\s\ \ . $$
This formula is used twice in the following calculation
\ba
  & &   \t^I(\S(e^1_\rho)) M^I(\C) \t^I(e^2_\rho)  \nn \\[1mm]
  &  =  &
    \t^I(\S(e^1_\rho)) \t^I(\S(e^{11}_{\s\t})) M^I (\C)
       \t^I(e^{12}_{\s\t}) \t^I(e^2_\rho) e^2_\s  \nn \\[1mm]
  & = &
    \t^I(\S(e^{11}_{\s\t})) M^I (\C)
     \t^I(e^{12}_{\s\t}) e^2_\s = M^(\C) \ \ . \nn
\ea
The second step simply follows from $ \D(e) \D(\xi) = \D (e \xi)
= \D(\xi) $. We conclude that $M^I(\C)$ contains only
$N^{I\bar I}_K \d_K \leq \d_I^2$ linear independent components.
The total dimension of the algebra $\L (M(\C))$ generated by
matrix components of $M^I(\C)$ is then $\sum_{I,K} N^{I
\bar I}_K \d_K$.

With respect to Lemma \ref{irredlemma}, this remark shows that
we cannot expect all linear maps in $\End(V^J)$ to appear
in the image of quantum monodromies under the representations $D^J$.
However, a careful reexamination of the proof shows that the
faithfulness of the representation theory is still guaranteed.
Indeed the same proof allows us to find $\sum N^{I \bar I}_K \d_K$
linear independent representation matrices on a space $V^I$. The
sum of these numbers over $I$ coincides with the dimension of the
the (quasi-associative) loop algebra $\L$ (cp. preceding paragraph).

All the proofs we gave in the Section on the genus 0 case were
designed so that they only use the faithfulness. They now carry
over to the truncated situation. Let us emphasize that the
irreducibility of representations is not completely lost. We
saw above counting arguments prevent the representations $D^I$
of the loop algebra $\L$ from being irreducible.
The same holds true for the multi-loop case. However, when we
descent to moduli algebras,
counting arguments allow to derive irreducibility of the
representation theory from the faithfulness much as this was
done in Section 6.3. So the general rule is that irreducibility
statements fail on the level of graph algebras but hold again
on the level of moduli algebras.

There is a similar story about handle algebras. We continue
to generate the representation space $\Re$ from a ground
state $\vac$ by application of the ``creation operators''
$ A^I $. From our discussion above, the resulting space
has the same dimension as the loop algebra $\L(A)$ associated
with monodromies $A^I$, i.e. we find ${\it dim\/}
(\Re) = \sum_{I,K} N^{I \bar I}_K \d_K $. In particular, the
space $\Re$ is no longer isomorphic to $\bigoplus_I V^I
\o V^{\bar I}$. We may, however, write $\Re \cong \bigoplus
(V^I \o V^{\bar I})'$ where $(V^I \o V^{\bar I})'$ is the
subspace of $V^I \o V^{\bar I}$ on which $D(e)^{I \bar I}$
acts as identity. Our proof at the end of Subsection 7.3
shows that the representations $\pi$ of the
handle algebra $\T$ is still
irreducible. In fact, the projectors $\chi_B^I$ now project
to the $\sum_K N^{I \bar I}_K \d_K$-dimensional subspace
$(V^I \o V^{\bar I})'$ of $\Re$. Monodromies $B^J$ and
$(A^J)^{-1} (B^J)^{-1} A^J$ give rise to $ (\sum_K
N^{I \bar I}_K \d_K)^2$ linear independent maps on
$(V^I \o V^{\bar I})'$. So we obtain a full matrix
algebra over this subspace. Irreducibility of $\pi$ is then
obtained as before.
Faithfulness of the representation follows from irreducibility
by the usual counting of dimensions. From this point on, the
rest of the discussion -- in particular the representation
theory at higher genera -- is straightforward.

These remarks suffice to generalize the above theory to
weak quasi Hopf algebras and hence to open the way for
its application to interesting symmetry algebras
associated with $U_q(\sg), q^p=1$.

\subsection{Open problems}

In this paper we essentially completed the program of operator
quantization of the Chern-Simons theory. We have defined the algebra
of observables equipped with the $*$-operation and the positive
integration functional and developed its representation theory.
It is shown that the list of irreducible unitary representations
matches with  the Hilbert spaces produced by Geometric Quantization
of the moduli space of flat connections.

Let us mention here several open problems related to the
framework of Combinatorial Quantization.

We mentioned in Section 4 that for generic values of $q$
the moduli algebra
$\MC_{g,m}^{\{ I_\nu\}}$ is isomorphic (as a linear space)
to the classical algebra of analytic functions on the
moduli space of flat connections. Eventually, these spaces
may be identified in many ways. However, we can restrict
the choice by the following condition. Pick up
some circle $x$ on the surface. Corresponding to this
circle one can construct (see Section 9) a set of
elements $c^I(x)\in \MC_{g,m}^{\{ I_\nu\}}$ via
\be \label{cIx}
c^I(x)=\k_I \tr_q M^I(x).
\ee
The representations of the quantum symmetry algebra
for generic $q$ being in one to one correspondence with
the representation of the group $G$, one can construct
classical counterparts of elements (\ref{cIx}):
\be
c^I_0(x)= \tr M^I(x).
\ee
These are already commuting analytic functions on the moduli
space.

So, one wishes to construct the map ${\cal Q}$ which maps
classical functions into the elements of the moduli algebra so that
\be        \label{calQ}
{\cal Q} (c^I_0(x)) = c^I(x)
\ee
for every $x$ and $I$. It is not {\em a priori} clear that such
a map exists. At least we wish to preserve  (\ref{calQ}) for
sufficiently many cycles on the surface.

Another conceptual problem related to the quantized moduli
space is the issue of a differential calculus. The classical
moduli space has a rich cohomology theory \cite{Wgt}, \cite{JK}.
{}From this perspective it would be interesting to define
the corresponding cohomology problem for quantized moduli spaces.
One can hope that the differential calculi on the quantum groups
may provide a good starting point.\\[1.5cm]
{\bf Acknowledgements:} We would like to thank S. Axelrod,
G. Felder, E. Frenkel, K. Gawedzki, T. Kerler, N. Reshetikhin,
and J. Roberts for comments and discussions. The work of V.S.
was partly supported by the Department of Energy  under DOE
Grant No. DE-FG02-88ER25065.


\begin{thebibliography}{99}
\bibitem{AFS} A. Yu. Alekseev, L.D.Faddeev, M.A.Semenov-Tian-Shansky
 {\em Hidden Quantum groups inside Kac-Moody algebras},
 Commun. Math. Phys. {\bf149}, no.2 (1992) p.335
\bibitem{AGS1} A. Yu. Alekseev, H. Grosse, V. Schomerus, {\em
 Combinatorial quantization of the Hamiltonian Chern Simons theory I},
 HUTMP 94-B336, Commun. Math. Phys., to appear
\bibitem{AGS2} A. Yu. Alekseev, H. Grosse, V. Schomerus, {\em
 Combinatorial quantization of the Hamiltonian Chern Simons theory II},
 HUTMP 94-B336, Commun. Math. Phys., to appear
\bibitem{ASh2} A. Yu. Alekseev, V. Schomerus, {\em Mapping class
 groups and moduli algebras}, in preparation
\bibitem{AA} A. Yu. Alekseev,
 {\em Integrability in the Hamiltonian Chern-Simons theory},
 hep-th/9311074, St.-Petersburg  Math. J. {\bf vol. 6} 2 (1994) 1
\bibitem{AlCo} D. Altschuler, A. Coste {\em Quasi-quantum groups,
 knots, three manifolds and topological field theory},
 Commun. Math. Phys. {\bf 150} (1992) 83
\bibitem{AMR} J. Andersen, J. Mattes, N. Reshetikhin {\em
 Poisson structure on the moduli space of flat connections
 and chord diagrams}, in preparation;  \\
 {\em Quantization of
 the algebra of chord diagrams}, in preparation
\bibitem{AtBo} M. Atiyah, R. Bott, {\em The Yang-Mills equation
 over Riemann surfaces}, Phil. Trans. of the Royal Soc. of
 London, Ser A {\bf 308} (1982) 523
\bibitem{ADW} S. Axelrod, S. Della Pietra, E. Witten, {\em
 Geometric quantization of Chern-Simons gauge theory},
 J. Diff. Geom. {\bf 33} (1991) 787
\bibitem{AS}  S.Axelrod, I.M.Singer,
 {\em Chern-Simons Perturbation Theory II},
 preprint HEP-TH/9304087.
\bibitem{Bou} D.V. Boulatov, {\em q-deformed lattice gauge theory
 and three manifold invariants}, Int. J. Mod. Phys. {\bf A8},(1993),
 3139
\bibitem{BuRo} E. Buffenoir, Ph. Roche, {\em Two dimensional
 lattice gauge theory based on a quantum group}, preprint
 CPTH A 302-05/94 and hep-th/9405126
\bibitem{KC} DeConcini, V. Kac {\em Representations of quantum
 groups at roots of 1}, in {\it Progress in Mathematics} vol.
 {\bf 92}, Birkh\"auser 1990.
\bibitem{Dri1} V.G. Drinfel'd, {\em Quantum groups}, Proc. ICM
 (1987) 798
\bibitem{Dri2} V.G. Drinfel'd, {\em Quasi Hopf algebras and Knizhnik
 Zamolodchikov equations}, in: Problems of modern quantum field
 theory, Proceedings Alushta 1989, Research reports in physics,
 Springer Verlag Heidelberg 1989 \\
 V.G. Drinfel'd, {\em Quasi-Hopf algebras}, Leningrad. Math. J.
 Vol. {\bf 1} (1990), No. 6
\bibitem{DJN} B. Durhuus, H.P. Jakobsen, R. Nest, {\em
 A construction of topological
 quantum field theories from 6j-symbols}, Nucl. Phys. B (Proc. Suppl.)
 (1991), 666
\bibitem{EM} S. Elitzur,G. Moore, A. Schwimmer, N. Seiberg, {\em
 Remarks on the canonical quantization of the Chern-Simons-Witten
 theory}, Nucl. Phys. {\bf 326}, (1989), 108
\bibitem{FRT} L.D. Faddeev, N. Yu. Reshetikhin, L.A.
 Takhtadzhyan, {\em Quantization of Lie Groups and Lie Algebras},
 Algebra and Analysis {\bf 1} (1989), 1 and Leningrad Math. J.
 Vol. 1 (1990), No. 1
\bibitem{Fed} B.V. Fedosov, {\em Formal quantization}, Some
 topics of modern Mathematics and their application to problems
 in Mathematical Physics, Moscow (1985), 129; \\
 {\em A simple  geometrical construction
 of deformation quantization}, J. Diff. Geom., to appear
\bibitem{Fin} M.Finkelberg, {\em Fusion Categories},
 Harvard University thesis, May 1993
\bibitem{FoRo} V.V. Fock, A.A. Rosly, {\em Poisson structures on
 moduli of flat connections on Riemann surfaces and r-matrices},
 preprint ITEP 72-92, June 1992, Moscow
\bibitem{FrGa} J. Fr\"ohlich, F. Gabbiani, {\em Braid statistics
 in local quantum theory}, Rev. Math. Phys. {\bf 2} (1991) 251
\bibitem{Gol} W. Goldman {\em Invariant functions on Lie groups
 and Hamiltonian flows of surface group representations}, Invent.
 Math. {\bf 85} (1986) 263; \\ {\em  The symplectic nature
 of fundamental groups of surfaces}, Adv. Math. {\bf 54} (1984) 200
\bibitem{JGH} F.M. Goodman, P. de la Harpe, V.F.R. Jones {\em
 Coxeter graphs and towers of algebras}, Springer Verlag, 1989
\bibitem{JK} L.Jeffrey, F.Kirwan, talk given at Newton Institute,
Cambridge, November 1994.
\bibitem{KaSc} M. Karowski, R. Schrader, {\em  A combinatorial approach
 to topological quantum field theories and invariants of graphs},
 Commun. Math. Phys. {\bf 151}, (1993), 355 \\
 M. Karowski, R. Schrader, {\em A quantum group version
 of quantum gauge theories in two dimensions}, J. Phys. {\bf A 25}, (1992),
 L1151
\bibitem{KaLu} D. Kazhdan, G. Lusztig, {\em Tensor structures
 arising from affine Lie algebras I/II}, J. Am. Math. Soc.
 {\bf vol. 6} 4 (1993) 905
\bibitem{Ker} T. Kerler, {\em Mapping class group actions on
 quantum doubles}, Commun. Math. Phys., to appear
\bibitem{Kir} R. Kirby, {\em A calculus for framed links in
 $S^3$}, Invent. Math. {\bf 45} (1978) 35
\bibitem{KZ} V. Knizhnik, A. B. Zamolodchikov, {\em Current algebra
 and Wess-Zumino model in two dimensions}, Nucl. Phys. {\bf B247}
 (1984) 83
\bibitem{Lus} Lusztig, {\em Finite dimensional Hopf algebras
 arising from quantized universal enveloping algebras}
 J. Am. Math. Soc. {\bf vol. 3} 1 (1990) 257
\bibitem{Lyb} V. Lyubashenko, {\em Tangles and Hopf algebras
 in braided tensor categories}, J. Pure Appl. Alg., to appear
\bibitem{MSIII} G. Mack, V. Schomerus, {\em Quasi Hopf quantum
 symmetry in quantum theory}, Nucl. Phys. {\bf B370}(1992),185
\bibitem{Maj} S. Majid, {\em Braided groups}, J. Pure Appl.
 Alg. {\bf 86} (1993) 187
\bibitem{MaPo} S. Matveev, M. Polyak, {\em A geometrical presentation
 of the surface mapping class group and surgery}, Commun. Math.
 Phys. {\bf 160} (1994)
\bibitem{ReSTS} N. Reshetikhin, M. Semenov-Tian-Shansky, {\em
  }, Lett. Math. Phys. {\bf 19} (1990) 133
\bibitem{ReTu} N. Reshetikhin, V. G. Turaev, {\em Ribbon
 graphs and their invariants derived from quantum groups}, Commun.
 Math. Phys. {\bf 127}, (1990),1
\bibitem{ReTu2} N. Reshetikhin, V. G. Turaev, {\em Invariants of
 3-manifolds via link polynomials and quantum groups}, Invent.
 Math. {\bf 103} (1991) 547
\bibitem{Sud} A. Sudbury, {\em Non-commuting coordinates and
 differential operators} in: Quantum groups, T. Curtright et al. (eds),
 World Scientific, Singapore 1991
\bibitem{Tur} V. Turaev, {\em Skein quantization of Poisson
 algebras of loops on surfaces}, Ann. scient. Ec. Norm. Sup.
 4e serie, t. 24 (1991) 635
\bibitem{TuVi} V.G. Turaev, O.Y. Viro, {\em  State sum of 3-manifolds
 and quantum 6j-symbols}, LOMI preprint (1990)
\bibitem{Ver} E. Verlinde, {\em Fusion rules and modular
 transformations in 2D conformal field theory}, Nucl. Phys.
 {\bf B 300} (1988) 360
\bibitem{Wit1} E.Witten,
 {\em Quantum field theory and the Jones  polynomial}, Commun. Math. Phys.
 {\bf 121}, (1989),351
\bibitem{Wgt} E.Witten,
{\em Two-dimensional gauge theory revisited}, J. Geom. Phys. {\bf 9}
(1992) 303
\end{thebibliography}
\end{document}